\documentclass[12pt]{article}
\usepackage{amsmath,amsfonts,amssymb,amsthm,amscd}
\title{On NIP and invariant measures}
\date{January 28, 2009}
\author{Ehud Hrushovski\thanks{Thanks for support from the ISF}\\Hebrew University of Jerusalem\and Anand Pillay\thanks{Supported
by a Marie Curie Chair EXC 024052 and EPSRC grant EP/F009712/1}\\University of Leeds}
\newtheorem{Theorem}{Theorem}[section]
\newtheorem{Proposition}[Theorem]{Proposition}
\newtheorem{Definition}[Theorem]{Definition} 
\newtheorem{Remark}[Theorem]{Remark}
\newtheorem{Lemma}[Theorem]{Lemma}
\newtheorem{Corollary}[Theorem]{Corollary}
\newtheorem{Fact}[Theorem]{Fact}

\newtheorem{Question}[Theorem]{Question}
\newtheorem{Problem}[Theorem]{Problem}
\newcommand{\R}{\mathbb R}   
\newcommand{\Q}{\mathbb Q}  
\newcommand{\Z}{\mathbb Z}  
\newcommand{\N}{\mathbb N}

\newcommand{\pf}{\noindent{\em Proof. }}

\begin{document}
\maketitle

\begin{abstract} We study forking, Lascar
strong types, Keisler measures and definable groups, under an assumption of
$NIP$ (not the independence property), continuing aspects of the 
paper \cite{NIP}.  Among key results are
(i) if $p = tp(b/A)$ does not fork over $A$ then the Lascar strong type of $b$ over $A$ coincides with the 
compact strong type of $b$ over $A$ and any global nonforking extension of $p$
is Borel definable over $bdd(A)$, (ii) analogous statements for Keisler 
measures and definable groups, including the fact that $G^{000} = G^{00}$ for $G$ definably amenable,
(iii) definitions, characterizations 
and properties of 
``generically stable" types and 
groups, (iv) uniqueness of invariant (under the group action) Keisler measures on groups with finitely satisfiable generics,
(v) a proof of the compact domination conjecture for (definably compact) commutative groups  in 
$o$-minimal expansions of real closed fields.  
\end{abstract}

\section{Introduction and preliminaries}
The general theme of this paper is to find and study stable-like behaviour in theories and
definable groups without the
independence property. This was
a theme in the model-theoretic analysis of algebraically closed valued fields \cite{HHM1}, \cite{HHM2}, 
\cite{metastable}. It was also an aspect of the paper \cite{NIP}, although there the main 
motivation
was to resolve some conjectures about definably compact groups in $o$-minimal structures. In fact a bonus in 
the current paper is a proof of the compact domination conjecture (formulated in \cite{NIP}), at least for commutative
groups, and various 
elaborations, which is fairly direct but also uses some of the general machinery we develop.

Both authors have been a bit slow to realize the
relevance of Shelah's notion of forking to theories with $NIP$. We rectify this in the current paper where we will be 
quite explicit about the role of forking.

Note that a theory $T$ is stable if and only if it is simple and has $NIP$. In stable
theories we have on the one hand the ``algebraic" behaviour of nonforking independence,
namely symmetry, transitivity, existence of nonforking extensions, as well
as local character (any type does not fork over some small set). On the other hand we have 
(again in stable theories) what can be called
broadly ``multiplicity theory", the structure and behaviour of nonforking extensions of a type.
Included in ``multiplicity theory" are alternative characterizations of nonforking, for example a global type $p$ does not fork 
over a model $M$ iff $p$ is definable over $M$ iff $p$ is finitely satisfiable in $M$
iff $p$ is $Aut({\bar M}/M)$-invariant (does not split over $M$) Included also is the finite equivalence
relation theorem: a type over an algebraically closed set is stationary, namely has a 
unique global nonforking extension. In the early 
texts on stability theory (the original papers and book of Shelah, but also the treatment by Lascar and Poizat), the 
proofs and development of the
algebraic properties of forking were tied up  with  multiplicity theory (in the form of heirs, coheirs for example). In the work on simple theories, the 
two strands were distinguished, the algebraic theory being valid in all simple theories, but not the latter. Multiplicity theory 
did make an appearance in simple theories, but in the
(weak) form of the ``Independence theorem over a model" and more generally over boundedly closed (rather than
algebraically closed) sets. As the validity of the algebraic theory of forking is characteristic
of simple theories, it will fail for unstable theories with $NIP$. One of the points of the current paper
is to recover aspects of the multiplicity part of stability theory for theories with $NIP$.
Invariance (rather than stationarity)  turns out to be important and already invariant types played a role in 
the analysis of algebraically closed valued fields. We show for example (in section 2), extending work of Shelah \cite
{Shelah715} and Adler \cite{Adler}  that assuming $NIP$ if $p$ is a 
global type and $A$ a small set then $p$ does not fork over $A$ iff $p$ is $Aut({\bar M}/bdd(A))$-invariant 
iff $p$ is ``Borel definable" over $bdd(A)$.  
Keisler measures figured a lot in the earlier paper \cite{NIP} and we entertained the possibility that 
replacing types by Keisler measures would give a smoother theory and better chance of recovering
stationary-like behaviour (uniqueness of nonforking extensions). In section 4 the results
for types are carried over to Keisler measures. This includes the ``Borel definability" over $bdd(A)$ of a measure which does not fork over $A$, the proof of which uses the Vapnik-Chervonenkis theorem. We reduce measures to types in two ways. 
(a) If the Keisler measure is $A$-invariant then it corresponds to a Borel measure on the
space of $bdd(A)$-invariant types (see 4.6), (b) An invariant Keisler measure is also the
"frequency average" of some sequence of invariant types (see 4.8). Though these representations
of Keisler measures might seem to make considerations of measures unnecessary, in fact some of our proofs of facts about types go through
measures. 
This point appears in section 5, where among other things we show that definable amenability (the existence of 
left invariant measures) of a definable group $G$
is equivalent to
the existence of generic types in the sense of forking.

In section 3 we study what we call {\em generically stable types} (and what Shelah has just called stable types), mainly in an NIP context. These are types whose behaviour vis-a-vis nonforking extensions is like that in stable theories, namely a global nonforking extension is both definable and finitely satisfiable. A special case of a generically stable type is a stably dominated type, as introduced in \cite{HHM2}, and
our results place those of \cite{HHM2} in the appropriate general context. Likewise {\em generically stable groups}, studied in section 6, generalize the stably dominated groups of \cite{metastable}.

In section 5 we recall the groups with $fsg$ (finitely satisfiable generics) which were introduced in \cite{NIP}, and
generalize some results from \cite{NIP} such as definable amenability, to definable groups with generic types in the sense of forking. 

Definably compact 
groups in a variety of settings ($o$-minimal, $P$-minimal, metastable) have either been proved to have
or are expected to have $fsg$. In section 7 we prove the {\em uniqueness} of translation invariant Keisler measures in $fsg$ groups. We see this as a kind of common generalization of the uniqueness of global translation invariant {\em types} for stable groups, and of the uniqueness of Haar measure for compact groups, and exemplifies our search for stable-like behaviour at the level of Keisler measures.

It is natural to try to both extend the notion of generic stability from types to measures, as well as to find group-free versions of the $fsg$ property. This is discussed briefly in sections 4 and 7, and the problems will be addressed in a future paper with P. Simon. 

Section 8 is devoted to a proof of compact domination for commutative definably compact groups in $o$-minimal expansions of 
real closed fields. In fact we prove a strong version, domination of $G$ by a ``semi-$o$-minimal" compact 
Lie group $G/G^{00}$.  Direct $o$-minimal style arguments as well as some of the general theory of 
invariant types, play a role in the proof. Our proof makes use of a theorem on the existence of definable Skolem functions for ``$o$-minimal subsets of finite-dimensional groups" which is proved in the appendix.

\vspace{5mm}
\noindent
Much inspiration for our work on measures comes from Keisler's work \cite{Keisler1} and \cite{Keisler2}. But our emphasis differs from Keisler's. Among the main points of Keisler is that if $\phi(x,y)$ is a stable formula, then any $\phi$-measure is a ``convergent weighted sum" of $\phi$-types. (Here a $\phi$-measure over $M$ is a finitely additive measure on the Boolean algebra generated by formulas of the form $\phi(x,a)$ for $a\in M$.) So all the properties of nonforking in a stable theory (symmetry, stationarity of types over models,...) pass automatically to measures.  Keisler does consider measures in (possibly unstable) theories with $NIP$. The good class of measures he identifies are so-called ``smooth" or ``minimally
unstable" measures.  Loosely (and possibly incorrectly) speaking, a  smooth measure over $M$ is a measure over $M$ which has a unique extension to the ``unstable part " of $\bar M$. He proves that such measures exist. So the only freedom in taking extensions of a smooth measure is with respect to stable formulas, and thus in a sense forking theory for smooth measures essentially reduces to forking theory for measures in stable theories. On the other hand our main focus will be global measures which are invariant over some
small set. Only in special cases will these also be smooth.

\vspace{5mm}
\noindent 
We will use  standard notation. Namely $T$ denotes a complete theory in a language $L$ and we work in 
a saturated model $\bar M$ of $T$. $x,y,z..$ usually denote finite tuples of variables, $A,B,C$ small subsets 
of $\bar M$ and $M,N,..$ small elementary submodels of $\bar M$.  There is no harm in working in ${\bar 
M}^{eq}$, but sometimes we will
assume our theory is one-sorted. A reader would benefit from having some familiarity with stability theory
and stable group theory, a reference for which is the first chapter of \cite{Pillay-book}. However the paper 
is reasonably self-contained and can on the whole be read  
independently of \cite{NIP}.

When it is convenient we denote the space of complete types over $A$ in free variable $x$ by $S_{x}(A)$. By a 
global (complete) type we mean some $p(x)\in S(\bar M)$. Recall that a partial type $\Sigma(x)$ (over 
some set of parameters and closed under finite conjunctions) is said to divide over a set $A$, if there is 
$\phi(x,b)\in \Sigma(x)$ and an $A$-indiscernible sequence $(b_{i}:i<\omega)$ of realizations of $tp(b/A)$ 
such that $\{\phi(x,b_{i}):i<\omega\}$ is inconsistent. $\Sigma(x)$ is said to fork over $A$ if $\Sigma$ 
implies some finite disjunction of formulas, each of which divides over $A$. Note that if $p(x)$ is a 
complete type over some $|A|^{+}$-saturated model $M$ and $A\subseteq M$ then $p(x)$ divides over $A$ iff 
$p(x)$ forks over $A$. Also if $A$ is any (small) set of parameters, and $p(x)\in S(A)$, then $p$ does not 
fork over $A$ if and only if $p(x)$ has a global extension which does not divide (equivalently does not fork) 
over $A$. 

A rather more basic notion is splitting: Let $p(x)\in S(B)$ and $A\subseteq B$. We say that $p$ does not split over $A$ if for any $L$-formula $\phi(x,y)$ and $b,c\in B$, if $tp(b/A) = tp(c/A)$ then $\phi(x,b)\in p$ iff $\phi(x,c)\in p$. This tends to be more meaningful when $B$ is some $|A|^{+}$-saturated model. 

We will be discussing various kinds of strong types, so let us fix notation. First by $Aut({\bar M}/A)$ we mean the group of automorphisms of $\bar M$ which fix $A$ pointwise. Two tuples (of 
the same length or elements of the same sort) $a$ and $b$ are said to have the same strong type over a set $A$ of 
parameters, if $E(a,b)$ for each finite $A$-definable equivalence relation $E$ on the relevant sort (where 
finite means having only finitely many classes). It is well-known that $stp(a/A) = stp(b/A)$ if and only if 
$a$ and $b$ 
have the same type over $acl(A)$ where $acl(-)$ is computed in ${\bar M}^{eq}$. To be consistent with later 
notation it might be better to call strong types, Shelah strong types or profinite strong types.

We say that $a$ and $b$ have the same compact strong type or $KP$ strong type over $A$ if $E(a,b)$ for each bounded 
equivalence relation $E$ on the relevant sort which is type-definable over $A$, that is defined by a possible 
infinite set of formulas over $A$. Here bounded means having strictly less than $|{\bar M}|$-many classes 
which is 
equivalent to having at most $2^{|T|+|A|}$-many classes. 
An ($\emptyset$-) hyperimaginary $e$ is by definition the class 
of 
a $\emptyset$-type-definable equivalence relation. We define $bdd(A)$ to be the set of hyperimaginaries which have 
small 
orbit under $Aut({\bar M}/A)$. It follows from \cite{Lascar-Pillay} (see 4.18 there) that $a$ and $b$ have 
the same $KP$ strong type over $A$ iff they have the same type over $bdd(A)$ (where types over 
hyperimaginaries are made sense of in \cite{HKP} for example). In particular if $a$ and $b$ have the same $KP$ 
strong type over $A$, then for any $A$-type definable set $X$ to which $a,b$ belong and bounded $A$
-type-definable equivalence relation $E$ on $X$, $a$ and $b$ are in the same $E$-class. All this makes sense 
for type-definable equivalence relations on sets of infinite tuples, although in \cite{Lascar-Pillay} we 
pointed out that a bounded infinitary hyperimaginary is ``equivalent" to a sequence of bounded finitary 
hyperimaginaries.

Finally we say that tuples $a$ and $b$ have the same Lascar strong type over $A$, if $E(a,b)$ for any 
bounded equivalence relation $E$ which is invariant under $Aut({\bar M}/A)$. The relation ``equality of 
Lascar strong type over $A$" is the transitive closure of the relation $R_{A}(-,-)$, where $R_{A}(x,y)$ says: $x$ and $y$ are elements of some infinite $A$-indiscernible sequence..
Following Lascar we let $Autf({\bar M}/A)$ denote the group of automorphisms of $\bar M$ which fix all Lascar strong types over $A$. Tuples $a,b$ from $\bar M$ (which are possibly infinite) will have the same Lascar strong type over $A$ if and only if they are in the same orbit under $Autf({\bar M}/A)$.

Note that Lascar strong type refines KP strong type refines strong type refines type (all over $A$). If $A = 
M$ is a model, Lascar strong type coincides with type. In a 
stable theory Lascar strong types over $A$ coincide with strong types over $A$. This is conjectured to be 
true in simple theories too, and was proved in the supersimple case. In the general simple case we only know 
that Lascar strong types coincide with $KP$ strong types. In any case in the current paper we will consider 
such 
questions for theories with $NIP$. 

\vspace{2mm}
\noindent
We will also be refereeing to various kinds of ``connected components" of a definable (or even type-definable 
group) $G$. Suppose $A$ includes the parameters over which $G$ is (type)-defined. Let $G^{0}_{A}$ be the
intersection of all (relatively) $A$-definable subgroups of $G$ of finite index. Let $G^{00}_{A}$ denote
the smallest type-definable over $A$ subgroup of $G$ of ``bounded" index (equivalently index at most 
$2^{|L|+|A|}$).  Let $G^{000}_{A}$ denote the 
smallest subgroup of $G$ of bounded index which is $Aut({\bar M}/A)$-invariant. If for example $G^{0}_{A}$ does
not depend on $A$, but only on $G$ we say that $G^{0}$ exists. Likewise for $G^{00}$ and $G^{000}$. 

There is an analogy between definable groups $G$ (and their quotients such as $G/G^{00}$ etc.) and automorphism 
groups (and their corresponding quotients) which permeates this paper. A basic construction from the first
author's thesis produces an automorphism group from a definable group (action). We recall this now, calling it
\newline 
{\bf Construction $C$}.

We start with a theory $T$ which includes a transitive $\emptyset$-definable group action of $G$ on $X$. Assume 
that $X$ is stably embedded and that there is some finite subset $\{a_{1},..,a_{m}\}$ of $X$ such that the 
pointwise stabilizer in $G$ of $\{a_{1},..,a_{m}\}$ is trivial
Fix a (saturated) model $M$ of $T$, and add a new sort $X'$ and a bijection $h:X\to X'$. Call the
resulting structure $M'$ and its theory  $T'$. Let $M''$ be the reduct of $M'$ which has all
existing relations on $M$ as well as all relations on $X'$ and $X\times X'$  which
are $\emptyset$-definable in $M'$ and $G$-invariant (under the action of $G$ induced by $h$
and the action of $G$ on $X$). Let $T_{X} = Th(M'')$. 
Then we have:
\newline
{\bf Conclusion.} In $M''$ the action of $G$ on $X'$ is by automorphisms, namely for every $\sigma\in G$,
$\sigma|X'$ extends to an automorphism of $M''$. In fact the restriction of
$Aut(M'')$ to $X'$ equals the semidirect product of $G$ and the restriction of $Aut(M)$ to $X$.

\vspace{2mm}
\noindent
Note that after naming $a_{1},..,a_{m},b_{1},..,b_{m}$ we recover $h$, and so $M'$ is ``internal" to $M''$:
Let $h^{m+1}:X^{m+1}\to X'^{m+1}$ be the bijection induced by $h$, and let $H$ be the union of the $G$-conjugates
of $h^{m+1}$. So $H$ is $\emptyset$-definable in $M''$, and note that $h(x) = y$ iff
$(x,a_{1},..,a_{m},y,b_{1},..,b_{m})\in H$.

\vspace{2mm}
\noindent
In the case where the action of $G$ on $X$ is regular we have the following equivalent construction.
Simply add a new sort $X'$ with a regular action of $G$ on $X'$. In this structure $G^{opp}$ acts on $X'$ as 
automorphisms and
the full automorphism group of the new structure is the semidirect product of $Aut(M)$ and $G$.

\vspace{2mm}
\noindent
We will repeatedly use the following characteristic property of theories with $NIP$.
\begin{Fact} Suppose $T$ has $NIP$. Then for any formula $\phi(x,y)$, there is $N<\omega$, such that if $(a_{i}:i<\omega)$ is an indiscernible sequence, then there does not exist $b$ such that 
$\neg(\phi(a_{i},b)\leftrightarrow\phi(a_{i+1},b))$ for $i=0,..,N-1$. 

In particular if $(a_{i}:i<\omega)$ is totally indiscernible (or an indiscernible set), then for any $b$, either
$|\{i<\omega:\models\phi(a_{i},b)\}| \leq N$ or 
\newline
$|\{i<\omega:\models\neg\phi(a_{i},b)\}| \leq N$.
\end{Fact}

At some point we will, assuming $NIP$, refer to $Av(I/M)$ where $I$ is some infinite indiscernible sequence (with no last element).
It is the complete type over $M$ consisting of formulas with parameters from $M$ which are true on a cofinal 
subset of $I$. This makes sense by Fact 1.1.

\vspace{5mm}
\noindent
The second author would like to thank H. Adler, A. Berarducci, C. Ealy and K. Krupinski for  
helpful conversations and 
communications around the topics of this paper. He would also like to thank the Humboldt Foundation for their support of a visit to Berlin in 
March-April 2007 when some of the work on this paper was done.

A first version of this paper was written in October 2007. Both authors would like to thank Itay Kaplan, Margarita Otero, Kobi Peterzil, Henryk Petrykowski, and Roman Wencel, for their detailed reading of parts of the manuscript and for pointing out errors, gaps, and possible improvements. In the new version several proofs are expanded and clarified, particularly in the current sections 2, 4, 7 and 8. Substantial changes include removing the old section 8 on generic compact domination, giving a more complete account of the Vapnik-Chervonenkis theorem and its applications in section 4, as well as adding an appendix proving the existence of definable Skolem functions in suitable $o$-minimal structures (which is needed for the proof of compact domination in the current section 8).

\section{Forking and Lascar strong types}
Forking in $NIP$ theories typically has a different character from forking in simple theories (although the 
definition, as in the introduction, is the same). 
In simple theories, forking is associated to a ``lowering of dimension". In $NIP$ theories forking can come
from just a lowering of order of magnitude within a given dimension. Although dimension is no
less important in $C$-minimal and $o$-minimal theories than in strongly minimal ones, we do not at the moment 
know 
the right $NIP$ based notion that specializes to
lowering of dimension in these cases. (Thorn forking is of course a very useful notion but does not apply
to the $C$-minimal case.)

This section builds on work of Poizat \cite{Poizat}, Shelah \cite{Shelah783} and Adler \cite{Adler}. Many of our key notions make an explicit or implicit appearance in Chapter 12 of the Poizat reference.
For completeness we will begin by restating some of the results by the above mentioned people. 
The first is a striking characterization of forking in $NIP$ theories from
\cite{Adler} but with roots in \cite{Shelah783}. 
\begin{Proposition} (Assume $NIP$) Let $p(x)\in S(\bar M)$ be a global type and $A$ a (small) set.
Then
\newline
(i) $p$ does not fork over $A$ iff $p$ is $Autf({\bar M}/A)$-invariant, in other words if
$p(x)$ is fixed by any automorphism of ${\bar M}$ which fixes all Lascar strong types over $A$.
\newline
(ii) In particular if $A = M$ is a model, then $p$ does not fork over $M$ iff $p$ is invariant under 
$Aut({\bar M}/M)$, in other words $p$ does not split over $M$.
\end{Proposition}
\pf
(i) Right implies left: suppose $\phi(x,y)\in L$ and $(b_{i}:i<\omega)$ is an $A$-indiscernible
sequence of realizations of $tp(b/A)$ where $\phi(x,b)\in p$. We may assume that $b = b_{0}$. As $Lstp(b_{i}/A)
= Lstp(b_{0}/A)$ for all $i$, $\phi(x,b_{i})\in p(x)$ for all $i$, so trivially
$\{\phi(x,b_{i}):i<\omega\}$ is consistent.
\newline
Left implies right: Suppose first that $b_{0},b_{1}$ are the first two members of an $A$-indiscernible
sequence $(b_{i}:i<\omega)$, and $\phi(x,y)\in L$. We claim that 
$\phi(x,b_{0})\in p$ iff $\phi(x,b_{1})\in p$. If not then without loss of generality $\phi(x,b_{0})\wedge\neg\phi(x,b_{1}) 
\in p$. But note that $((b_{i},b_{i+1}):i=0,2,4,..)$
is also an $A$-indiscernible sequence. So as $p$ does not divide over $A$, 
$\{\phi(x,b_{i})\wedge\neg\phi(x,b_{i+1}):i = 0,2,4,..\}$ is consistent, but this contradicts $NIP$ (see Fact 
1.1.) So our claim is proved. Now if $Lstp(b/A) = Lstp(c/A)$ then we can find $b = b_{0}, b_{1},..,b_{n} = c$,
such that $(b_{i},b_{i+1})$ are the first two members of an $A$-indiscernible sequence, for each $i=0,..,n-1$.
So by our claim, $\phi(x,b)\in p$ iff $\phi(x,c)\in p$. This completes the proof of (i).

(ii) is immediate because types over models and Lascar strong strong types over models coincide.

\vspace{2mm}
\noindent
\begin{Definition} Let $p(x)\in S({\bar M})$ be a global type. 
\newline
(i) We say that $p$ is invariant over the small subset $A$ of ${\bar M}$ if $p$ is $Aut({\bar M}/A)$
-invariant. 
\newline
(ii) We say that $p$ is invariant if it is invariant over some small set.
\end{Definition}

Invariant types were studied by Poizat as ``special" types. 
By Proposition 2.1, if $T$ has $NIP$ then the invariant global types coincide with the global types which do not fork over 
some small set. If the global type $p$ is $A$-invariant then we have a kind of defining 
schema for $p$, namely for each $\phi(x,y)\in L$ we have some family $D_{p}{\phi}$ of complete $y$-types over 
$A$ 
such that for any $b\in {\bar M}$, $\phi(x,b)\in p$ iff $tp(b/A)\in D_{p}{\phi}$.
So we can apply the schema $D_{p}$ to not only supersets $B$ of $A$ living in ${\bar M}$ but also
to sets $B\supseteq A$ living in a proper elementary extension ${\bar M}'$ of ${\bar M}$. In any case for any such set $B$, 
by $p|B$ we mean the complete type over $B$ resulting from applying the schema $D_{p}$ to $B$. We will see 
subsequently that under the $NIP$ hypothesis the defining schema $D_{p}$ will be ``Borel". Given
invariant global types $p(x)\in S(\bar M)$, $q(y)\in (\bar M)$ we can form the product $p(x)\otimes q(y)\in 
S_{xy}(\bar M)$ as follows: Let $\phi(x,y)$ be over ${\bar M}$. We may assume $\phi(x,y)$ to be over small $A$ where both $p,q$ are $A$-invariant. We put $\phi(x,y)\in p(x)\otimes q(y)$ if for some (any) $b$ realizing $q(y)|A$, $\phi(x,b)\in p(x)$.

Alternatively, if we are willing to consider elements of some $|{\bar M}|^{+}$-saturated model containing 
${\bar M}$, define $p(x)\otimes q(y)$ to be $tp(a,b/{\bar M})$ where $b$ realizes $q(y)$ and $a$ realizes 
$p|({\bar M}b)$.

Note that if the global types $p(x)$, $q(y)$ are invariant, then so is $p(x)\otimes q(y)$. We see easily that $\otimes$ is associative. However it need not be commutative. Namely considering both $p(x)\otimes q(y)$ and $q(y)\otimes p(x)$ as elements of $S_{xy}({\bar M})$, they may not be equal. 

For an invariant global type $p(x)$, and disjoint copies $x_{1},..,x_{n}$ of the variable $x$ we define $p^{(n)}(x_{1},..,x_{n})$  inductively by: $p_{1}(x_{1}) = p(x_{1})$ and  
$p_{n}(x_{1},..,x_{n}) = p(x_{n})\otimes p_{n-1}(x_{1},..,x_{n-1})$. We let $p^{(\omega})(x_{1},x_{2},....)$ be the union of the $p_{n}(x_{1},..,x_{n})$ which will be a complete infinitary type over ${\bar M}$.

Assuming that $p(x)\in S(\bar M)$ is $A$-invariant, then by a {\em Morley sequence} in $p$ over $A$, we mean
a realization $(a_{1}, a_{2}, a_{3},...)$ in ${\bar M}$ of $p^{(\omega)}|A$. 

\begin{Lemma} Let $p(x)\in S({\bar M})$ be invariant. Then
\newline 
(i) Any realization $(b_{1}, b_{2},....)$ of $p^{(\omega)}$ (in an elementary extension of ${\bar M}$) is an indiscernible sequence over ${\bar M}$.
\newline
(ii) Suppose $A\subset {\bar M}$ is small and $p$ is $A$-invariant. If $a_{1}, a_{2},...$ from ${\bar M}$ are such that $a_{n+1}$ realizes $p|(Aa_{1},..,a_{n})$, then $(a_{1},a_{2},....)$ is a Morley sequence in $p$ over $A$. In particular $tp(a_{1},a_{2},.../A)$ depends only on $p$ and $A$.
\end{Lemma}
\pf Straightforward and left to the reader. 

\vspace{2mm}
\noindent
\begin{Remark} (Assume $NIP$) More generally we can define a Morley sequence of $p\in S({\bar M})$ over $A$, assuming just that $p$ does not fork over $A$, to be a realization in ${\bar M}$ of $Lstp(b_{1},b_{2},..../A)$ where $(b_{1},b_{2},..)$ realizes $p^{(\omega)}$ (in a model containing ${\bar M}$).  This is consistent with the previous definition.
\end{Remark}

\begin{Lemma} (Assume $NIP$) (i) Suppose $p(x)$, $q(x)$ are $A$-invariant global types. Then $p = q$ iff 
$p^{(\omega)}|A = q^{(\omega)}|A$ iff for all $n$ and realization $e$ of $p^{(n)}|A$, $p|Ae = q|Ae$. 
\newline
(ii) Suppose $Q(x_{0},x_{1},...)$ is the type over $A$ of some $A$-indiscernible sequence. Then
$Q = p^{(\omega)}|A$ for some $A$-invariant global type $p(x)$ if and only if whenever $I_{j}$ for $j\in J$ are
realizations of $Q$ then there is an element $c$ such that $(I_{j},c)$ is $A$-indiscernible for all $j\in J$.
\end{Lemma}
\pf (i) It suffices to prove that if $p|Ae = q|Ae$ for any realization $e$ of any $p^{(n)}|A$, then $p=q$.
Supposing for a contradiction that $p\neq q$ there is $\phi(x,b)\in p$, $\neg\phi(x,b)\in q$. Let $a_{1}, a_{2},...$ in ${\bar M}$ be such that $a_{i}$ realizes $p|(Aa_{1}..a_{i-1}b)$ for $i$ odd, and $a_{i}$ realizes 
$q|(Aa_{1}..a_{i-1}b)$ for $i$ even. Our assumption, together with Lemma 2.3(ii), implies that $(a_{1},a_{2},..)$ is a Morley sequence in $p$ over $A$, hence by 2.3(i) indiscernible over $A$. But $\phi(a_{i},b)$ holds iff $i$ is odd, contradicting Fact 1.1. 
\newline 
(ii) Left implies right is clear and does not require $NIP$. (Let $c$ realize $p|(A\cup\cup_{j}I_{j})$.)
For the other direction, assume $Q$ has the given property. Define the global type $p$ by: $\phi(x,b)\in p$ 
iff any realization $I$ of $Q$ extends to an indiscernible sequence $I'$ such that $\phi(x,b)$ is eventually 
true on $I'$. Then $NIP$ and our assumptions on $Q$ yield that $p$ is consistent, complete and $A$-invariant 
and
that $Q$ is the type over $A$ of its Morley sequence.

\vspace{5mm}
\noindent
We continue with some newer material. 
We first give a rather better and more general result on ``Borel definability" than that in \cite{NIP}. Given 
a 
(small) subset $A$ of $\bar M$, by a closed set over $A$ we mean the set of realizations in ${\bar M}$ of a 
partial type over $A$. An open set over $A$ is the complement (in the relevant ambient sort) of a closed set 
over $A$. From these we
build in the usual way the Borel sets over $A$. Alternatively these correspond to the Borel subsets of the 
relevant Stone space of complete types over $A$. A global type $p(x)$ will be called  ``Borel definable over 
$A$" if for any $L$-formula $\phi(x,y)$, the set of $b$ in ${\bar M}$ such that $\phi(x,b)\in p(x)$ is a Borel 
set over $A$. So if $p(x)$ is definable over $A$ in the usual sense then $p$ will be 
Borel definable over $A$ and if $p$ is Borel definable over $A$ then $p$ is $A$-invariant. In fact, we will be 
proving {\em strong} Borel definability over $A$, in the sense that for any $\phi(x,y)\in L$, the set of $b$ 
such that $\phi(x,b)\in p$ is a finite Boolean combination of closed sets over $A$. 

\begin{Proposition} (Assume $NIP$) Suppose that $p(x)\in S(\bar M)$
is a global type which is $A$-invariant. Then $p$ is strongly Borel definable over
$A$. 
\end{Proposition}
\pf
Let $\phi(x,y)\in L$. Let $N < \omega$ be as given for $\phi(x,y)$ by Fact 1.1. 
\newline
{\em Claim.} 
For any $b$, $\phi(x,b)\in p$ if and only if for some $n\leq N$ [there is
$(a_{1},..,a_{n})$ realising $p^{(n)}|A$ such that
\newline 
$\models \phi(a_{i},b)\leftrightarrow \neg\phi(a_{i+1},b)$ for $i=1,..,n-1$,  $(*)_{n}$
\newline
and
$\models \phi(a_{n},b)$,
\newline
but there is no $(a_{1},..,a_{n+1})$ such that $(*)_{n+1}$ holds].
\newline
{\em Proof of claim.} Suppose $\phi(x,b)\in p$. By Fact 1.1 choose any realization $(c_{i}:i<\omega)$
of $p^{(\omega)}|A$ with a maximal finite alternation (at most $N_{\phi}$) of truth values
of $\phi(c_{i},b)$ for $i<\omega$.  
Hence, eventually $\phi(c_{i},b)$ holds: for if not, let $c_{\omega}$ realize 
$p|(A\cup\{c_{i}:i<\omega\}\cup\{b\})$, and we contradict maximality. 

The converse holds by the above proof applied to $\neg\phi(x,b)$. So the claim is proved and clearly yields a strongly Borel definition of the set of $b$ such that $\phi(x,b)\in p$. 

\begin{Remark} (i) Define an $A$-invariant global type $p(x)$ to have $NIP$ if its Morley sequence
over $A$, $(b_{i}:i<\omega)$ (which has a unique type over $A$), has the property that for every 
$\phi(x,y)\in L$ there is $n_{\phi}< \omega$ such that for any $c$, there at most $n_{\phi}$
alternations of truth values of $\phi(b_{i},c)$. Then Proposition 2.6 goes through for $p$.
\newline
(ii) Proposition 2.6 also holds when $A$ is a set of hyperimaginaries, such as $bdd(B)$ for some
set $B$ of imaginaries.
\end{Remark}

We now consider (assuming $NIP$ still)  global types which do not fork over $A$, and Lascar strong types over $A$.

\begin{Remark} (Assume $NIP$) Let $p,q$ be global types which do not fork over $A$. Then $p\otimes q$ does not fork over $A$.
In particular $p^{(n)}$ and $p^{(\omega)}$ do not fork over $A$.
\end{Remark}
{\em Proof.} This follows from the well-known fact (valid for any theory $T$) that
if $tp(a/B)$ does not fork over $A$ and $tp(b/Ba)$ does not fork over $Aa$ then
$tp(a,b/B)$ does not fork over $A$. We will give a quick proof of this fact for completeness:
\newline
First, using the hypotheses,  find a saturated model $M$ containing $B$ such that $tp(a/M)$ does not fork over $A$, and $tp(b/Ma)$ does not fork over $Aa$. It is enough to prove:
\newline
{\em Claim.} $tp(ab/M)$ does not divide over $A$.
\newline
{\em Proof of Claim.}  Let $c_{0}\in M$, and $q(x,y,c_{0}) = tp(ab/Ac_{0})$. Let $(c_{i}:i<\omega)$ be an $A$-indiscernible sequence. We must show that $\cup\{q(x,y,c_{i}):i<\omega\}$ is consistent.
Let $p(x,c_{0}) = tp(a/Ac_{0})$. As $tp(a/M)$ does not fork over $A$, $\cup\{p(x,c_{i}):i<\omega\}$ is consistent, and we may assume (after applying an automorphism which fixes pointwise $Aac_{0}$) that $(c_{i}:i<\omega)$ is $Aa$-indiscernible. But now, as $tp(b/Ma)$ does not fork over $Aa$, we may find $b'$ realizing
$\cup\{q(a,y,c_{i}):i<\omega\}$. So then $(a,b')$ realizes $\cup\{q(x,y,c_{i}):i<\omega\}$, and we are finished.

\begin{Lemma} (Assume $NIP$) Suppose that $p(x)$ is a global type which does not fork over $A$.
Let $c,d$ realize $p|A$. Then $Lstp(c/A) = Lstp(d/A)$ iff there is an (infinite) sequence $a$ such that both 
$(c,a)$ and $(d,a)$ realize $p^{(\omega)}|A$.
\end{Lemma}
\pf Right implies left is immediate (for, as remarked earlier, elements of an infinite $A$-indiscernible 
sequence have the same Lascar strong type over $A$).
\newline
Left to right: Note first that for some $\sigma\in Aut({\bar M}/A)$, $\sigma(p)(x)$ implies $Lstp(c/A)$. 
As $\sigma(p)^{(\omega)}|A = p^{(\omega)}|A$, we may assume that already $p(x)$ implies $Lstp(c/A)$. Let $(a_{0},a_{1},a_{2},....)$ realize $p^{(\omega)}$ (in a
bigger saturated model). So $(a_{1},a_{2},....)$ also realizes $p^{(\omega)}$ and
does not fork over $A$. Hence by 2.1 whether or not $\phi({\bar x},c)$ is in $p^{(\omega)}$
depends on $Lstp(c/A)$. Hence $tp(c,a_{1},a_{2},.../A) = tp(a_{0},a_{1},..../A) =
tp(d,a_{1},a_{2},.../A)$ = $p^{(\omega)}|A$, as required.

\vspace{5mm}
\noindent
Lemma 2.9 says that on realizations of $p|A$, having the same Lascar strong type over
$A$ is a type-definable (over $A$) equivalence relation, hence by our discussion of Lascar strong types
in the introduction we see:
\begin{Corollary} (Assume $NIP$). (i) Suppose $p(x)\in S(A)$ does not fork over $A$.
Then on realisations of $p$, Lascar strong type over $A$ coincides with compact (KP) strong type over
$A$.
\newline
(ii) Suppose that $T$ is $1$-sorted (namely $T^{eq}$ is the $eq$ of a $1$-sorted theory). Suppose that
any complete $1$-type over any set $A$ does not fork over $A$. Then over any set $A$,
Lascar strong types coincide with compact strong types, hence $T$ is ``$G$-compact" over any
set of parameters.
\end{Corollary}
\pf (ii) The assumption, together with the discussion in the proof of 2.8,  implies that any 
complete type over any set $A$ does not fork over $A$. So we can apply (i).

\vspace{5mm}
\noindent

We can now strengthen Proposition 2.1. 
\begin{Proposition} (Assume $NIP$) Suppose that $p(x)$ is a global type. Then $p$ does
not fork over $A$ if and only if $p$ is $bdd(A)$-invariant.
\end{Proposition}
\pf Right to left is clear and does not use the $NIP$ assumption. For left to right, assume $p$ does not fork over $A$. Let $\sigma$
be an automorphism of $\bar M$ fixing $bdd(A)$ pointwise, and we have to show that $\sigma(p) = p$. By Remark 2.8 and Corollary 2.10(i) we have
\newline
{\em Claim I.} For any realization $\bar a$ of $p^{(\omega)}|A$, $Lstp({\bar a}/A) = 
Lstp(\sigma({\bar a})/A)$. 

\vspace{2mm}
\noindent
By Claim I and Proposition 2.1(i) it follows that
\newline
{\em Claim II.}  For any realization ${\bar a}$ of $p^{(\omega)}|A$, $p|A{\bar a} = \sigma(p)|A{\bar a}$.

\vspace{2mm}
\noindent
Now let $M$ be a small model containing $A$, and let $\bar a$ realize $p^{(\omega)}|M$.
\newline
{\em Claim III.} For any $c$ realizing either $p|A{\bar a}$ or $\sigma(p)|A{\bar a}$, ${\bar a}c$ is an $A$-indiscernible sequence.
\newline 
{\em Proof.}  By Claim II, it is enough to prove that ${\bar a}c$ is $A$-indiscernible for
$c$ realizing $p|A{\bar a}$. Note that $p$ does not fork over $M$ hence by 2.1(i) is $M$-invariant. So by Lemma 2.3, if $c$ realizes $p|M{\bar a}$ then ${\bar a}c$ is $M$-indiscernible, hence also $A$-indiscernible. But $c$ realizes 
$p|A{\bar a}$ and the latter is a complete type over $A{\bar a}$. So for {\em any} $c$ realizing $p|A{\bar a}$, ${\bar a}c$ is an $A$-indiscernible sequence. 

\vspace{2mm}
\noindent
{\em Claim IV.}  Let ${\bar a'}$ be an indiscernible sequence (of realizations of $p|A$) extending ${\bar a}$. Then for $c$ realizing either $p|A{\bar a'}$ or $\sigma(p)|A{\bar a'}$, ${\bar a'}c$ is $A$-indiscernible.
\newline
{\em Proof.} This can be seen in various ways. For example it can be deduced from Claim III, using the fact that each of $p$, $\sigma(p)$ is invariant under $Autf({\bar M}/A)$, and the fact that any two increasing $n$-tuples from ${\bar a'}$ have the same Lascar strong type over $A$.

\vspace{2mm}
\noindent
Now suppose for a contradiction that $p\neq \sigma(p)$. So for some $\psi(x,y)\in L$ and $e\in {\bar M}$,
$\psi(x,e)\in p$ and $\neg\psi(x,e) \in \sigma(p)$. Let $c_{i}$ realize
$p|A{\bar a}c_{0}...c_{i-1}e$ for $i$ even and realize $\sigma(p)|A{\bar a}c_{0}...c_{i-1}e$ for $i$ odd.
By Claim IV, $(\bar a,c_{i})_{i}$ is $A$ indiscernible. But $\models \psi(c_{i},e)$ iff $i$ is even, contradicting
$NIP$. This concludes the proof of Proposition 2.11.

\vspace{5mm}
\noindent
Finally we will give an analogue of 2.10 (ii) for strong types (which is closely related
to material in \cite{Ivanov-Macpherson} and \cite{Ivanov}). First a preparatory lemma.
\begin{Lemma} (Assume $NIP$). Suppose $A$ is algebraically closed, $tp(a/A)$ has a global $A$-invariant extension, 
and $e\in acl(Aa)$. Then $tp(ae/A)$ has a global $A$-invariant extension.
\end{Lemma} 
\pf Let $p(x) = tp(a/A)$ and fix some global $A$-invariant (so nonforking) extension $p'(x)$ of 
$p(x)$. Let $q(x,y) 
=tp(ae/A)$,
and let $q'(x,y)$ be any global extension of $q(x,y)$ whose restriction to $x$ is $p'(x)$. We will show that
$q'(x,y)$ is $A$-invariant. Let $\delta(x,y)$ be a formula over $A$ such that $\delta(a,y)$ isolates 
$tp(e/Aa)$.
\newline
{\em Claim I.} $q'$ does not fork over $A$. 
\newline
{\em Proof.} Let $\phi(x,y,b)\in q'(x,y)$ and let $(b_{i}:i<\omega)$ be $A$-indiscernible with
$b_{0} = b$. We may assume that $\models \phi(x,y,b)\rightarrow \delta(x,y)$. As $p'$ does not fork over $A$, there 
is $a'$ realizing $\{\exists y\phi(x,y,b_{i}):i<\omega\}$. For each $i$, let $e_{i}$
realize $\phi(a',y,b_{i})$. As there are finitely many possible choices for the $e_{i}$, there is an infinite 
subset $I$ of $\omega$ such that $e_{i} = e_{j}$ for $i,j\in I$. So $\{\phi(x,y,b_{i}):i\in I\}$ is 
consistent, which is enough.

\vspace{2mm}
\noindent
As $p'$ is $A$-invariant so is $p'^{(n)}$ for any $n$, hence:
\newline
{\em Claim II.} For all $n$, $p'^{(n)}|A$ implies a ``complete" Lascar strong type
over $A$. 

\vspace{2mm}
\noindent
{\em Claim III.} For any $n$,  $q'^{(n)}|A$ implies a ``complete" Lascar strong type over $A$.
\newline
{\em Proof.} As $q'^{(n)}$ does not fork over $A$, by Lemma 2.9 the relation of having the same Lascar strong 
type over $A$, on realizations of $q'^{(n)}|A$ is type-definable over $A$. But by Claim II clearly there are
only finitely many Lascar strong types over $A$ extending $q'^{(n)}|A$: If $((a_{1},e_{1}),..,(a_{n},e_{n}))$, 
and $((a_{1}',e_{1}'),..,(a_{n}',e_{n}'))$ are realizations of $q'^{(n)}|A$ with distinct Lascar strong types, we 
may by Claim II assume that $a_{i} = a_{i}'$ for $i=1,..,n$, so there are only finitely many possibilities
for the sequence $(e_{i})_{i}$. Hence equality of Lascar strong type on realizations of $q'^{(n)}|A$ is the 
restriction to $q'^{(n)}|A$ of a finite $A$-definable equivalence relation. As $A$ is algebraically closed, 
there is just one Lascar strong type over $A$ extending $q'^{(n)}$, proving Claim III. 

\vspace{2mm}
\noindent
The $A$-invariance of $q'$ follows from Claim III, as in the proof of 2.11.

\begin{Proposition} (Assume $NIP$) Let $T$ be $1$-sorted and work in $T^{eq}$. The following are equivalent:
\newline
(i) For any algebraically closed set $A$ and complete $1$-type $p$ over $A$ in the home sort, $p$ has a
global $A$-invariant extension.
\newline
(iii) For any complete type $p$ over any algebraically closed set $A$, $p$ has a global
$A$-invariant extension.
\newline
(iii) For any $A$ (a) any $p(x)\in S(A)$ does not fork over $A$, and (b) Lascar strong types over $A$ 
coincide with strong types over $A$.
\end{Proposition}
\pf (i) implies (ii). It is enough to prove that for any $n$ any complete $n$-type (in the home sort) over any
algebraically closed set $A$ has a global $A$-invariant extension. We prove it by induction on $n$. Suppose 
true for $n$. Let $p(x_{1},..,x_{n},x_{n+1})$ = $tp(a_{1},..,a_{n},a_{n+1}/A)$ with the $a_{i}$'s elements of
the home sort, and $A$ algebraically closed. Let (the infinite tuple) $e$ be an enumeration of 
$acl(Aa_{1},..,a_{n})$. By the induction 
hypothesis and 2.12, $tp(e/A)$ has a global $A$-invariant extension, realized by $e'$ say (in a bigger saturarated model
${\bar M}'$). Let $a_{1}',..,a_{n}'$ denote the copies of the $a_{i}$ in $e'$. By the hypothesis over the 
algebraically closed base $e'$, there is $a_{n+1}'$ (in the bigger model) with $tp(e',a_{n+1}'/A) = 
tp(e,a_{n+1}/A)$ such that $tp(a_{n+1}'/{\bar M}e')$ does not split over $e'$. It follows easily that 
$tp(e',a_{n+1}'/{\bar M})$ is $A$-invariant, hence also $tp(a_{1}',..,a_{n}',a_{n+1}'/{\bar M})$ is the 
$A$-invariant extension of $p$ we are looking for. 
\newline
(ii) implies (iii), and (iii) implies (i) are clear.

\begin{Corollary} If $T$ is $o$-minimal or $C$-minimal then (i) to (iii) of Proposition 2.13 hold.
\end{Corollary}
\pf Condition (i) holds in $C$-minimal theories through the existence of
``generic" $1$-types (see \cite{HHM2}). In the $o$-minimal case (i) holds without 
even the condition that $A$ is algebraically closed.

\section{Generically stable types}

Here we make a systematic study of what Shelah has called ``stable types" in \cite{Shelah715}. We discuss our choice of language a 
bit later. We begin with some preliminary remarks.
\begin{Lemma} (Assume $NIP$) Let $p(x)$ be a global type which does not fork over a small set $A$.
\newline
(i) Suppose $p$ is definable. Then $p$ is definable over $acl(A)$. In particular $p$ is $acl(A)$-invariant.
\newline
(ii) Suppose that $p$ is finitely satisfiable in some small model. Then $p$ is finitely
satisfiable in any model which contains $A$.
\end{Lemma}
\pf (i) By 2.11 $p$ is $bdd(A)$-invariant. So if $p$ is definable, then for any
$\phi(x,y)\in L$ the $\phi$-definition of $p$ is over $bdd(A)$ hence over $acl(A)$.
\newline
(ii) Let $M_{1}$ be a small model in which $p$ is finitely satisfiable.
Let $M$ be an arbitrary (small) model containing $A$. Let $\phi(x,c)\in p$. Let $M_{1}'$
realize a coheir of $tp(M_{1}/M)$ over $Mc$. As $p$ is $M$-invariant, $p$ is finitely satisfiable in $M_{1}'$
so there is $a'\in M_{1}'$ such that $\models\phi(a,c)$. So there is $a'\in M$ such that $\models\phi(a,c)$.

\vspace{5mm}
\noindent
Among our main results is:
\begin{Proposition} (Assume $NIP$) Let $p(x)\in S({\bar M})$, and let $A$ such that $p$ is $A$-invariant.  
\newline
Consider the conditions (i), (ii), (iii), (iv) and (v) below.
\newline
(i) $p(x)$ is definable (hence $A$-definable), and also finitely satisfiable in some/any small model containing $A$.
\newline
(ii) $p^{(\omega)}|A$ is totally indiscernible. That is, if $(a_{i}:i<\omega)$ is a Morley sequence
in $p$ over $A$, then $(a_{i}:i<\omega)$ is an indiscernible set (not just sequence) over $A$.
\newline 
(iii) For any formula 
$\phi(x,y)$ there is $N$ such that for any Morley sequence $(a_{i}:i<\omega)$ of $p$ over $A$, and $c$, 
$\phi(x,c)\in p$ if and only if 
\newline
$\models\vee_{w\subset 2N,|w|=N}\wedge_{i\in w}\phi(a_{i},c)$.
\newline 
(iv) For all small $B\supseteq A$, $p$ is the unique global nonforking extension of $p|B$,
\newline
(v) For all $n$, $p^{(n)}$ is the unique global nonforking extension of $p^{(n)}|A$.

\vspace{2mm}
\noindent
Then (i), (ii), (iii) and (iv) are equivalent, and imply (v). 
Moreover if $A$ has the additional property that every complete type over $A$ does not fork
over $A$, then (v) implies each of (i),(ii),(iii),(iv).
\end{Proposition}
\pf 
(i) implies (ii): Fix a small model $M\supseteq A$ such that $p$ is finitely satisfiable in $M$ (and of course definable over $M$), and 
there is no harm in proving (ii) with $M$ in place of $A$. Let $(a_{i}:i<\omega)$ be a Morley sequence in $p$
 over $M$.
We will show  
\newline 
(*) for any $n$ and $i\leq n$, $a_{i}$ realizes the restriction of $p'$
to $Ma_{0}..a_{i-1}a_{i+1}...a_{n}$. 
\newline 
Note that (*) will be enough to show by induction that for any $n$ and permutation $\pi$ of $\{0,..,n\}$,
$tp(a_{0},..,a_{n}/M) = tp(a_{\pi(0)},..,.a_{\pi(n)}/M)$ which, using indiscernibility of 
the {\em sequence} $(a_{i}:i<\omega)$ will prove its total indiscernibility.
\newline
So let us prove (*). Note that $tp(a_{i+1},..,a_{n}/Ma_{0},..,a_{i})$ is finitely satisfiable in $M$. As 
$tp(a_{i}/a_{0},..,a_{i-1},M)$ is definable over $M$, it follows that
$tp(a_{0}..a_{i-1}a_{i+1}..a_{n}/Ma_{i})$ is finitely satisfiable in $M$, whence 
$tp(a_{i}/Ma_{0}..a_{i-1}a_{i+1}..a_{n})$ is an heir of $p|M$, so (as $p$ is definable), realizes
$p|(Ma_{0}..a_{i-1}a_{i+1}..a_{n})$. This proves (*).

\vspace{2mm}
\noindent
(ii) implies (iii): This is by the ``in particular" clause of Fact 1.1.

\vspace{2mm}
\noindent
(iii) implies (i).  Clearly $p$ is definable. But it also follows from (iii) that $p$ is finitely satisfiable
in any model $M$ containing $A$. For suppose $\phi(x,c)\in p$. Let $I$ be a Morley sequence in $p$ over $A$
such that $tp(I/Mc)$ is finitely satisfiable in $M$. By (iii) $\phi(a,c)$ for some $a\in I$ hence $\phi(a,c)$ 
for some $a\in M$.

\vspace{2mm}
\noindent
(ii) implies (v). As (ii) for $p$ clearly implies (ii) for $p^{(n)}$, it suffices
to prove (v) for $n=1$, namely that $p$ is the unique global nonforking extension of $p|A$. 

Let us first note that because $p$ is $A$-invariant $p|A$ implies $p|bdd(A)$, and thus, by 2.11, {\em any} global nonforking extension of $p|A$ is $A$-invariant. 

Now let $q$ be an arbitrary global nonforking extension of $p|A$. Let $I = (a_{i}:i<\omega)$ be a Morley sequence in 
$p$ over $A$. We will prove that $I$ realizes $q^{(\omega)}|A$ which will be enough, by Lemma 2.5, to conclude 
that $p = q$. 
So let $b$ realize
$q|(A\cup I)$. We prove inductively that 
\newline
(**) $tp(a_{0},..,a_{n},b/A) = tp(a_{0},..,a_{n},a_{n+1}/A)$ for all $n$.
\newline
Note that of course $tp(a_{0}/A) = tp(b/A) = p|A$ which is in a sense the pre-base step.
\newline
So assume (**) is true for $n-1$. Suppose that $\models\phi(a_{0},..,a_{n},b)$. As $q$ does not fork over 
$A$ we have, by 2.1 and indiscernibility of $I$ that $\models\phi(a_{0},..,a_{n-1},a_{i},b)$ for all 
$i\geq n$. By (iii)
(for $p$) we see that $\phi(a_{0},..,a_{n-1},x,b)\in p(x)$. By the induction hypothesis and the 
$A$-invariance of $p$, we conclude that 
\newline $\phi(a_{0},..,a_{n-1},x,a_{n})\in p$, so
$\models\phi(a_{0},..,a_{n-1},a_{n+1},a_{n})$. Finally total indiscernibility of $I$ yields that
$\models\phi(a_{0},..,a_{n-1},a_{n},a_{n+1})$ as required.

\vspace{2mm}
\noindent
As condition (i) is preserved after replacing $A$ by any $B\supseteq A$, it 
follows from what we have proved up to now that each of (i),(ii), (iii) implies (iv).
\newline
\newline
(iv) implies (iii). Let $I = (a_{i}:i<\omega)$ be any Morley sequence in $p$ over $A$.
Note that $Av(I/{\bar M})$ is an $A\cup I$-invariant extension of $p|AI$ hence equals $p$. It follows easily
that for each $\phi(x,y)$ there is $N$ such that for any $c$ either
$|\{i<\omega:\models\phi(a_{i},c)\}| < N$ or $|\{i<\omega:\models\neg\phi(a_{i},c)\}| < N$. So we obtain (iii). 

\vspace{2mm}
\noindent
Now assume the additional hypothesis on $A$, and we prove
\newline
(v) implies (ii).
With some abuse of notation, let $p^{(\omega^{*})}(x_{0},x_{1},..)$ denote $tp(a_{0},a_{1},a_{2}...)/{\bar M})$
where for each $n$, $(a_{n-1},a_{n-2},..,a_{0})$ realizes $p^{(n)}$. Let $Q(x_{i})_{i \in\omega}$ be the 
restriction 
of $p^{(\omega^{*})}$ to $A$. So (iv) implies that 
\newline
(***) $p^{(\omega^{*})}$ is the unique global
nonforking extension of $Q$. 
\newline
Note that if $a$ realizes $p|A$ and $(a_{i})_{i}$ realizes 
the restriction of $p^{(\omega^{*})}$ to $Aa$ then $(a_{0},a_{1},...a_{n},.....,a)$ is $A$-indiscernible.
We claim that $Q$ satisfies the right hand side condition of 2.5 (ii), namely whenever $I_{1}$, $I_{2}$  
are realizations of $Q$ then there is $c$ such that $(I_{j},c)$ is indiscernible, for $j=1,2$. For, by our 
hypothesis on $A$, let $I_{1}',I_{2}'$ realize a global nonforking extension of $tp(I_{1},I_{2}/A)$. By (***),
each of $I_{1}', I_{2}'$ realizes $p^{(\omega^{*})}$. So choosing $c\in {\bar M}$ realizing $p|A$, realizing 
in ${\bar M}$ the restriction of $tp(I_{1}',I_{2}'/{\bar M})$ to $Ac$, and using an automorphism, gives the 
claim. By 2.5 (ii), $Q = q^{(\omega)}|A$ for some $A$-invariant global type $q$, which note must extend $p|A$.
Hence, by (iv), $q = p$. So if $I = (a_{0},a_{1},....)$ realizes $Q$ and $a$ realizes $q|AI$ then $a$ also 
realizes $p|AI$. Hence both $(a,I)$ and ($I,a)$ are indiscernible sequences over $A$, which easily implies 
that $I$ is an indiscernible set over $A$, giving (ii). The proof is complete.

\begin{Remark} (i) Assuming $NIP$, we will call a global type $p(x)$ {\em generically stable} if it is both definable and finitely satisfiable in some small model, namely $p$ satisfies 3.2(i) for some $A$.
We may also want to talk about generically stable types without a $NIP$ assumption, in which case we will mean a global $p$ such that for some $A$ (i), (ii), (iii) and (iv) of 3.2 are satisfied. We leave the reader to study the implications between (i)-(iv) in the absence of $NIP$.
\newline 
(ii) (Assume $NIP$.) Suppose the global type $p$ is generically stable, and $p$ does not fork over $B$. Then by 3.1, $p$ satisfies 3.2(i) with $A = acl(B)$. Hence by Proposition 3.2, we recover the finite equivalence relation theorem: any two global nonforking extensions of $p|B$ are distinguished by some finite $B$-definable equivalence relation.
\newline 
(iii) Proposition 3.2 goes through assuming only that $p$ has $NIP$ (as in Remark 2.7).
\newline
(iv) In $ACVF$, for any set $A$, any complete type over $A$ has a global nonforking extension. (See \cite{HHM2} or Proposition 2.13.) Hence (i)-(v) of 3.2 are equivalent in $ACVF$, for any $A$. 
\newline
(v) In 3.2(iv) it is not enough to require just that $p|A$ has a unique global nonforking extension.
\end{Remark}
\noindent
{\em Explanation of (v).} So we give an example of a $NIP$ theory and a type $p(x)\in S(A)$ with a unique global
$A$-invariant extension which is {\em not} definable. We consider the basic $C$-minimal theory consisting of 
a dense linear ordering $(I,<)$ with greatest element $\infty$, and another sort on which there are 
equivalence relations $E_{i}$ indexed (uniformly) by $i\in I$ and with $E_{j}$ infinitely refining $E_{i}$ if 
$i<j$ (plus some other axioms, see \cite{Hrushovski-Kazhdan}). By a definable ball we mean an $E_{i}$-class 
for some $E_{i}$. By a type-definable ball we mean a possibly infinite intersection of definable balls. We 
can produce a model $M$ and a type-definable over $M$ ball $B$ such that $B$ contains no proper $M$-definable 
ball. Let the global type $p_{B}(x)\in S({\bar M})$ be the ``generic type" of $B$, namely $p$ says that $x\in 
B$ and $x$ is not in any proper definable sub-ball of $B$. Then $p$ is the unique $M$-invariant extension of 
$p|M$ but is not definable. 

\vspace{5mm}
\noindent
In the $NIP$ context, our generically stable types coincide with what Shelah \cite{Shelah715} calls {\em stable} 
types. However there is already another meaning for a complete or even partial type $\Sigma(x)$ to be stable. 
It is that {\em any} extension of $\Sigma(x)$ to a complete global type is definable (over some set, not 
necessarily the domain of $\Sigma$). This notion is also sometimes called ``stable, stably embedded" (although mainly in the case where the partial type is a single formula). One family of examples of 
generically stable types come through stable domination in the sense of \cite{HHM2}. Recall that $q(x)\in S(A)$ ($A$ algebraically closed) is said to be stably dominated if there is stable partial type $\Sigma(x)$ over $A$ (stable in the strong sense mentioned above), and an $A$-definable function $f$ from the set of realizations of $p$ to the realizations of $\Sigma$ such that, if $q\in S(A)$ is $f(p)$, and $a$ realizes $p$, then whenever $f(a)$ is independent from $B$ over $A$ then $tp(a/A,f(a))$ has a unique extension over $Bf(a)$. A stably dominated type is generically stable, as is easily verifiable. Generically stable but unstable types occur in algebraically closed valued fields through stable domination, namely via the stable part of the structure which is essentially the residue field. In ``pure" unstable $NIP$ theories such as $o$-minimal and weakly $o$-minimal theories, or $p$-adically closed fields there are no (nonalgebraic) generically stable types. On the other hand in simple theories any stationary type (type with a unique nonforking extension) is easily seen to be generically stable.

\vspace{2mm}
\noindent
We finish this section with some remarks on invariant types and symmetry. 

\begin{Lemma} Suppose $p(x), q(y)$ are global types such that $p$ is finitely satisfiable in some small model, and $q$ is definable. (So both $p$ and $q$ are invariant.) Then $p_{x}\otimes q_{y} = q_{y}\otimes p_{x}$.
\end{Lemma}
\pf
Let $M$ be a small model such that $q$ is definable over $M$ and $p$ is finitely satisfiable
in $M$. Let $(a,b)$ realize $(p(x)\otimes q(y))|M$, namely $b$ realizes $q|M$ and $a$ realizes $p|Mb$. We want to show that $(a,b)$ realizes $(q(y)\otimes p(x))|M$. Suppose not. Then there is a formula $\phi(a,y)\in q$ such that $\models\neg\phi(a,b)$. Let $\psi(x)$ (a formula over $M$) be the $\phi(x,y)$-definition for $q(y)$. So 
$\models\neg\psi(a)$. As $tp(a/Mb)$ is finitely satisfiable in $M$, there is $a'\in M$ such that 
$\models\psi(a')\wedge\neg\phi(a',b)$, which is a contradiction as $tp(b/M) = q|M$. (The reader should note that this is just a restatement of uniqueness of heirs for definable types.)

\vspace{5mm}
\noindent
Finally we return to generically stable types:
\begin{Proposition} Suppose $T$ has $NIP$, and that $p(x)$, $q(y)$ are global invariant types such that $p(x)$ is generically stable. Then $p_{x}\otimes q_{y} = q_{y}\otimes p_{x}$
\end{Proposition}
{\em Proof.} Suppose $\phi(x,y)\in L_{\bar M}$ and $\phi(x,y)\in p(x)\otimes q(y)$. Let $M$ be a small model such that each of $p$ and $q$ do not fork over $M$ and $\phi$ is over $M$. By assumption, $p$ is definable over $M$. Moreover 
\newline
(*) if $\psi(y)$ is the $\phi(x,y)$-definition of $p$, then for any realisation $c$ of $\psi$, and any Morley sequence $(a_{i}:i<\omega)$ of $p$ over $M$, $\models\psi(a_{i},c)$ for all but finitely many $i<\omega$. 

\vspace{2mm}
\noindent
From our assumption that $\phi(x,y)\in p(x)\otimes q(y)$, it follows that $\psi(y)\in q$. Now suppose for a contradiction that $\neg\phi(x,y)\in q(y)\otimes p(x)$. Let $(a_{i}:i<\omega)$ be a Morley sequence in $p$ over $M$. So then $\neg\phi(a_{i},y)\in q$ for all $i$, so there is $b$ realizing $q|(M\cup\{a_{i}:i<\omega\})$ such
that $\models\neg\phi(a_{i},b)$ for all $i$. As $b$ realizes $\psi(y)$ this contradicts (*), and completes the proof.

\section{Measures and forking}
In this section we generalize some of the results on types in the previous sections to Keisler measures.

Recall that a Keisler measure $\mu$ on a sort $S$ over a set of parameters $A$, is a finitely additive probability measure on $A$-definable subsets of $S$ (or on formulas over $A$ with free variable in sort $S$), namely $\mu(X) \in [0,1]$ for all $A$-definable $X$, $\mu(S) = 1$ and the measure of the union of two disjoint $A$-definable sets $X$ and $Y$ is the sum of the measures of $X$ and of $Y$.  A complete type over $A$ (in sort $S$) is a special case of a measure. We sometimes write $\mu(x)$ or $\mu_{x}$ to mean that the measure is on the sort ranged over by the variable $x$. A Keisler measure over ${\bar M}$ is called a global (Keisler) measure. 

\vspace{2mm}
\noindent
Let us emphasize that a Keisler measure $\mu_{x}$ over $A$ is the same thing as a
regular Borel probability measure on the compact space $S_{x}(A)$.
Regularity means that for any Borel subset $B$ of $S(A)$,
and $\epsilon > 0$, there are closed $C$ and open $U$ such that $C\subseteq B \subseteq U$ and
$\mu(U\setminus C) < \epsilon$. Note that $\mu_{x}$ defines a finitely additive probability measure (still called $\mu$)  on the
algebra of clopens of $S(A)$. Theorem 1.2 of \cite{Keisler1} extends $\mu$ to a Borel probability measure 
$\beta$ on $S_{x}(A)$ using the Loeb measure construction in a mild way. And Lemma 1.3(i) of \cite{Keisler1}
says that this $\beta$ is regular. On the other hand a Borel probability measure $\beta$ on $S_{x}(A)$ gives, by
restriction to the clopens, a Keisler measure $\mu_{x}$ over $A$. Moreover if $\beta$ is regular then 
for any closed set $C$, $\mu(C)$ will be the infimum of the $\mu(C')$ where $C'$ is a clopen containing $C$.
Hence $\beta$ is determined by $\mu$. 

\vspace{2mm}
\noindent 
An important, and even characteristic, fact about $NIP$ theories is that for any global Keisler measure $\mu_{x}$ there are only boundedly many definable sets up to $\mu_{x}$-equivalence. (See Corollary 3.4 of \cite{NIP}.)

\vspace{2mm}
\noindent
As in \cite{NIP} notions relating to types generalize naturally to measures. There we discussed the notions of a measure being definable, and of being finitely satisfiable: For example if $\mu$ is a global Keisler measure then $\mu$ is definable over $A$, if for each closed subset of $[0,1]$ and $L$-formula $\phi(x,y)$,
$\{b\in {\bar M}: \mu(\phi(x,b)) \in C\}$ is type-definable over $A$. We say $\mu$ is finitely satisfiable in $M$ if every formula with positive measure is realized by a tuple from $M$.
But we also have the notion of forking:
\begin{Definition} Let $\mu$ be a Keisler measure over $B$, and $A\subseteq B$. We say that $\mu$ does not divide over 
$A$ if whenever $\phi(x,b)$ is over $B$ and $\mu(\phi(x,b)) > 0$ then $\phi(x,b)$ does not divide over $A$. Similarly we say $\mu$ does not fork over $A$ if every formula of positive $\mu$-measure does not fork over $A$.
\end{Definition}

\begin{Remark} Suppose $\mu$ is a global Keisler measure.
\newline
(i) $\mu$ does not fork over $A$ iff $\mu$ does not divide over $A$. 
\newline 
(ii) If $\mu$ is either definable over $M$ or finitely satisfiable in $M$, then $\mu$ does not fork over $M$.
\end{Remark}
\pf  (i) If $\mu(\phi(x)) > 0$ and $\phi$ forks over $A$ then $\phi$ implies a finite disjunction of formulas
each of which divides over $A$. One of those formulas must have positive measure by finite additivity. 
\newline 
(ii) If $\mu$ is either definable over or finitely satisfiable in $M$, then for any $\phi(x,y)\in L$ and $b$, $\mu(\phi(x,b))$ depends on $tp(b/M)$. (In the case that $\mu$ is definable over $M$ this is immediate. If $\mu$ is finitely satisfiable in $M$, then $tp(b_{1}/M) = tp(b_{2}/M)$ implies that the measure of the symmetric difference of $\phi(x,b_{1})$ and $\phi(x,b_{2})$ is $0$, so again $\phi(x,b_{1})$ and $\phi(x,b_{2})$ have the same measure.) So if $(b_{i}:i<\omega)$ is an $M$-indiscernible sequence and $\mu(\phi(x,b_{0})) = r> 0$, then for all $i$, $\mu(\phi(x,b{_i})) = r$. By Lemma 2.8 of \cite{NIP}, $\{\phi(x,b_{i}):i<\omega\}$ is consistent.

\vspace{5mm}
\noindent
Proposition 2.1 readily generalizes to measures.
\begin{Proposition} (Assume $NIP$)  Suppose $\mu$ is a global Keisler measure, and $A$ a small set. Then the following are equivalent:
\newline
(i) $\mu$ does not fork over $A$,
\newline
(ii) $\mu$ is invariant under $Autf({\bar M}/A)$, 
\newline
(iii) whenever $Lstp(b_{1}/A) = Lstp(b_{2}/A)$, then $\mu(\phi(x,b_{1})\Delta\phi(x,b_{2})) = 0$.
\end{Proposition}
\pf  (i) implies (iii). Suppose that $\mu$ does not fork over $A$. Let $b_{0},b_{1}$ begin an indiscernible sequence $(b_{i}:i<\omega)$. We claim that
$\mu(\phi(x,b_{0})\Delta\phi(x,b_{1})) = 0$. If not then without loss 
$\mu(\phi(x,b_{0})\wedge\neg\phi(x,b_{1})) > 0$. As $((b_{i},b_{i+1}):i=0,2,4,..)$ is $A$-indiscernible and $\mu$ does not fork over $A$,
$\{\phi(x,b_{i})\wedge\neg\phi(x,b_{i+1}):i=0,2,...\}$ is consistent, contradicting $NIP$. So we clearly 
obtain (iii). 
\newline
(iii) implies (ii) is immediate.
\newline
(ii) implies (i). Assume (ii) and suppose that $\mu(\phi(x,b)) = \epsilon > 0$ and that
$(b = b_{0},b_{1},.......)$ is $A$-indiscernible. So the $b_{i}'$ have the same Lascar strong type over $A$, hence $\mu(\phi(x,b_{i})) = \epsilon$ for all $i$. By Lemma 2.8 of \cite{NIP} again, 
$\{\phi(x_, b_{i}):i<\omega\}$
is consistent.

\vspace{5mm}
\noindent
So we see that the global Keisler measures which do not fork over some small set coincide with those which are invariant over some small model, and we can just call them invariant global measures (as we did for types). There is an obvious notion of Borel definability for a measure. Namely we say that (global) $\mu$ is Borel definable over $A$ if for any $\phi(x,y)$ and closed subset $C$ of $[0,1]$ the set of $b\in \bar M$ such that $\mu(\phi(x,b))\in C$ is Borel over $A$. We will prove at the end of this section that (assuming $NIP$) any global invariant Keisler measure is Borel definable.

We can therefore define a product on invariant measures (global Keisler measures invariant over some small set of model) by integration: Suppose $\mu_{x}$ is a Borel definable (over $M$) global Keisler measure and $\lambda_{y}$ is another global
Keisler measure. Fix $\phi(x,y)\in L_{\bar M}$. For any $b\in \bar M$ let $f(b) = \mu(\phi(x,b))$ Then $f$ is a Borel function over $M$ on the sort to which the $b$'s belong, so we can form
the integral $\int f(y) d\lambda_{y}$, and we call this $(\mu\otimes \lambda)(\phi(x,y))$. When $\mu$, $\lambda$ are
global complete types, this agrees with the product as defined in section 2.

\begin{Remark} It is natural to attempt to generalize the material from section 3 to measures. So 
(assuming $NIP$) we call a global Keisler measure $\mu$ {\em generically stable} if $\mu$ is finitely satisfiable in some small model, and also definable. The analogues of Proposition 3.2 and Proposition 3.5 will be treated in a future work joint with P. Simon. Likewise for the measure analogue of the symmetry Lemma 3.4.
\end{Remark}

\vspace{2mm}
\noindent
For now, we continue the generalization of the results from section 2 to measures. We also begin to relate invariant measures (those measures which do not fork over a small set) to invariant types. We first generalize Proposition 2.11.

\begin{Proposition} (Assume $NIP$) Let $\mu_{x}$ be a global Keisler measure which does not fork over $A$. Then $\mu$ is
$bdd(A)$-invariant.
\end{Proposition}
\pf  Let $B$ be the Boolean algebra of formulas $\phi(x,b)$ over ${\bar M}$ quotiented by the equivalence 
relation $\phi(x,b) \sim \psi(x,c)$ if $\mu(\phi(x,b)\Delta\psi(x,c)) = 0$. Let $[\phi(x,b)]$ denote the class
of $\phi(x,b)$, i.e. the image of $\phi(x,b)$ in $B$. 
Note that $B$ is small, because
by Proposition 4.3, 
\newline
(*) if $Lstp(b/A) = Lstp(c/A)$ then $[\phi(x,b)] = [\phi(x,c)]$. 
\newline
For any ultrafilter $U$ on 
$B$, we obtain a complete global type $p_{U}(x)$ by putting $\phi(x,b)\in p_{U}$ iff $[\phi(x,b)]\in U$.
Note that by (*) $p_{U}$ is invariant under $Autf({\bar M}/A)$, so does not fork over
$A$, so is in fact $bdd(A)$-invariant.

We claim that in fact if $b$ and $c$ have the same compact strong type over $A$ (namely
the same type over $bdd(A)$) then $\mu(\phi(x,b)\Delta\phi(x,c)) = 0$, which will show
that $\mu$ is $bdd(A)$-invariant. If not we have, without loss, $\mu(\phi(x,b)\wedge\neg\phi(x,c)) > 0$,
so there is an ultrafilter $U$ containing $[\phi(x,b)\wedge\neg\phi(x,c)]$, and we obtain
the global type $p_{U}$ which as pointed out above is $bdd(A)$-invariant. As $\phi(x,b)\in p_{U}$, also 
$\phi(x,c)\in p_{U}$, a contradiction.

\vspace{5mm}
\noindent
Note that given small $A$ the collection of global $x$-types which do not fork over $A$ is closed
in $S_{x}({\bar M})$ and coincides with the space of global types invariant over $bdd(A)$. 
As in section 8 we call this space $S_{x}^{bdd(A)}({\bar M})$. In the same vein as Proposition 4.5 we have.
\begin{Proposition} Fix small $A$. Then there is a natural bijection between global Keisler measures which do 
not fork over $A$ and regular Borel probability measures on $S_{x}^{bdd(A)}({\bar M})$.
\end{Proposition}
\pf We have already observed the bijection between global Keisler measures and regular Borel probability
measures on the space of global types. Let $\mu_{x}$ be a global Keisler measure which does not fork over 
$A$, and let $\beta$ be the regular Borel measure on $S_{x}({\bar M})$ corresponding to it. It is enough to 
show that the restriction of $\beta$ to $S_{x}^{bdd(A)}({\bar M})$ determines $\mu$. But this is clear. For 
if $\phi(x), \psi(x)$ are formulas over ${\bar M}$ and $\mu(\phi\Delta\psi) > 0$, then there is some 
$p\in S_{x}^{bdd(A)}({\bar M})$ containing (without loss) $\phi\wedge\neg\psi$.

\begin{Proposition} (Assume $NIP$) Let $p(x)$ be a complete type over $A$. Then the following are equivalent:
\newline
(i) $p$ does not fork over $A$ (that is $p$ has a global nonforking extension),
\newline
(ii) there is a global Keisler measure $\mu$ which extends $p$ and is $A$-invariant. 
\end{Proposition}
{\em Proof.} (ii) implies (i): Let $\mu$ be a global $A$-invariant Keisler measure which extends $p$. Suppose, for a contradiction, that $p(x)$ forks over $A$. So there are $\phi(x)\in p(x)$, and formulas
$\psi_{1}(x),..,\psi_{n}(x)$  (with parameters), each dividing over $A$ such that $\models \phi(x)\rightarrow\vee_{i}\psi_{i}(x)$. As $\mu(\phi(x)) = 1$, $\mu(\psi_{i}(x)) > 0$ for some
$i=1,..,n$. But as $\mu$ is $A$-invariant, $\mu$ does not divide over $A$, contradiction. 
\newline 
(i) implies (ii). Let $p'(x)$ be a global nonforking extension of $p$. By 2.11 $p'$ is $bdd(A)$-invariant. By Proposition 2.6, $p'$ is Borel definable over $bdd(A)$. Let $G$ be the compact Lascar group over 
$A$, namely $Aut(bdd(A)/A)$. $G$ is a compact group which thus has a (unique left and right invariant) Haar measure $h$. Let $S$ be the space of all complete types over $bdd(A)$. 
So $S$ is a compact space acted on continuously by $G$. 

We now define the $A$-invariant global measure $\mu$ extending $p(x)$. Let $\phi(x,b)$ be a formula over $\bar M$. 
Let $q(y) = tp(b/A)$, and let $Q\subset S$ be the space of complete types over $bdd(A)$ extending $q$ (a closed subset of $S$). The Borel definability of $p'$ over $bdd(A)$ says that the set $X$ of $tp(b'/bdd(A))$ such that $\phi(x,b')\in p'$ is a Borel subset of $S$. Hence $X\cap Q$ is also a Borel subset of $Q$. The compact space $Q$ is acted on continously and transitively by $G$, hence has the form $G/H$ for some closed subgroup of $G$. Let $\pi$ be the canonical projection from $G$ onto $Q$. As $X\cap Q$ is Borel in $Q$, $\pi^{-1}(X\cap Q)$ is a
Borel subset of $G$, hence $h(\pi^{-1}(X\cap Q))$ is defined and we define this to be $\mu(\phi(x,b))$. Note that if $tp(b'/A) = tp(b/A) ( = q)$ then by construction $\mu(\phi(x,b')) = \mu(\phi(x,b))$, so $\mu$ is $A$-invariant.
Also if $\phi(x,b) \in p(x)$ (namely $b\in A$), then $Q$ is a singleton ($\{tp(b/A)\}$) and $X\cap Q = Q$, whereby $\mu(\phi(x,b)) = 1$. So $\mu$ extends $p$. We leave it to the reader to check finite additivity of $\mu$.

\vspace{2mm}
\noindent
In the stable case it is not so hard to see, via the finite equivalence relation theorem for example,  that 
the $A$-invariant measure $\mu$ extending $p$ in (ii) above is unique. Likewise if $p$ has a generically 
stable global nonforking extension. It would be interesting to find examples of uniqueness of $\mu$, when $p$ 
does not have a generically stable nonforking extension. We will come back to the issue of when $\mu$ is 
unique later in the paper.

\vspace{2mm}
\noindent
We now return to the question of the Borel definability of invariant measures.
We will make use of the Vapnik-Chervonenkis theorem \cite{Vapnik-Chervonenkis} which we now describe.

Let $(X,\Omega,\mu)$ be a probability space, that is a set $X$ equipped with a $\sigma$-algebra
$\Omega$ of subsets (events) of $X$ and a countably additive probability measure $\mu$ with values in $[0,1]$. Namely
$\mu(\emptyset) = 0$, $\mu(X) = 1$, $\mu(E)$ is defined for any $E\in \Omega$ and if $E_{i}\in \Omega$ for $i<\omega$ are pairwise disjoint and $E = \cup_{i}E_{i}$ then $\mu(E) = \sum_{i}\mu(E_{i})$. For $k> 0$, let $\mu^{k}$ be the product measure on $X^{k}$. Also, for $k > 0$, $A\in {\Omega}$ and ${\bar p} = (p_{1},..,p_{k})\in X^{k}$, let
$fr_{k}(A, {\bar p})$ be the proportion of $p_{i}$'s which are in $A$, namely $|\{i:p_{i}\in A\}|/n$. 
For any $A\in \Omega$, let $g_{A,k}: X^{k}\to [0,1]$ be defined by 
$g_{A,k}({\bar p}) = |fr_{k}(A,{\bar p}) - \mu(A)|$. Also let $h_{A,k}:X^{2k}\to [0,1]$ be defined by
$h_{A,k}({\bar p},{\bar q}) = |fr_{k}(A,{\bar p}) - fr_{k}(A,{\bar q})|$. Note that $g_{A,k}$ is $\mu^{k}$-measurable and $h_{A,k}$ is $\mu^{2k}$-measurable. We will say (somewhat nonstandardly) that a family $\cal C$ of events has $NIP$ if there is a natural number $d$ such that for no subset $F$ of $X$ of cardinality $d$ do we have that $\{F\cap A: A\in {\cal C}\}$ is the full power set of $F$. 
For ${\cal C}$ with $NIP$ the smallest such $d$ will be called the $VC$-dimension of ${\cal C}$.

With this notation, Theorem 2 of \cite{Vapnik-Chervonenkis} gives:

\vspace{2mm}
\noindent
{\em VC Theorem.} Suppose that the family $\cal C$ of events has $NIP$. Suppose also that for every $k>0$,
$sup_{A\in{\cal C}}g_{A,k}$ is $\mu^{k}$-measurable and $sup_{A\in{\cal C}}h_{A,k}$ is $\mu^{2k}$-measurable. Then there is a function $f(-,-)$ such that
for any $\epsilon > 0$, and $k < \omega$ 
$\mu^{k}(\{{\bar p}: sup_{A\in{\cal C}}g_{A,k}({\bar p}) > \epsilon\}) \leq f(k,\epsilon)$, and for any given 
$\epsilon$, $f(k,\epsilon) \to 0$ as $k\to\infty$.
\newline
Moreover the function $f$ depends only on the $VC$-dimension of $\cal C$.

\vspace{3mm}
\noindent
It is convenient to state a version or consequence of the VC Theorem in which the two measurability assumptions are dropped.

First let us note that any family $\cal F$ of measurable functions from $X$ to $[0,1]$ say, has an {\em essential supremum}, which is by definition measurable, and can be chosen as the outright supremum of some countable subfamily of ${\cal F}$. By an essential supremum of ${\cal F}$ we mean a measurable function $g:X\to [0,1]$ such that for each 
$f\in {\cal F}$, $g\geq f$ on a set of $\mu$-measure $1$, and whenever $g':X\to [0,1]$ has the same property, then $g\leq g'$ on a set of measure $1$. Such a $g$ can be found
as follows. We construct measurable $g_{\alpha}$ for $\alpha$ a countable ordinal in the following way: let $g_{0}$ be some member of ${\cal F}$. At limit stages, take suprema. Given $g_{\alpha}$, if it is already an essential supremum of ${\cal F}$, stop. Otherwise there is some $f\in {\cal F}$ such that $f>g_{\alpha}$ on a set of positive measure. Put $g_{\alpha + 1} = sup\{g_{\alpha},f\}$. Note that $\int g_{\alpha}$ is strictly increasing with $\alpha$, hence the construction has to stop at some countable ordinal $\beta$. Then $g_{\beta}$ is measurable, is an essential supremum of ${\cal F}$, and is, by construction, the supremum of some countable subfamily of ${\cal F}$.

Note that any two essential suprema $g, g'$ of ${\cal F}$ are equivalent in the sense that they are equal almost everwhere. We will just write $ess sup {\cal F}$ for a representative of the equivalence class which we will assume to be a supremum of a countable subfamily of ${\cal F}$. 

\vspace{3mm}
\noindent
{\em VC Theorem}*. Suppose that the family $\cal C$ of events has $NIP$.  
Then for any $\epsilon > 0$,  
$\mu^{k}(\{{\bar p}: ess sup_{A\in{\cal C}}g_{A,k}({\bar p}) > \epsilon\}) \to 0$ as $k\to\infty$.

\vspace{2mm}
\noindent
This *-version is immediately deduced from the first version: Let us fix $\epsilon$. We can find a  countable subfamily ${\cal C}_{0}$ of ${\cal C}$ with the same $VC$-dimension, and such that for each $k$,
$ess sup_{A\in{\cal C}}g_{A,k}$ = $sup_{A\in{\cal C}_{0}}g_{A,k}$. Note that $sup_{A\in{\cal C}_{0}}h_{A,k}$ is measurable. Hence we can apply the $VC$ Theorem with ${\cal C}_{0}$ in place of ${\cal C}$ and we obtain the *-version.



\vspace{5mm}
\noindent
We will apply the above theorems to the following situation. Let $M$ be a model of a $NIP$ theory, 
$\phi(x,y)$ an $L$-formula, and $\mu_{x}$ a Keisler measure over $M$. Take $X$ to be the Stone space $S_{x}(M)$, 
$\Omega$ the $\sigma$-algebra of its Borel subsets, and identify $\mu_{x}$ with the measure induced on $X$ as described at the beginning of this section. Let $\cal C$ be the family of clopen subsets of $X$ given by the formulas $\phi(x,c)$ as $c$ varies over $M$. It is easy to see that as $T$ has $NIP$ the family ${\cal C}$ has $NIP$.

\begin{Lemma} (Assume $NIP$.) 
\newline
Let $M$ be a model and $\mu_{x}$ a Keisler measure over $M$.
\newline
(i) Let $\phi(x,y)\in L$, and let $\epsilon > 0$. Then there is $k$, and a Borel subset $B$ of $S_{x}(M)^{k}$ of positive $\mu^{k}$ measure such that for all $(p_{1},..,p_{k})\in B$ and any $c\in M$, $\mu(\phi(x,c))$ is within $\epsilon$ of the proportion of the $p_{i}$ which contain $\phi(x,c)$.
\newline
(ii) If moreover $M$ is saturated (for example equals ${\bar M}$) and $A\subseteq M$ is small such that $\mu$ does
not fork over $A$, then in (i) we can choose $B$ such that for any $(p_{1},..,p_{k})\in B$, each $p_{i}$ does not fork over $A$.
\end{Lemma}
{\em Proof.} (i) We will identify $\phi(x,c)$ with the clopen subset of $X$ it defines, and as above write $fr_{k}(\phi(x,c), {\bar p})$ for the proportion of $p_{i}$ containing $\phi(x,c)$. We fix our $\epsilon > 0$. By {\em VC Theorem}*, choose $k$ and $\delta> 0$ such that 
$\mu^{k}(\{{\bar p}: ess sup_{A\in{\cal C}}g_{A,k}({\bar p}) < \epsilon\}) > \delta$.
It follows that for each finite  set $C = \{c_{1},..,c_{m}\}$ of parameters from $M$,
the Borel set $B_{C}$ = $\{(p_{1},..,p_{k})\in S(M)^{k}: fr_{k}(\phi(x,c_{j}),{\bar p})-\mu(\phi(x,c_{j})) < \epsilon$, for $j=1,..,m \}$ has $\mu^{k}$-measure $>\delta$. Note that $B_{C'}\subseteq B_{C}$ when $C\subseteq C'$. As Borel sets are approximated in measure by closed sets, we can use compactness to find a Borel (in fact closed) subset $B$ of $S_{x}(M)^{k}$ of $\mu^{k}$-measure $\geq \delta$ such that {\em all} $(p_{1},..,p_{k})\in B$ satisfy the conclusion of (i).

\vspace{2mm}
\noindent
(ii) This is proved exactly as (i) but working in the closed set $S_{x}^{bdd(A)}(M)$ of $S_{x}(M)$ consisting of, types over $M$ which do not fork over $A$, and bearing in mind 4.6.

\vspace{2mm}
\noindent
\begin{Corollary} (Assume $NIP$.) Suppose that $\mu$ is a global Keisler measure which does not fork over $A$. Then $\mu$ is Borel definable over $bdd(A)$.
\end{Corollary}
{\em Proof.} Fix $\phi(x,y)\in L$, and $\epsilon > 0$. Let $p_{1},..,p_{n}$ be given by Lemma 4.8 (ii) (for $M = {\bar M}$).
By Proposition 2.6, each $p_{i}$  is strongly Borel definable over $bdd(A)$. Namely for each $i$ there is some finite Boolean 
combination $Y_{i}$ of type-definable over $bdd(A)$-sets such that for all $c$, $\phi(x,c)\in p_{i}$ iff 
$c\in Y_{i}$. It follows that if $c,c'$ lie in exactly the same $Y_{i}$ for $i=1,..,n$,
then $\mu(\phi(x,c))$ and $\mu(\phi(x,c'))$ differ by $< 2\epsilon$. This suffices to prove Borel definability of $\mu$.

\vspace{5mm}
\noindent
A thorough account of the model-theoretic results related to VC-type theorems will appear in a forthcoming paper.

\section{Generics and forking in definable groups}
In \cite{NIP} the $fsg$ (finitely satisfiable generics) property for definable groups was 
introduced. Strengthenings of the $fsg$ will be discussed in section 6. In this section we
discuss weakenings of the $fsg$, relate the $fsg$ to forking, and in general try to obtain equivariant 
versions of results from section 2. 

\vspace{2mm}
\noindent
Let $G$ be a $\emptyset$-definable group in ${\bar M}$. $S_{G}(A)$ denotes the set of complete types over $A$ 
extending ``$x\in G$". Recall from \cite{NIP} and \cite{Shelah-groups} that if $T$ has $NIP$ then $G^{00}$ 
(the smallest type-definable bounded index subgroup of $G$) exists. Our notion of ``generic" will be from \cite{NIP}. However in \cite{metastable} ``generic" is used differently (there it is a 
translation invariant definable type). Also we distinguish it (as in \cite{NP}) from
``$f$-generic":
\begin{Definition} (i) A definable subset $X$ of $G$ is said to be left generic if finitely many left translates of $X$ by elements of $G$ cover $G$. 
\newline
(ii) $p(x)\in S_{G}(A)$ is left generic if every formula in $p$ is left generic. 
\newline
(iii) A definable subset $X$ of $G$ is said to be left $f$-generic if for each $g\in G$, $gX$ does not fork over
$\emptyset$. 
\newline
(iv) $p(x)\in S_{G}(A)$ is left $f$-generic if every formula in $p$ is left $f$-generic. 
\newline
(v) $G$ has the $fsg$ if there is a global type $p(x)\in S_{G}(\bar M)$ such that every left translate $gp$ of $p$ is almost finitely satisfiable over $\emptyset$, namely for every model $M_{0}$, every left translate of $p$ is finitely satisfiable in $M_{0}$.
\end{Definition}

Note that if $G$ is the additive group in $RCF$ (or just the underlying group in $DOAG$) then there are no complete generic types, but there are two $f$-generic types, at $+\infty$ and $-\infty$. (Both these types are definable over $\emptyset$ and invariant under the group operation). 

In \cite{NIP} we pointed out:
\begin{Fact} $G$ has the $fsg$ if and only if
\newline
(a) The class of left generic definable subsets of $G$ coincides with the class of right generic sets (we call them simply generic), and
\newline
(b) The class of non generic definable subsets of $G$ is a (proper) ideal, and
\newline
(c) If $X$ is a generic definable subset of $G$, then $X$ meets every model $M_{0}$.
\end{Fact}

The $fsg$ is a rather strong property for a definable group. It is a kind of abstract ``definable compactness".
The following was observed in \cite{NIP}. (Theorem 8.1 and Remark 3 of  Section 8.)
\begin{Remark} In $o$-minimal expansions of real closed fields, the definable groups with $fsg$ are precisely the 
definably compact groups. 
\end{Remark}

\vspace{2mm}
\noindent
Definable groups in 
algebraically closed valued fields have been studied by the first author (\cite{metastable}). From the 
analysis there, one can with additional work, obtain a characterization of the definable groups $G$ with $fsg$,
at least when $G$ is abelian: they are precisely the (abelian) groups $G$ such that there is a definable 
homomorphism $h:G\to A$ where $A$ is internal to the value group and is definably compact, and $Ker(h)$ is 
stably dominated. The $p$-adically closed field case is considered in \cite{OP}, where among other things it is
pointed out that definably compact groups defined over the standard model have $fsg$.

\begin{Corollary} If $G$ has the $fsg$, then there exists a (simultaneously left and right) global generic type $p\in S_{G}(\bar M)$ and moreover any such global generic type is also (left and right) $f$-generic. 
\end{Corollary}

\begin{Problem} Suppose $G$ has the $fsg$ and assume $NIP$ if one wants. Is it the case that every left $f$-generic type is generic?
\end{Problem}

In the case of definably compact groups in $o$-minimal expansions of real closed fields, Problem 5.5 has a 
positive answer. This follows from the results in \cite{Dolich} and \cite{Peterzil-Pillay} which give that: 
if $X$ is definable, closed and bounded then $X$ does not fork over $M_{0}$ if and only if $X$ meets $M_{0}$. 
So assume $X$ to be left $f$-generic and we want to prove that $X$ is generic. We may assume $X$ is closed. 
Let $M_{0}$ be any small model. So for any $g\in G$, $gX$ meets $M_{0}$. By compactness $X$ is right generic 
in $G$ (so generic in $G$).

\vspace{2mm}
\noindent
We now work without the $fsg$ assumption, improving a result from \cite{NIP}.
\begin{Proposition} (Assume $NIP$) Suppose that $G$ has a global left $f$-generic type $p$. Then 
\newline
(i) $Stab_{l}(p) = G^{00} = G^{000}$, and
\newline
(ii) $G$ is definably amenable, that is $G$ has a left invariant global Keisler measure.
\end{Proposition}
\pf  Recall that $Stab_{l}(p)$, the left stabilizer of $p$, is by definition the set of $g\in G$
such that $gp = p$. It is a subgroup of $G$ but on the face of it has no definability properties.
\newline 
(i) Let us fix a small model $M_{0}$ for convenience.
\newline 
{\em Claim I.}   If $tp(g_{1}/M_{0}) = tp(g_{2}/M_{0})$ then $g_{1}^{-1}g_{2}\in Stab_{l}(p)$. 
\newline
{\em Proof of Claim I.} By $f$-genericity of $p$ and the $NIP$ assumption, both $p$ and $g_{1}p$ are 
$Aut_{M_{0}}(\bar M)$-invariant. Let $f\in Aut_{M_{0}}(\bar M)$ such that $f(g_{1}) = g_{2}$.
So $g_{1}p = f(g_{1}p) = f(g_{1})f(p) = g_{2}p$. Hence $g_{1}^{-1}g_{2}\in Stab_{l}(p)$.

\vspace{2mm}
\noindent
{\em Claim II.} Suppose $h\in Stab_{l}(p)$. Then $h = g_{1}^{-1}g_{2}$ for some $g_{1},g_{2}$ with the
same type over $M_{0}$. 
\newline
{\em Proof of Claim II.} By definition $ha = b$ for some realizations of $p$ (in a model containing $\bar M$). So we can find realizations $a,b$ of $p|M_{0}$ in $\bar M$ such that $ha = b$. So $h = ba^{-1}$. Put $g_{1} = b^{-1}$ and $g_{2} = a^{-1}$ and note they have the same type over $M_{0}$.

\vspace{2mm}
\noindent
From Claims I and II we deduce that $Stab_{l}(p)$ is type-definable (over $M_{0}$) as the set
of $g_{1}^{-1}g_{2}$ for $g_{1},g_{2}\in \bar G$ having the same type over $M_{0}$. On the other hand by Claim I, the index of $Stab_{l}(p)$ in $G$ is bounded by the number of types over $M_{0}$. It follows that $Stab_{l}(p)$ has to be $G^{00}$, because $p$ determines a coset of $G^{00}$. The same reasoning yields that 
$Stab_{l}(p) = G^{000}$.

\vspace{2mm}
\noindent
(ii) This follows by the same proof as that of Proposition 6.2 of \cite{NIP}, using part (i) together with
Corollary 2.6: We can reduce to the case where $T$ is countable. Given $X$ a definable subset of $G$, let $M_{0}$ be a countable model over which $X$ is defined. Note that $p$ does not fork over $M_{0}$. The proof of Corollary 2.3 yields that not only is $p$ Borel definable over $M_{0}$, but (using countability of the language and of $M_{0}$) that any ``$\phi(x,y)$-definition of $p$" is a countable union of closed sets over $M_{0}$. In particular $Y = \{g\in G: X\in gp\}$ is a countable union of closed sets over $M_{0}$. By part (i) whether or not $g\in Y$ depends only on the coset $g/G^{00}$. Let $\pi:G\to G/G^{00}$ be the canonical surjection. Then $\pi(Y)$ is a Borel subset of $G/G^{00}$ and we define $\mu(X)$ to be ${\bf h}(\pi(Y))$ where ${\bf h}$ is the unique normalized Haar measure on the compact group $G/G^{00}$. $\mu$ is clearly additive and is also left-invariant:
Given $h\in G$, $\mu(hX) = {\bf h}(\{g/G^{00}:hX\in gp\})$. But $hX\in gp$ iff $g = hg'$ for some $g'$ such that $X\in g'p$. Hence $\{g/G^{00}:hX\in gp\} = h/G^{00} \cdot\{g'/G^{00}:X\in gp\}$. So by left-invariance
of ${\bf h}$, $\mu(hX) = {\bf h}(\{g'/G^{00}:X\in gp\}) = \mu(X)$. 

\begin{Remark} Proposition 5.6 also holds when $G$ is type-definable (with appropriate definitions). But for 
$G$ definable (over $\emptyset$) , we could also deduce 5.6 from 4.7 and 2.10(i), together with Construction C from
the introduction, in the easy form where the group action is regular. (Assume 
$NIP$). 
Let $T'$ the theory $T$ together with an additional sort $S$ and a $\emptyset$-definable regular action of $G$ 
on $S$. Let $x$ be a variable of sort $S$. Then the following are equivalent: (i) In $T$ there is a global 
left $f$-generic type, (ii) in $T'$, $x=x$ does not fork over $\emptyset$.
\end{Remark}

We now aim towards a kind of converse of 5.6. First:
\begin{Lemma} (Assume $NIP$) Let $G$ be a definable (or even type-definable) group, and let $\mu$ be a global 
left invariant Keisler
measure on $G$. Let $M_{0}$ be any small model over which $G$ is defined. Then there is a global left 
invariant Keisler measure $\mu'$ on $G$ such that $\mu'|M_0 = \mu|M_0$ and $\mu'$ is definable 
(hence $Aut({\bar M}/M)$-invariant for some small $M$).
\end{Lemma}
\pf 
The proof follows Keisler's construction in \cite{Keisler1} of ``smooth measures" assuming $NIP$. 
\newline 
{\em Step I.} We construct $M$ extending $M_{0}$ and extension $\mu_{1}$ of $\mu$ to a left $G(M)$-invariant measure over $M$ such that $\mu_{1}$ has a unique extension to a global left-invariant measure $\mu'$.
\newline
Suppose that already $\mu$ and $M_{0}$ do not work. So $\mu$ has two distinct extensions $\lambda, \beta$ to 
left invariant global measures. So there is a definable subset $X$ of $G$ and some positive real number
$r$ such that $|\lambda(X)-\beta(X)| > r$. Let $\mu'$ be the average of $\lambda$ and $\beta$. Note that $\mu'$
has the property
\newline 
(*) for any $M_{0}$-definable subset $Y$ of $G$, $\mu'(X\Delta Y) > r/4$.

Let $M_{1}$
be a small model containing $M_{0}$ such that $X$ is over $M_{1}$, and let $\mu_{1}$ be the restriction of 
$\mu'$ to $M_{1}$. Continue with $\mu_{1}$ in place of $\mu$. We claim that at some point we arrive at 
$\mu_{\alpha}$ over a model $M_{\alpha}$ which is left $G(M_{\alpha})$-invariant and has a unique global left 
invariant extension. Otherwise by (*) we find a formula $\phi(x,y)$, a positive real number $\epsilon$ 
and $(c_{i}:i<\omega)$, and a global measure $\mu''$ such that $\mu''(\phi(x,c_{i})\Delta\phi(x,c_{j})) > 
\epsilon$ for all $i\neq j$. This (together with $NIP$) is a contradiction. (See 3.14 of \cite{Keisler1} or 
3.2 of \cite{NIP}.) So Step I is accomplished.
\newline
{\em Step II.} $\mu'$ from Step I is definable (over $M$). 
\newline
This is proved (via Beth's theorem) exactly as in 
the proof of Lemma 2.6 (definability of minimal measures) in \cite{NIP}.

\vspace{5mm}
\noindent
The following consequence of Lemma 5.8 was also noticed by Krzysztof Krupinski.

\begin{Proposition} (Assume $NIP$) Suppose that the definable group $G$ is definably amenable (has a global left 
invariant Keisler measure). Then, after possibly adding constants for some small model over which $G$ is defined, 
there is a global left $f$-generic type of $G$.
\end{Proposition}
\pf  By Lemma 5.8, there is a left invariant global Keisler measure $\mu$ on $G$ which is moreover definable, 
over some small model $M$, hence does not fork over $M$. Let $p$ be a global type of $G$ such that $\mu(X) > 0$ for every $X\in p$. Then
for every $X\in p$ and $g\in G$, $\mu(gX) = \mu(X) > 0$, hence $gX$ does not fork over $M$.
So $p$ is a left $f$-generic of $G$ after adding names for the model $M$.

\vspace{5mm}
\noindent
So we obtain the following:
\begin{Corollary} (Assume $NIP$) $G$ is definably amenable if and only if $G$ has a global left $f$-generic type 
after possibly adding constants. Moreover under these circumstances $G^{000}$ exists and coincides with 
$G^{00}$.
\end{Corollary}

The corollary goes through for type-definable $G$ too. Shelah \cite{Shelah-groups} has recently proved that assuming
$NIP$, then for any abelian type-definable group $G$, $G^{000}$ exists. As abelian groups are amenable, the 
corollary above improves Shelah's theorem. On the other hand Gismatullin \cite{Gismatullin} has recently 
been able to drop the abelian hypothesis from Shelah's theorem. 

The following clarifies the question of the existence of left $f$-generic types and addition of constants.
In particular the proof of (ii) gives another example where forking does not equal dividing. 

\begin{Proposition} (Assume $NIP$) (i) Suppose $G$ is a definable group and has a global left $f$-generic type after 
naming elements of some model. Then $G$ has a global left $f$-generic type after naming elements of any model.
\newline
(ii) There is a $\emptyset$-definable group $G$ in a theory $T$ with $NIP$, such that $G$ is definably
amenable, but there is no global $f$-generic type of $G$.
\end{Proposition}
\pf  (i) We assume $G$ to be $\emptyset$-definable in $T$. We will use Remark 5.7. That is we work in the 
theory $T'$ which is $T$ together with a new sort $S$ and a regular right action of $G$ on $S$. 
${\bar M}$ is our monster model of $T$ and then ${\bar M'}$ denotes the corresponding (monster) model of $T'$,
that is ${\bar M}$ together with the new sort $S$. Let $x$ be a variable of sort $S$. By 5.7 it is enough to 
prove:
\newline
{\em Claim.} Suppose $x=x$ forks over the small elementary substructure $M_{0}$ of ${\bar M}$. Then $x=x$
forks over any small elementary substructure of ${\bar M}$ containing $M_{0}$.
\newline
{\em Proof of claim.} Assume for simplicity that $x=x$ implies $\phi(x,a)\vee\psi(x,b)$ where 
each of $\phi(x,a)$, $\psi(x,b)$ divides over $M_{0}$ (in ${\bar M'}$). Let $(a_{i}:i<\omega)$ and
$(b_{i}:i<\omega)$ be infinite $M_{0}$-indiscernible sequences in $tp(a/M_{0})$ and $tp(b/M_{0})$ respectively, such that $\{\phi(x,a_{i}):i<\omega\}$ and $\{\psi(x,b_{i}):i<\omega\}$ are each inconsistent.
Let $M_{1}$ be a small elementary substructure of ${\bar M}$ containing $M_{0}$. By an automorphism
argument we may assume that $tp(M_{1}/M_{0}\cup\{a_{i}:i<\omega\}\cup \{b_{i}:i<\omega\})$ is finitely 
satisfiable in $M_{0}$. But then clearly each of the sequences $(a_{i}:i<\omega)$ and $(b_{i}:i<\omega)$ is 
also $M_{1}$-indiscernible. So $x=x$ forks over $M_{1}$.
\newline
(ii) Let $K$ be a (saturated) real closed field. Let $G$ be $K\times K$, equipped with its addition, and with 
all
relations which are $\emptyset$-definable in $(K,+,\cdot)$ and invariant under the action of 
$SL(2,K)$.  Of course the theory of $G$ with this structure has $NIP$,
and $G$ is definably amenable as it is abelian.  
Suppose for a contradiction that $G$ had a global $f$-generic type $p(x)$. In particular $p(x)$ does not fork 
over $\emptyset$, so is invariant under automorphisms fixing $bdd(\emptyset)$ by section 2. As $SL(2,K)$ acts 
on $G$ by automorphisms
and $SL(2,K)^{00} = SL(2,K)$, it follows that $p$ is $SL(2,K)$-invariant. But then $p$ induces a $PSL(2,K)$
-invariant global type of the projective line over $K$,
which contradicts Remark 5.2 (iv) of \cite{NIP} and its proof.

\vspace{5mm}
\noindent
One might think that any global left generic type of any definable group is automatically left $f$-generic ($NIP$ or no $NIP$). But we need some assumptions to prove it.
\begin{Proposition} (Assume $NIP$). Let $G$ be a definable group which is definably amenable. Then any global left generic type of $G$ is also left $f$-generic.
\end{Proposition}
{\em Proof.} It is clearly sufficient to prove that any left generic definable set $X$ does not divide over 
$\emptyset$. And for this it is enough to find a small model $M_{0}$ such that any left generic definable set $X$ does not divide over $M_{0}$. (Given left generic $X_{a}$ defined over $a$, and an indiscernible sequence $(a_{i}:i<\omega)$ in $tp(a/\emptyset)$ we can find a sequence $(a_{i}':i<\omega)$, indiscernible over $M_{0}$ and with the same type over $\emptyset$ as $(a_{i}:i<\omega)$.)

The assumption that there is a global left generic type implies easily that the class $\cal I$ of non left generic definable sets is a proper ideal of the Boolean algebra of all definable subsets of $G$. By Proposition 6.3(i) of \cite{NIP}, there is a small model $M_{0}$ such that for every left generic definable set $X$ there is some $M_{0}$-definable subset $Y$ of $G$ such that the symmetric difference $X\Delta Y$ of $X$ and $Y$ is in $\cal I$, namely is non left generic. 

Now suppose that $X_{a}$ is a left generic definable subset of $G$ (defined over $a$). We want to show that $X_{a}$ does not divide over $M_{0}$. Let $Y$ be definable over $M_{0}$ such that the symmetric difference of $X_{a}$ and $Y$ is nongeneric. Replacing $X_{a}$ by $X_{a}\cap Y$ we may assume that $X_{a}\subseteq Y$. 
Let $(a_{i}:i<\omega)$ be an $M_{0}$-indiscernible sequence with $a_{0} = a$. Suppose by way of contradiction 
$\cap_{i=1,..n}X_{a_{i}}$ is inconsistent for some $n$. Then (as $Y$ is defined over $M_{0}$), the union of the $Y\setminus X_{a_{i}}$ for $i=1,..,n$ covers $Y$. But $Y$ is left-generic, and each $Y\setminus X_{a_{i}}$ is non left generic, contradicting $\cal I$ being an ideal. This proves the proposition.

\section{Generically stable groups}
In this section we introduce various strengthenings of the $fsg$ leading to the notion of a generically stable group (analogous to a generically stable type). Assuming $NIP$ these strengthenings will be equivalent. 

\vspace{2mm}
\noindent
We first give a natural definition of the $fsg$ for type-definable groups.
\begin{Definition} Let $G$ be a type-definable group. We say that $G$ has $fsg$ if there is some small model $M_{0}$ and global type $p(x)$ of $G$ such that $p$ is finitely satisfiable in $G(M_{0})$, namely for every formula $\phi(x)\in p$, $\phi$ is realized by some $g\in G(M_{0})$.
\end{Definition}

The basic results on definable groups with $fsg$ go through easily for type-definable groups $G$. Namely we define a relatively definable subset $X$ of $G$ to be left generic if finitely many left translates of $X$ by elements of $G$ cover $G$.
\begin{Lemma} Suppose $G$ is type-definable with $fsg$. Then
\newline
(i) the class of relatively definable left generic subsets of $G$ coincides with the class of relatively definable
right generic subsets of $G$ and forms a proper ideal of the class of relatively definable subsets of $G$.
\newline
(ii) There is a global generic type of $G$, namely a global type $p$ extending $x\in G$ such that any relatively
definable subset of $G$ which is in $p$ is generic.
\newline
(iii) There is a small model $M_{0}$ such that any global generic type of $G$ is finitely satisfiable in $G(M_{0})$, and such that moreover if $X$ is a relatively definable generic subset of $G$ then finitely many (left or right) translates of $X$ by elements of $G(M_{0})$ cover $G$. 
\newline
(iv) $G^{00}$ exists and is both the left and right stabilizer of any global generic type of $G$.
\end{Lemma}

\begin{Definition}  Let $G$ be a definable group.
\newline
(i) We say that $G$ has $fsg^{+}$ if $G$ has $fsg$ and some global generic type is definable.
\newline
(ii) We say that $G$ has $fsg^{++}$ if $G^{00}$ exists and has $fsg$ as a type-definable group in its own right.
\newline
(iii) We say that $G$ is generically stable if $G$ has $fsg$ and some global generic type is generically stable.
\end{Definition}

Note that (ii) seems to be a rather minor variation on $G$ having $fsg$ but in fact it is a substantial strengthening. Note in particular that if $G = G^{00}$ and $G$ has the $fsg$ then it has $fsg^{++}$. Definably compact groups in $o$-minimal structures satisfy none of the above properties, although they have $fsg$. We will see soon that (iii) implies (ii) implies (i) and that there are all equivalent assuming $NIP$.

\begin{Lemma} Suppose $G$ has $fsg^{++}$ (namely $G^{00}$ exists and has $fsg$).
Then
\newline
(i) $G$ has $fsg$,
\newline
(ii) Any global generic type of $G^{00}$ is a global generic type of $G$.
\newline 
(iii) There is a unique global generic type of $G$ which implies $``x\in G^{00}$, say $p_{0}$.
\newline
(iv) $p_{0}$ is definable and is also the unique global generic type of $G^{00}$.

\end{Lemma}
{\em Proof.} 
(i) Let $M_{0}$ be a small model witnessing that $G^{00}$ has $fsg$. We may assume that $M_{0}$ contains representatives of each coset of $G^{00}$ in $G$ (as there are boundedly many). Let $p$ be a global generic type of $G^{00}$. Let $g$ realize $p$. Then $g_{1}^{-1}g \in G^{00}$ for some $g_{1}\in G(M_{0})$. By assumption 
$g_{1}^{-1}gp$ is finitely satisfiable in $G(M_{0})$ (in fact in $G^{00}(M_{0})$. Hence so is $gp$. This shows that $G$ has $fsg$. 
\newline
(ii) We have shown in (i) that if $p$ is a global generic of $G^{00}$ then also $gp$ is finitely satisfiable in $G(M_{0})$ for any $g\in G$ which implies  that $p$ is a global generic type of $G$. (Alternatively if $X$ is a definable subset of $G$ and $X\cap G^{00}$ is generic in $G^{00}$ then finitely many translates of $X$ cover $G^{00}$ whereby finitely many translates of this set of finitely many translates covers $G$.)
\newline
(iii) Suppose for a contradiction that $p$ and $q$ are distinct global generic types of $G$ each of which implies 
that $``x\in G^{00}$. Let $X$ be a definable set in $p$ such that $Y = G\setminus X$ is in $q$. Let $X_{0} = X\cap G^{00}$ and $Y_{0} = Y\cap G^{00}$. So $X_{0}$ and $Y_{0}$ are relatively definable subsets of $G^{00}$ which partition $G^{00}$. By Lemma 6.2 we may assume that $X_{0}$ is generic in $G^{00}$. So finitely many translates of $X_{0}$ by elements of $G^{00}$ cover $G^{00}$. In particular (as $q(x)$ implies $x\in G^{00}$), there is $g\in G^{00}$ such that $gX\in q$. But by Corollary 4.3 of \cite{NIP} the symmetric difference of $X$ and $gX$ is nongeneric (in $G$), so as $q$ is a generic type of $G$, we see that $X$ is also in $q$, a contradiction. 
\newline
(iv) Let $p_{0}$ be the unique type from (iii). By part (ii), $p_{0}$ is also the unique global generic type of $G^{00}$. Let $G^{00}$ be the intersection of the directed family $(Y_{i})_{i\in I}$ of definable subsets of $G$.
For a given $L$-formula $\phi(x,y)$ and $c\in {\bar M}$, $\phi(x,c)\in p_{0}$ iff for some $i\in I$ finitely many translates of (the set defined by) $\phi(x,c)$ by elements of $G^{00}(M_{0})$ cover $Y_{i}$. The same is true for $\neg\phi(x,y)$. By compactness the set of $c$ such that $\phi(x,c)\in p_{0}$ is definable. So $p_{0}$ is definable.

\begin{Proposition} Let $G$ be a definable group. Then
\newline
(i) If $G$ has $fsg^{++}$ then $G$ has $fsg^{+}$.
\newline
(ii) If $G$ is generically stable then $G$ has $fsg^{++}$.
\newline 
(iii) Assume $NIP$. If $G$ has $fsg^{+}$ then $G$ is generically stable. (Hence assuming $NIP$ the three properties $fsg^{+}$, $fsg^{++}$ and generic stability are equivalent).
\end{Proposition}
{\em Proof.} (i) follows from the previous lemma.
\newline
(ii) By translating we obtain a generically stable global generic type $p(x)$ of $G$ such that
$p(x)\models ``x\in G^{00}"$. Fix a small model $M_{0}$ such that $p$ is definable over and finitely satisfiable in $M_{0}$ (and of course $G$ and $G^{00}$ are defined over $M_{0}$). Let $M$ be a small $|M_{0}|^{+}$-saturated extension of $M_{0}$. Let $(a_{i}:i<\omega)$ be a Morley sequence in $p$ over $M_{0}$. Then $a_{i}\in G^{00}(M)$ for all $i$. Suppose $\phi(x,c)\in p$. The generic stability of $p$ implies that $\phi(x,c)$ is satisfied by some $a_{i}$, so $\phi(x,c)$ is satisfied in $G^{00}(M)$. If $g\in G^{00}$ then $gp = p$ hence any left translate of $p$ by an element of $G^{00}$ is finitely satisfiable in $G^{00}(M)$, so $G^{00}$ has $fsg$ and $G$ has $fsg^{++}$.
\newline
(iii) is immediate because assuming $NIP$ any global type which is both definable over and finitely satisfiable in some small model, is generically stable.

\vspace{5mm}
\noindent
We now develop some consequences of the ``weakest" of the new properties, $fsg^{+}$.
\begin{Proposition} Suppose that the definable group $G$ has $fsg^{+}$. Then
\newline
(i) $G^{00} = G^{0}$.
\newline 
(ii) There is a unique global generic type $p_{0}$ of $G$ extending $``x\in G^{00}"$ (and $p_{0}$ is definable).
\newline 
(iii) Every global generic type of $G$ is definable and the set of global generics is in one-one 
correspondence with $G/G^{00}$. 
\newline 
(iv) For any $L$-formula $\phi(x,y)$ there is $N<\omega$ such for any $c$, $\phi(x,c)$ is generic in $G$ if and
only if at most $N$ left (right) translates of $\phi(x,c)$ cover $G$.
\end{Proposition}
{\em Proof.} (i) Fix a definable global generic type $p$. By 4.3 from \cite{NIP} we know that $Stab(p) = G^{00}$.
For a fixed formula $\phi(x,y)$ let (just for now) $p|\phi$ be the set of formulas of the form $\phi(g^{-1}x,c)$ which are in $p$ for $g\in G$ and $c\in {\bar M}$. Then $Stab_{\phi}(p)$ denotes 
$\{g\in G: g(p)|\phi = p|\phi\}$ and is a definable subgroup of $G$. As $Stab(p) = Stab_{\phi}(p)$ for $\phi\in L$ it follows that each $Stab_{\phi}(p)$ has finite index and so $G^{00}$ is the intersection of a family of definable subgroups of finite index, so equals $G^{0}$. 
\newline
(ii) By translating the given definable generic $p$ we can find a definable global generic type which extends 
$``x\in G^{00}"$ and we call it $p_{0}$. Note that if $q(x)$ is a global generic extending $``x\in G^{00}"$ then so is $q^{-1}$. Hence it suffices to prove:
\newline
{\em Claim.} If $q$ is any global generic type extending $``x\in G^{00}"$ then $q^{-1} = p_{0}$.
\newline
{\em Proof of claim.} Let $M_{0}$ be any model. So $q$ is finitely satisfiable in $M_{0}$ and $p_{0}$ is definable over $M_{0}$. Let $a$ realize $p|M_{0}$  and $b$ realize $q|(M_{0},a)$. As $a\in G^{00}$ and $G^{00} = Stab(q)$ it follows that $c= ab$ also realizes $q|(M_{0},a)$. By Lemma 3.4, $a$ realizes $p|(M_{0},c)$. Again as $c^{-1}\in G^{00}$ we see that $c^{-1}a = b^{-1}$ also realises $p|(M_{0},c)$. In particular $tp(b^{-1}/M_{0}) =
p_{0}|M_{0}$. As $M_{0}$ was arbitrary it follows that $q^{-1} = p_{0}$. The claim is proved as well as part (ii).
\newline 
(iii) Now any global generic type of $G$ is a translate of a global generic which implies $``x\in G^{00}"$. By (ii) every global generic type of $G$ is a translate of the definable type $p_{0}$, hence is also definable.
\newline 
(iv) If $\phi(x,y)\in L$ and $c\in {\bar M}$, then $\phi(x,c)$ is generic in $G$ if for some $g\in G$, $\phi(g^{-1},c)\in p_{0}$. So by definability of $p_{0}$ the set of $c$ such that $\phi(x,c)$ is generic is definable, which is enough.

\vspace{5mm}
\noindent
Let us remark that there are groups with $G^{00} = G^{0}$, with the $fsg$ but without $fsg^{+}$ (even in a $NIP$ environment). The easiest example is simply the additive group $R^{+}$ of the valuation ring $R$ in a saturated $p$-adically closed field. Here $(R^{+})^{00}$ is the intersection of the finite index definable subgroups $v(x)\geq n$ for $n\in {\Z}^{+}$, but there are many generic types extending $(R^{+})^{00}$.

\vspace{2mm}
\noindent
\begin{Remark} Suppose the definable group $G$ is generically stable. Let $X$ be a definable subset of $G$. Then $X$ is generic if and only if every left (right) translate of $X$ does not divide (or even fork) over $\emptyset$. 
\end{Remark}
{\em Proof.} We know left to right by just the $fsg$. Now suppose that $X$ is nongeneric, and defined over some model $M_{0}$. Let $p$ be some generically stable global generic type of $G$, and let $(a_{i}:i<\omega)$ be a Morley sequence of $p$ over $M_{0}$. Then for any $g\in G$, as $gX\notin p$, we see that $a_{i}\notin gX$ for some $i<\omega$. It follows that $\{a_{i}^{-1}X:i<\omega\}$ is inconsistent. As $(a_{i}^{-1}:i<\omega)$ is indiscernible over $M_{0}$ we see that $X$ divides over $M_{0}$ so over $\emptyset$. 

\vspace{5mm}
\noindent
Section 7 of this paper is devoted to establishing the uniqueness of invariant (under the group action) measures for definable groups with $fsg$ (assuming $NIP$). In the case of generically stable groups this can be seen quickly and we do it now.

\begin{Lemma} Suppose the definable group $G$ is generically stable, and $\mu$ is a left invariant global Keisler measure on $G$. Then $\mu$ is generic, in the sense that for any definable subset $X$ of $G$, if $\mu(X) > 0$ then $X$ is generic.
\end{Lemma}
{\em Proof.} Suppose that $X$ is a nongeneric definable subset of $G$, defined over a model $M_{0}$. The argument in the proof of Remark 6.7 gives an $M_{0}$-indiscernible sequence $(a_{i}:i<\omega)$ such that $\{a_{i}X:i<\omega\}$ is inconsistent. But $\mu(a_{i}X) = \mu(X)$ for all $i$. By Lemma 2.8 of \cite{NIP}, $\mu(X) = 0$. 

\begin{Lemma} Suppose $G$ is a definable group with $fsg^{+}$. Then any definable subset of $G$ is a Boolean combination of translates (cosets) of definable subgroups of $G$ of finite index and nongeneric definable sets.
\end{Lemma}
{\em Sketch of proof.} Fix a formula $\phi(x,y)\in L$. By a $\phi$-formula we mean a Boolean combination of
formulas $\phi(gx,c)$ ($g\in G, c\in {\bar M}$). By a global $\phi$-type of $G$ we mean a maximal consistent 
set of $\phi$-formulas over ${\bar M}$ (containing in addition $x\in G$). A global $\phi$-type will be called 
generic iff it contains only generic formulas (iff it extends to a global generic type of $G$). As in the 
proof of 6.6, there is a definable subgroup $G^{0}_{\phi}$ of finite index which is the (say left) 
stabilizer of every global generic $\phi$-type. It follows that each of the finitely many cosets of 
$G^{0}_{\phi}$ contains a unique global generic $\phi$-type. Thus every $\phi$-set (set defined by a $\phi$
-formula) is a Boolean combination of cosets of $G^{0}_{\phi}$ and nongeneric definable sets.

\begin{Corollary} Suppose the definable group $G$ is generically stable. Then there is a unique left 
invariant (under ther group action) Keisler measure on $G$ which is also the unique right invariant Keisler 
measure.
\end{Corollary} 
{\em Proof.} If $\mu$ is a left invariant Keisler measure on $G$, then $\mu$ is determined on definable 
subgroups of finite index and their translates, and by Lemma 6.8 is $0$ on all nongeneric sets. By Lemma 6.9 
there is only one possibility for $\mu$. Lemma 6.9 is clearly still true replacing definable subgroups of 
finite index by normal definable subgroups of finite index. Hence we see that $\mu$ is also the unique right 
invariant Keisler measure on $G$. 

\vspace{5mm}
\noindent
Note that a special case is: if $G$ is a connected stable group, and $p$ is its unique global generic type, 
then not only is $p$ the unique left (right) invariant global type of $G$, but it is also the unique left 
(right) invariant global Keisler measure on $G$.

\section{Uniqueness of invariant measures for groups with $fsg$}
Given a definable group $G$ with $fsg$, and assuming $NIP$ we have from \cite{NIP} that there is a left 
invariant global Keisler measure $\mu$ on $G$ (namely $G$ is definably amenable). In fact we constructed such 
$\mu$ which is generic (the definable sets of positive measure are precisely the generics). Clearly $\mu^{-1}$ 
is a right invariant generic global Keisler measure. Our uniqueness theorem (Theorem 7.7 below) generalizes the uniqueness 
of invariant types for connected stable groups, as well in a sense generalizing the uniqueness of Haar measure
for compact groups. (If one is willing to pass to continuous model theory, a compact group is like something 
finite, hence stable.) The main point is to prove a Fubini (or symmetry) theorem for suitable measures. A. Berarducci also 
raised the uniqueness issue (in the $o$-minimal context) in \cite{Berarducci} and pointed out the relevance
of Fubini. Our methods go through variants of Grothendieck groups.

\vspace{5mm}
\noindent
We begin with a lemma on homomorphisms from vector spaces to the reals. First recall some notation: if $A$ is 
a $\Q$-vector space, then a subset $P$ of $A$ is said to be a {\em cone} if $P$ is closed under addition and 
under scalar multiplication by positive elements of $\Q$. 
\begin{Lemma}  Let $A$ be a $\Q$-vector space, $B$ a $\Q$-subspace of $A$, $P$ a cone in $A$, and $\mu_{B}$ a 
nonzero homomorphism from $B$ to $\R$ such that $\mu_{B}(B\cap P) \geq 0$. Assume that for any $a\in A$ there 
is $b\in B$ such that $b-a\in P$.
Let $c\in A$ and let $e = inf\{\mu_{B}(y):y\in B, y-c\in P\}$.
\newline
Then (i) if $e = -\infty$ then there is no homomorphism $\mu:A\to \R$ extending $\mu_{B}$ with $\mu(P)\geq 0$.
\newline
(ii) If $e\in \R$, then there is a homomorphism $\mu:A\to\R$ extending $\mu_{B}$ such that $\mu(P)\geq 0$ and
$\mu(c) = e$.
\end{Lemma}
{\em Proof.} First there is no harm in assuming that $0\in P$.
\newline
(i) It is enough to prove the stronger statement that if $\mu:A\to \R$ is any homomorphism
extending $\mu_{B}$, and $\mu(P)\geq 0$ then $\mu(c) \leq e$. But this is immediate: if $e'>e$ then there
is $y\in B$ such that $y-c\in P$ and $\mu_{B}(y)\leq e'$, but then $\mu_{B}(y)-\mu(c)\geq 0$ so
$\mu (c) \leq \mu_{B}(y) \leq e'$. 
\newline
(ii) Let us assume $e\in\R$, and by the stronger statement we have just proved all we need to do is to find a homomorphism $\mu:A\to\R$ extending $\mu_{B}$ such that $\mu(P)\geq 0$ and 
$\mu(e)\geq c$. 

As $\mu_{B}$ is nonzero, let $b_{1}\in B$ be such that $\mu(b_{1}) > 0$, and there is no harm
in assuming $\mu_{B}(b_{1}) = 1$. Replacing $c$ by $mb_{1} + c$ for sufficiently large $m$ we
may also assume that $e>0$.

Let $P_{B} = \{b\in B: \mu_{B}(b)\geq 0\}$. Note that if $\mu:A\to\R$ is a homomorphism such that $\mu(b_{1}) = 1$ and $\mu$ is nonnegative on $P_{B}$ then $\mu$ extends $\mu_{B}$. So it suffices to find some homomorphism $\mu$ from $A$ to $\R$
such that $\mu(b_{1}) = 1$, $\mu$ is nonnegative on $P + P_{B} = \{a+b:a\in P, b\in P_{B}\}$, and
$\mu(c) \geq e$. (Because then $\mu$ is nonnegative on $P_{B}$, and also on $P$.)
\newline
Let $P_{c} = \{\alpha c - \beta b_{1}: \alpha, \beta \in \Q$, $\alpha, \beta \geq 0$, $\beta < e\alpha\}\cup\{0\}$.
\newline
{\em Claim.} $-b_{1} \notin P' = P+P_{B} + P_{c}$.
\newline
{\em Proof of claim.}  Suppose not.
\newline
{\em Case (i).} $-b_{1} = a + b + \alpha c - \beta b_{1}$ for $a\in P$,
$b\in P_{B}$, and $\alpha,\beta$ as in the 
definition of $P_{c}$. 
\newline
Multiplying by $\alpha^{-1}$ and letting $a' = \alpha^{-1}a$, $b' = \alpha^{-1}b$,
and $\gamma = \alpha^{-1}\beta$
we have
\newline
$-\alpha^{-1} b_{1} = a' + b' + c - \gamma b_{1}$, whence
\newline
$\gamma b_{1} - (b' + \alpha^{-1}b_{1}) - c  = a' \in P$.
\newline
But then $e \leq \mu_{B}(\gamma b_{1} - (b' + \alpha^{-1} b_{1})) \leq \gamma < e$ a contradiction.
\newline
{\em Case (ii).} $-b_{1} = a + b$ for $a\in P$, $b\in P_{B}$.
\newline
Then $-a = b_{1} + b$, so $-a \in B$ and $\mu_{B}(-a) > 0$. But then also $a\in B$, so $a\in B\cap P$, so by assumption
$\mu_{B}(a)\geq 0$ whereby $\mu_{B}(-a) < 0$, again a contradiction.
 
\vspace{2mm}
\noindent
By the claim we can extend the cone $P'$ to a maximal cone $P''$ not containing $-b_{1}$, and note that
$b_{1}\in P''$, and also each of $P+P_{B}$, $P_{c}$ are contained in $P''$. Note also that $P''$
defines a linear preorder on $A$, namely for each $d\in A$, $d\in P''$ or $-d\in P''$. (If $d\notin P''$ then 
$-b_{1} \in P'' + \Q ^{+} d$, so $-d\in P'' + \Q ^{+} b_{1}$, so $-d \in P''$.) Our assumptions on $A$, 
together with the definition of $P_{B}$ imply that for any $a\in A$ there is $n>0$ such that 
$nb_{1} - a \in P+P_{B}$ so also in $P''$. Thus $\Q ^{+}b_{1}$ is cofinal in $A$. Now it is clear that the 
homomorphism from $\Q b_{1}$ to $\R$ which sends $b_{1}$ to $1$ extends to an preorder preserving homomorphism 
$\mu$ from $A$ to $\R$ (which sends $A_{0} = \{x\in A: b_{1} - nx \in P''$ for all $n\in \Z\}$ to $0$). Then 
$\mu(P+P_{B}) \geq 0$, as $P+P_{B}\subseteq P''$. But also, as $P_{c}\subseteq P''$ we have that $c-\gamma 
b_{1} \in P''$ for all positive $\gamma < e$. Hence $\mu(c) \geq e$, and we have found the required $\mu$.

\vspace {5mm}
\noindent
We now consider a certain variant of the Grothendieck (semi)-group introduced in section 5 of \cite{NIP}
We will also work at the more general level of definable group actions rather than just definable 
groups. So we will fix a definable group action $(G,X)$ and a small model $M_{0}$ over which $(G,X)$ is 
definable. Using the notation analogous to that in \cite{NIP}, we will take the relevant semigroup 
$K_{semi}(X)$ to be the 
collection of nonnegative cycles $\sum_{i}k_{i}X_{i}$ in $X$ up to piecewise translation by members of 
$G(M_{0})$. $K_{0}(X)$ will be the corresponding Grothendieck group. (Recall that we define 
$x_{1}\sim_{0}x_{2}$ for $x_{1},x_{2}\in K_{semi}(X)$ if there is $y\in K_{semi}(X)$ such that $x_{1}+y = 
x_{2}+y$ in $K_{semi}$. Then $K_{0}$ is the collection of $\sim_{0}$-classes together with formal inverses.)
When we apply Lemma 7.1 to this situation, the $\Q$-vector space $A$ will be 
the tensor product of $\Q$ with $K_{0}(X)$. Define $P_{0}$ to be the image of $K_{semi}$ in $K_{0}$ (under the 
canonical map) and then $P$ will be $\{\alpha x: \alpha\in \Q^{+}, x\in P_{0}\}$, a cone in $A$. Define 
$B_{0}$ to be the image in $K_{0}(X)$ of the ``$M_{0}$-definable" members of $K_{semi}(X)$, and
then $B$ will be the tensor product of $\Q$ with $B_{0}$.

 We now give a somewhat more concrete representation of the objects defined in the previous paragraph. 
It will be convenient both notationally and conceptually to introduce a certain category  
$\cal D$ in place of the semigroup
of nonnegative cycles on $X$ (namely before identification via piecewise translations). 
First we think of $\Z$ as living in our monster model ${\bar M}$ as the directed union of the finite sets (or 
sorts)
$\{-m,...,m\}$. We also have the group structure on $\Z$ definable in ${\bar M}$, piecewise. Let $\tilde X$ be 
$X\times \Z$. A definable subset of $\tilde X$ is by definition a definable subset of $X\times\{-m,..,m\}$ for
some $m$. Likewise for $\tilde G = G\times\Z$, which is now a ``definable" group. The natural action of 
$\tilde G$ on $\tilde X$ is also ``definable", namely for each $m,n$ we have the map
$(G\times \{-m,..,m\})\times (X\times\{-n,..,-n\})\to (X\times\{-m-n,..,m+n\})$ is definable. The objects of 
the category are definable subsets of $\tilde X$. If $Y, Y'$ are in $\cal D$ then a morphism $f$ from $Y$ to 
$Y'$ is an injective map $f$ from $Y$ into $Y'$ such that there is a partition $Y = Y_{1}\cup ...\cup Y_{n}$ of $Y$
and elements $g_{i}\in {\tilde G}(M_{0}) = G(M_{0})\times \Z$ such that for each $i$ and $y\in Y_{i}$,
$f(y) = g_{i}y$. Equivalently (by compactness) $f$ is a definable embedding of $Y$ in $Y'$ such that for each 
$y\in Y$ there is $g\in {\tilde G}(M_{0})$ such that $f(y) = gy$. So $Y$ and $Y'$ will be isomorphic in $\cal 
D$ if there is $f$ is above which is a bijection between $Y$ and $Y'$. Then as a set $K_{semi}(X)$ is the set 
of isomorphism classes of members of $\cal D$. The addition on $K_{semi}(X)$ is induced by the following
(non well defined) addition on $\cal D$: if $Y, Y'$ are definable subsets of $X\times \{-m,..,m\}$ then $Y + Y'
= Y \cup (2m+1)Y'$. (Likewise for $nY$.) Connecting with earlier notation, $B_{0}$ is the image in $K_{0}(X)$ of the set of 
isomorphism 
classes of $M_{0}$-definable elements of $\cal D$ (and $B$ the tensor product of $B_{0}$ with $\Q$).
For $Y$ a definable subset of $\tilde X$, $[Y]$ denotes its image in $A = K_{0}(X)\otimes \Q$.

We will be interested in Keisler measures $\mu$ on $X$ which are $G(M_0)$-invariant. Note 
that any such $\mu$ extends uniquely to (and is determined by) a $\tilde G_{0}$-invariant finitely additive 
measure on the definable 
subsets of $\tilde X$. 

\begin{Lemma} The global $G(M_{0})$-invariant Keisler measures on $X$ correspond to the homomorphisms
$h:A\to \R$ such that $h(P)\geq 0 $ and $h[X] = 1$. Namely if $\mu$ is a global $G(M_{0})$-invariant Keisler 
measure
on $X$, then there is a unique homomorphism $h:A\to \R$ such that $h([Y]) = \mu(Y)$ for any definable subset 
$Y$ of $X$. Moreover this $h$ satisfies $h([X]) = 1$ and $h(P)\geq 0$. Conversely if $h:A\to \R$ is such that 
$h([X]) = 1$ and $h(P)\geq 0$ then defining $\mu(Y) = h([Y])$ for any definable subset $Y$ of $X$ gives a 
global $G(M_{0})$-invariant Keisler measure $\mu$ on $X$.
\end{Lemma}
{\em Proof.} Clear.

\vspace{2mm}
\noindent
\begin{Corollary} Suppose that $\mu$ is a global $G(M_{0})$-invariant Keisler measure on $X$ which
is moreover the unique global $G(M_{0})$-invariant Keisler measure on $X$ extending
$\mu|M_{0}$. Then for any definable subset $D$ of $X$ (or $\tilde X$) with $\mu(D) = \beta$ and $\epsilon > 0$, there are 
$n\in \N$,
$M_{0}$-definable $E_{0}, E_{1} \in {\cal D}$ and $D'\in {\cal D}$ such that
\newline
(i) $\mu(E_{1}) - \mu(E_{0}) < n(\beta + \epsilon)$ and
\newline
(ii) there is a morphism $f$ (in ${\cal D}$) from $E_{0} + nD + D'$ to $E_{1} + D'$. 
\end{Corollary}
{\em Proof.} With earlier notation let $h:A\to\R$ be the homomorphism corresponding to $\mu$ as given by 
Lemma 7.2.  So by Lemma 7.1, $h$ is the unique homomorphism such that $h(P)\geq 0$, and which extends
$h|B$. Let $D$ be a definable subset of $X$ (or $\tilde X$) and put $c = [D]$, and let
$\mu(D) = h([D]) = \inf\{h(y):y\in B, y-c \in P\}$. Fix $\epsilon > 0$. So for some $y\in B$ we have 
$h(y) < \mu(D) + \epsilon$ and $y-c\in P$.  Now for large enough positive $n$ we have
$n(y-c) = [D'']$ for some definable subset $D''$ (of $\tilde X$), and $ny = [E_{1}] - [E_{0}]$ with
$E_{0}, E_{1}$ $M_{0}$-definable subsets of $\tilde X$. So we have:
\newline
(a) $\mu(E_{1}) - \mu(E_{0})$ = $h([E_{1}] - [E_{0}]) = nh(y) < n(\mu(D) + \epsilon)$, which gives (i).
\newline
We also have:
\newline
(b) $[E_{1}] - [E_{0}] - n[D] = [D'']$. 

We then obtain (after possibly multiplying everything by some $m>0$) some definable $D'$
such that in $\cal D$, $E_{0} + D'' + nD + D'$ is isomorphic to $E_{1} + D'$. Ignoring $D''$ this gives a 
morphism $f$ in $\cal D$ from $E_{0} + nD + D'$ to $E_{1} + D'$ which is (ii).

\vspace{5mm}
\noindent
Let us give a more explicit statement.
\begin{Remark} The conclusion of Corollary 7.3 can be restated as:
\newline
For any definable subset $D$ of $X$ with $\mu(D) = \beta$ and for any $\epsilon > 0$
there are $n,m,m',m'' \in \N$, $M_{0}$-definable sets $E_{0}, E_{1}$ in $\cal D$, and some
$D'\in \cal D$,  as well as
$D_{i}\in \cal D$ and $g_{i}\in \tilde G_{0}$ for $i=1,,.m$, $D_{j}'\in \cal D$ and $g_{j} \in \tilde G_{0}$
for $j = 1,..,m'$,
and $D_{k}''\in \cal D$ and $g_{k''}\in\tilde G_{0}$ for $k=1,..,m''$ such that
\newline
(i) $\mu(E_{1}) - \mu(E_{0}) < n(\beta + \epsilon)$,
\newline
(ii) $nD = \cup_{i=1,..,m}D_{i}$, $D' = \cup_{j=1,,.m'}D_{j}'$, 
$E_{0} = \cup_{k=1,..,m''}D_{k}''$, the $D_{i}, D_{j}', D_{k}''$ are pairwise disjoint,
and $D'$ is disjoint from both $D$ and $E_{1}$,
\newline
(iii) the sets $g_{i}D_{i}$, $g_{j'}D_{j}'$ and $g_{k}''D_{k}''$ are pairwise disjoint 
subsets of $E_{1}\cup D'$.
\end{Remark}

We can now prove the sought after symmetry (or Fubini) theorem. Recall that if
$\mu_{x}, \lambda_{y}$ are global Keisler measures on definable sets $X,Y$ respectively,
and $\mu$ is definable (or even Borel definable), then we can form the global Keisler measure
$\mu\otimes\lambda$ on $X\times Y$: for $D$ a definable subset of $X\times Y$,
$(\mu\otimes\lambda) (D) = \int \mu(D(y))d\lambda$, where $D(y) = \{x:(x,y)\in D\}$.
We may also write $(\mu\otimes\lambda)(D)$ as $\int_{D}d\mu d\lambda$.

\begin{Proposition} Suppose $(G,X)$ is a definable group action, and $Y$
a definable set, all defined over $M_{0}$. Suppose $\mu$ is a global Keisler
measure on $X$ which is definable and satisfies the hypothesis of Corollary 7.3.
Suppose $\lambda$ is a global definable Keisler measure on $Y$. Then $\mu_{x}\otimes\lambda_{y}
= \lambda_{y}\otimes\mu_{x}$.
\end{Proposition}

Before beginning the proof, let us define an action of $G$ on
$X\times Y$ by $g(x,y) = (gx,y)$ and we claim that both measures $\mu\otimes\lambda$ and 
$\lambda\otimes\mu$ are $G(M_{0})$-invariant for this action: Let $D$ be a definable subset of
$X\times Y$, and $g\in G(M_{0})$. So for any $y\in Y$, $(gD)(y) = g(D(y))$, so $\mu((gD)(y))
= \mu(D(y))$, and thus $(\mu\otimes\lambda((gD) = (\mu\otimes\lambda)(D)$. On the other hand,
let $f(x) = \lambda(D(x))$. Now clearly $D(g^{-1}(x)) = g(D(x))$, so $f(g^{-1}x) = \lambda(gD(x))$. 
As $\mu$ is $G(M_{0})$-invariant, $\int f(x) d\mu  = \int f(g^{-1}x) d\mu$, so 
$(\lambda\otimes \mu)(D) = (\lambda\otimes\mu((gD)$. All this of course extends to the actions of 
$\tilde G$  on $\tilde X$ and ${\tilde X\times Y}$ and the relevant measures on ${\cal D}_{X}$ and ${\cal D}_{X\times Y}$.

\vspace{2mm}
\noindent
{\em Proof of Proposition 7.5.} Let $D$ be a definable subset of $X\times Y$. We have to show that
$(\mu\otimes\lambda)(D) = (\lambda\otimes\mu)(D)$. By considering also the complement of $D$, it suffices to 
prove that for any $\epsilon > 0$, $(\lambda\otimes\mu)(D) \leq (\mu\otimes\lambda)(D) + \epsilon$.

Fix $\epsilon > 0$. By Corollary 7.3, for each $y\in Y$ we find $n_{y},m_{y},m'_{y}m, m''_{y}$, $(E_{0})_{y}$, $(E_{1})_{y}$, $D'_{y}$ etc.
such that (i), (ii) and (iii) of Remark 7.4 hold for $D_{y} = \{x:(x,y)\in D\}$. By compactness we may partition $Y$ into definable sets
$Y_{\nu}$, such that the $n_{y},m_{y}, m'_{y}, m''_{y}$, 
$(E_{0})_{y}, (E_{1})_{y}$, 
$(g_{i})_{y}, (g_{j}')_{y},
(g_{k}'')_{y}$ are constant on each $Y_{\nu}$. Focus attention on one $Y_{\nu}$. Let 
$D_{\nu} = D\cap (X\times Y_{\nu}$. Let $D'_{\nu}\subseteq X\times Y_{\nu}$ be such that for $y\in Y_{\nu}$, 
$D'_{y} = \{x:(x,y)\in D'_{nu}\}$. So clearly we have:
\newline
{\em Claim I.} There is a morphism in ${\cal D}_{X\times Y_{\nu}}$, from 
$(E_{0}\times Y_{\nu}) + nD_{\nu} + D'_{\nu}$ into $(E_{1}\times Y_{\nu}) + D'_{\nu}$. 

\vspace{2mm}
\noindent
As $\mu(E_{1}) - \mu(E_{0}) < n(\mu(D_{y}) + \epsilon)$ for all $y\in Y_{\nu}$, we obtain, on integrating
with respect to $\lambda$ over $Y_{\nu}$ that
\newline
{\em Claim II.} $(\mu(E_{1}) - \mu(E_{0})\lambda(Y_{\nu}) \leq n(\mu\otimes\lambda)(D_{vu}) + 
\epsilon\lambda(Y_{\nu})$.

\vspace{2mm}
\noindent
But the left hand side in Claim II coincides with $(\lambda\otimes\mu)(E_{1}\times Y_{\nu})-
(\lambda\otimes\mu)(E_{0}\times Y_{\nu)}$. So denoting $\lambda\otimes\mu$ by $r$ we rewrite Claim II as
\newline
{\em Claim IV.} $r(E_{1}\times Y_{\nu}) - r(E_{0}\times Y_{\nu}) \leq n(\mu\otimes\lambda)(D_{vu}) + 
\epsilon\lambda(Y_{\nu})$.

\vspace{2mm}
\noindent
We have already noted that 
$r$ on $X\times Y$ is 
$G(M_{0})$-invariant, so applying
$r$ to Claim I and using the disjointness there, we obtain:
\newline
{\em Claim IV.} $r(E_{0}\times Y_{\nu}) + nr(D_{\nu}) + r(D'_{\nu}) \leq r(E_{1}\times Y_{\nu}) +
r(D'_{\nu})$. 
\newline
From Claims II and IV we obtain
\newline
{\em Claim V.} $r(D_{\nu}) \leq (\mu\otimes\lambda)(D_{\nu}) + \epsilon\lambda(Y_{\nu})$.

\vspace{2mm}
\noindent
Summing over $\nu$ in Claim V, we obtain
\newline
$\lambda\otimes\mu)(D) \leq (\mu\otimes\lambda)(D) + \epsilon$ which is what we wanted. The proposition is 
proved.

\begin{Lemma} (Assume $NIP$) Suppose $G$ is a definable group, defined over a small model $M_{0}$, and $\mu$ is a 
Keisler
measure on $G$ over $M_{0}$, which is left $G(M_{0})$-invariant. Then
\newline
(i) There is a global left $G(M_{0})$-invariant Keisler measure $\mu'$ extending $\mu$, and a small model
$M$ containing $M_{0}$ such that $\mu'$ is the unique left $G(M_{0})$-invariant global Keisler measure
extending $\mu'|M$. Again any such $\mu'$ is definable.
\newline
(ii)  Suppose in addition that $G$ has $fsg$. Then any $\mu'$ as in (i) is left invariant, namely left $G({\bar M})$-invariant. 
\end{Lemma}
\pf
(i) is proved in exactly the same way as 5.8.
\newline
(ii) By the definability of $\mu'$, $H =Stab(\mu')$ is a type-definable subgroup of $G$. We want to show 
that 
it is all of $G$. If not, it is clear that the complement of $H$ contains a generic definable subset $X$ of 
$G$. By the $fsg$, $X\cap G(M_{0})$ is nonempty. But then $X$ contains an element of $H$ (as $\mu'$ is 
$G(M_{0})$-invariant). Contradiction.

\vspace{5mm}
\noindent
Combined with the material of the previous section we can now obtain our main result.
\begin{Theorem} (Assume $NIP$). Suppose $G$ is a definable group with $fsg$. Then $G$ has a unique left invariant 
global Keisler measure, which is also the unique right invariant global Keisler measure of $G$. This measure 
is both definable and generic.
\end{Theorem}
{\em Proof.} We already know from \cite{NIP} that $G$ has some left invariant 
global Keisler
measure. Let $\mu, \lambda$ be left invariant global Keisler measures. We will first show that $\mu = 
\lambda^{-1}$.
Let $D$ be any definable subset of $G$. 
We want to prove that $\mu(D) = \lambda(D^{-1})$. Let $M_{0}$ be a small model over which both $G$ and $D$ are
defined. Let $\mu''$ be a global Keisler measure satisfying (i) of 7.6 for some small $M$ containing $M_{0}$.
Namely $\mu''$ extends $\mu|M_{0}$ and is the unique left $G(M_0)$-invariant Keisler measure extending $\mu''|M$. By (ii) of 7.6, 
$\mu''$ is left $G(\bar M)$-invariant, in particular left $G(M)$-invariant, so is also the unique left $G(M)$
-invariant extension of $\mu''|M$. So renaming $M$ as $M_{0}$, $\mu''$ satisfies the hypothesis of 7.3. By 
7.6 (ii) $\mu''$ is definable, and as already mentioned $\mu''$ is left invariant. So as we are currently just 
interested in $\mu(D)$ we may assume that $\mu = \mu''$. Likewise we may assume that $\lambda$ is definable. 
By Proposition 7.5, $\mu\otimes\lambda = \lambda\otimes\mu$ (on $G\times G$). Let $Z = \{(x,y)\in G\times 
G:x\in yD\}$ which is 
equal to $\{(x,y)\in G\times : y\in xD^{-1}\}$. 
So $(\mu\otimes\lambda)(Z) = \int \mu(yD) d\lambda = \mu(D)$ as $\mu$ is left invariant.
Likewise $(\lambda\otimes\mu)(Z) = \int \lambda(xD^{-1})d\mu = \lambda(D^{-1})$. 

So we have succeeded in proving that $\lambda = \mu^{-1}$. This applies in particular when $\lambda = \mu$
yielding that $\mu = \mu^{-1}$ (so $\mu$ is also right invariant). Hence $\mu = \lambda$ and is also the 
unique right invariant global Keisler measure on $G$. 

\vspace{5mm}
\noindent
Finally we point out the extension of Proposition 7.7 to homogeneous spaces.
\begin{Proposition} (Assume $NIP$) Suppose $(G,X)$ is a definable homogeneous space, and $G$ has $NIP$.
Then there is a unique $G$-invariant global Keisler measure on $X$.
\end{Proposition}
{\em Proof.} Let $\mu$ be the unique invariant global; Keisler measure on $G$ given by the previous 
proposition. So it is definable and satisfies the hypothesis of 7.3. Note that we obtain from $\mu$ a $G$
-invariant global Keisler measure $\lambda$ on $X$: given $x\in X$, the map taking $g\in G$ to $gx$ gives a 
($G$)-invariant surjection $\pi_{x}:G\to X$. For $Y$ a definable subset of $X$, define $\lambda_{x}(Y) = 
\mu(\pi_{x}^{-1}(Y))$. Then $\lambda_{x}$ is $G$-invariant and does not depend on the choice of $x$ so we call 
it $\lambda$. Note that for any $x\in X$, and definable subset $Y$ of $X$, $\lambda(Y) = \mu(\{g\in G:gx\in 
Y\})$.
Clearly $\lambda$  is also definable and satisfies the hypothesis of 7.3. 

Let $\lambda'$ be another $G$-invariant global Keisler measure on $X$. We want to show $\lambda' = \lambda$. As in 
the proof of 7.7 we may assume $\lambda'$ is definable. By 7.5, $\mu\otimes \lambda' = \lambda'\otimes \mu$ 
on $G\times X$. Let $Y$ be a definable subset of $X$ and let $Z = \{g,x):gx\in Y\}$.
Then $(\mu\otimes\lambda')(Z) = \int_{x\in X}(\mu\{g\in G:gx\in Y\})d\lambda'$ = 
$\int_{x\in X}\lambda(Y)d\lambda'$ = $\lambda(Y)$. 

And $(\lambda'\otimes\mu)(Z) = \int_{g\in G}\lambda'(\{x\in X:gx\in Y\})d\mu$ = $\int_{g\in G}\lambda'(Y)d\mu$ = $\lambda'(Y)$. So $\lambda = \lambda'$.

\vspace{5mm}
\noindent
There is an obvious formulation of $X$ having $fsg$ where $X$ is a definable homogeneous space (for a
definable group $G$):
there is a global type $p$ of $X$ such that for every $g\in G$,
$gp$ is finitely satisfiable in $M_{0}$ for every model $M_{0}$. 

If $G$ has $fsg$ then clearly $X$ does too.

\begin{Question} (Assume $NIP$) (i) Does Proposition 7.8 holds if only $X$ has $fsg$?
\newline
(ii) Are there examples of transitive group actions $(G,X)$ such that $X$ has $fsg$, $G$ acts faithfully, and $G$ does not have $fsg$?
\end{Question}

\begin{Remark}  Let us say that a definable group $G$ has the {\em weak fsg} if there is some small model $M_{0}$ over which $G$ is defined such that $G$ has no proper type-definable subgroup containing $G(M_{0})$. The following results follow from the proofs above.
\newline
(i) (Assume $NIP$.) If $G$ is definably amenable with weak $fsg$ then $G$ has a unique left invariant Keisler measure which is also its unique right invariant Keisler measure, and is definable.
\newline
(ii) (Assume $NIP$.) Suppose the definable group is definably amenable and has no proper nontrivial type definable subgroups. Then $G$ has a global left invariant type, which is definable, and is moreover the unique left invariant Keisler measure and the unique right invariant Keisler measure on $G$.
\end{Remark} 
{\em Proof.} (i) is given by Proposition 7.5, and the proofs of Lemma 7.6 and Theorem 7.7.
\newline
The assumptions of (ii) imply that $G$ has weak $fsg$, and so let $\mu$ be the unique (left, right) invariant measure on $G$ given by (i). As $\mu$ is definable, $\mu$ does not fork over some small $M_{0}$. Let $p$ be some random global type for $\mu$, namely every formula in $p$ has positive $\mu$-measure. Then clearly every left translate of $p$ does not fork over $M_{0}$. By Proposition 5.6, $Stab_{l}(p) = G^{00}$. But our assumptions imply that $G^{00} = G$, hence $p$ is left invariant. By part (i) $\mu = p$. 

\vspace{5mm}
\noindent
Note that statement (ii) in  Remark 7.10 above mentions only types, but the only proof we have of it involves measures. 
In fact, using 7.5 one can conclude that, assuming $NIP$, any group with weak $fsg$ has $fsg$ (and so the groups in Remark 7.10 (ii) are generically stable groups in the sense of section 6) but the proof uses the theory of generically stable measures which will be treated in a future paper.

\vspace{5mm}
\noindent
Finally in this section we speculate on group free analogues of the $fsg$ property. Here is a possible definition, in which the definable group $G$ is replaced by a complete type $p(x)\in S(A)$. For simplicity we make a blanket assumption of $NIP$.
\begin{Definition} We say that $p(x)\in S(A)$ has $fsg$ if $p$ has a global extension $p'(x)$ such that
for every $(|L|+|A|)^{+}$-saturated model $M$ containing $A$, $p'$ is finitely satisfiable in $p(M)$, that is 
for each formula $\phi(x,b)\in p'$ there is $a\in M$ which realizes $p$ and satisfies $\phi(x,b)$).
\end{Definition}

\begin{Remark} (i) In Definition 7.11 the global extension $p'$ of $p$ is necessarily a nonforking extension of $p$.
\newline
(ii) Suppose $A = bdd(A)$ and $p(x)\in S(A)$  has $fsg$. Then $p$ has a unique global nonforking extension $p'$, and $p'$ is generically stable. 
\newline
(iii) Suppose the global type $p'$ is generically stable and $p'$ does not fork over $A$. Let $p = p'|A$. Then $p$ has $fsg$. 
\end{Remark}
{\em Proof.}  (i) If $(b_{i}:i<\omega)$ is an $A$-indiscernible sequence (in ${\bar M}$) then we can find
a $(|L|+|A|)^{+}$-saturated model $M$ containing $A$ such that $(b_{i}:i < \omega)$ remains indiscernible over $M$. 
If $\phi(x,b_{0})\in p'$ then $\models \phi(a,b_{0})$ for some $a\in M$ hence $\phi(a,b_{i)}$ for all $i<\omega$. So $p'$ does not fork over $A$.
\newline
(ii) By part (i) and 2.11 (and the assumption that $A = bdd(A)$), $p'$ is $A$-invariant. We will show that some Morley sequence in $p'$ over $A$ is totally indiscernible (and apply 3.2). Let $M$ be a small $(|L|+|A|)^{+}$-saturated  model containing $A$. Let $(a_{i}:i<\omega)$ be a Morley sequence in $p'$ over $M$  (so also one over $A$). It is enough to
prove that for any $n$, $tp(a_{0},..,a_{n-1},a_{n}/A) = tp(a_{n},a_{0},..,a_{n-1}/A)$. So fix $n\geq 1$. Suppose that
$\phi(x_{0},..,x_{n-1},x)$ is over $A$ and $\models \phi(a_{0},..,a_{n-1},a_{n})$. Then $\phi(a_{0},..,a_{n-1},x)\in p'$, so by hypothesis is realized by some $c\in M$ which realizes $p$. But then $(c,a_{0},a_{1},..,a_{n-1})$ also begins a Morley sequence in $p'$ over $A$, so $tp(c,a_{0},..,a_{n-1}/A) = tp(a_{0},..,a_{n-1},a_{n}/A)$, and we see that $\models \phi(a_{n},a_{0},..,a_{n-1})$. This suffices. 
\newline
(iii) follows from the material in section 3 and is left to the reader.

\vspace{5mm}
\noindent
In the context of Remark 7.12(iii) we say that the fsg type $p(x)\in S(A)$ comes from a generically stable type.
So (by 7.12(ii)) group actions, namely the action of the compact Lascar group $Aut(bdd(A)/A)$ on the set of extensions of $p$
over $bdd(A)$, will enter the picture whenever an $fsg$ type
$p(x)\in S(A)$ does {\em not} come from a generically stable type.

Another source of $fsg$ types is through transitive group actions and construction $C$ from the
section 1. Namely suppose that $(G,X)$ is a $\emptyset$-definable
group action in $T$ (with $NIP$), and we form $T_{X}$ with additional sort $X'$. Then in $T_{X}$ there is a 
unique (so isolated) $1$-type $p(x)$ over $\emptyset$ in $X'$. If $X$ has $fsg$ in $T$ then it is
not hard to see that $p$ has $fsg$ in $T_{X}$.


In fact, types $p\in S(A)$ with $fsg$ can be characterized as types of the form $\mu|A$ where $\mu$ is a global $A$-invariant generically stable measure (in the sense of Remark 4.4). This was proved by P. Simon and will again appear in a joint work with the authors.

\section{The Compact Domination Conjecture}

We will prove
\begin{Theorem} Assume $\bar M$ to be an $o$-minimal expansion of
a real closed field. Let $G$ be a definably connected definably compact commutative group definable
in $\bar M$. Then $G$ is compactly dominated. That is, let $\pi:G\to G/G^{00}$ be
the canonical surjective homomorphism. Then  for any definable subset $X$ of $G$, 
$Y_{X} = \{b\in G/G^{00}:\pi^{-1}(b)\cap X\neq\emptyset$ and $\pi^{-1}(b)\cap G\setminus X \neq \emptyset\}$ has Haar measure $0$.
\end{Theorem}

In fact the proof (and subsequent elaborations) will yield a bit more: in the structure ${\bar M}^{*}$ obtained by expanding $\bar M$
by a predicate for $G^{00}$, $G/G^{00}$ will be ``semi-$o$-minimal" with dimension that of $G$, and moreover
for any definable subset $X$ of $G$ the set $Y_{X}$ above (which is now definable in 
${\bar M}^{*}$) has dimension $< dim(G/G^{00})$. So in a sense $G$ is $o$-minimally dominated (by a
definable $o$-minimal compact Lie group). This of course suggests many problems and issues for future work. 
Also, in the paper \cite{centralextensions} joint with Peterzil, the full compact domination conjecture (i.e. for not necessarily commutative definably compact groups) is deduced from Theorem 8.1 and results in \cite{NIP}, using a structure theorem for definably compact groups in $o$-minimal expansions of real closed fields. Of course Theorem 8.1 builds on and uses the work and contributions of a number of people, 
including Berarducci, Dolich, Edmundo, Otero, and Peterzil.

Until we say otherwise  ${\bar M}$ denotes a (saturated) $o$-minimal expansion of a real closed field $K$, 
$G$ is a definably compact definably connected definable group in $\bar M$ of $o$-minimal
dimension $n$, and $\pi$ is the canonical surjective homomorphism from $G$ to $G/G^{00}$. Without loss of 
generality $G$ is defined over $\emptyset$. $M_{0}$ will denote a fixed small model. We will make use in a few places of the fact that $G^{00}$ can be defined by a countable collection of formulas. (This is by the DCC result in \cite{BOPP}.)
The overall proof has several steps and ``new" 
ingredients, including
a beautiful result of Otero and Peterzil (\cite{O-P}). Some of the lemmas go through or can be formulated at various
greater levels of generality but we tend to concentrate on the case at hand. By Lemma 10.5 of \cite{NIP} in order to prove Theorem 9.1 it suffices to show that if $X$ is a definable subset of $G$ of dimension $<n$
then the Haar measure of $\pi(X)\subseteq G/G^{00}$ is $0$. We will aim towards this.

\vspace{2mm}
\noindent
The first step is to show that $G^{00}$ is definable in some weakly $o$-minimal expansion of $\bar M$.
Recall that weak $o$-minimality means that every definable set of elements (rather than tuples) of the 
universe is a finite union of convex sets (with respect to the underlying ordering). Let 
${\bar M}^{*}$ be the expansion of $\bar M$ obtained by adjoining a predicate for the trace on $\bar M$ of
any set definable with parameters in some elementary extension of $\bar M$. In other words, 
for each $L$-formula $\phi(x,y)$ ($x,y$ tuples) and complete type $q(x)$ over $\bar M$, adjoin a predicate
for $\{b:\phi(x,b)\in q(x)\}$. By the results of Baisalov and Poizat \cite{B-P}, or alternatively of Shelah 
\cite{Shelah783}, 
$Th({\bar M}^{*})$ is weakly $o$-minimal, with quantifier elimination. The weak $o$-minimality of a theory 
means that every model is weakly $o$-minimal, namely that for some (any) model there is a bound on the number of convex 
components of definable subsets of $1$-space in definable families.

\begin{Lemma} $G^{00}$ is definable in ${\bar M}^{*}$.
\end{Lemma}
\pf  Let $p$  be a global generic type of $G$ (which exists by 5.2 and 5.3). By Proposition 5.6 (and 5.12) we have 
\newline
(i) $Stab(p) = G^{00}$,
\newline
and moreover as mentioned in the proof of 5.6, 
\newline
(ii) $G^{00}= \{g_{1}g_{2}^{-1}:tp(g_{1}/M_{0}) = tp(g_{2}/M_{0})\}$.

\vspace{2mm}
\noindent
Given a formula $\psi(x)$ (over some parameters) defining a subset $X$ of $G$, let
$Stab_{\psi}(p) = \{g\in G:$ for all $h\in G$, $hX\in p$ iff $ghX\in p\}$. Then clearly $Stab_{\psi}(p)$ is a subgroup of $G$. Moreover  $Stab(p) = \cap_{\psi}Stab_{\psi}(p)$ and
\newline
(iii) $Stab_{\psi}(p)$ is definable in ${\bar M}^{*}$.

\vspace{2mm}
\noindent
We will show that $G^{00}$ is the intersection of finitely many $Stab_{\psi}(p)$, even with the $\psi$ over $M_{0}$.

\vspace{2mm}
\noindent
For $\psi(x)$ over $M_{0}$, let $S_{\psi}(p)$ be the smallest type-definable over $M_{0}$ (in $\bar M$) set 
containing
$Stab_{\psi}(p)$. (Note that as $p$ is $Aut({\bar M}/M_{0})$-invariant, so is $Stab_{\psi}(p)$.) We don't necessarily know that $S_{\psi}(p)$ is a subgroup of $G$, but:
\newline
{\em Claim (iv)}. ($\psi$ over $M_{0}$.) $S_{\psi}(p)\cdot G^{00}$ is a (type)-definable subgroup of $G$ (clearly containing
$Stab_{\psi}(p)$). 
\newline
{\em Proof.} Clearly $\pi(S_{\psi}(p))\cdot G^{00}$ is the closure of $\pi(Stab_{\psi}(p))$ in $G/G^{00}$, and 
hence is a closed subgroup of $G/G^{00}$. Its preimage in $G$ is thus a subgroup of $G$ and coincides
with $S_{\psi}(p)\cdot G^{00}$. 

\vspace{2mm}
\noindent
{\em Claim (v).}  $\cap_{\psi}(S_{\psi}(p)\cdot G^{00}) = G^{00}$ (where the $\psi$'s in the left hand side are taken only over $M_{0}$). 
\newline
{\em Proof.} By (i), the left hand side contains the right hand side. Let us show the reverse inclusion. Fix $g\in G\setminus G^{00}$. So for every 
$h\in G^{00}$, $gh^{-1}\notin G^{00}$. So by (ii), for every $h\in G^{00}$ and $a$ realizing $p|M_{0}$,
$gh^{-1}a$ does not realize $p|M_{0}$. By compactness there is a  formula $\psi(x)\in p|M_{0}$,
such that for any $h\in G^{00}$ and $a$ satisfying $\psi(x)$, $gh^{-1}a$ does not satisfy $\psi(x)$.
However, for any $g_{1}\in Stab_{\psi}(p)$ clearly there is $a$ such that both $a$ and $g_{1}a$ satisfy
$\psi$. Thus for every $g_{1}\in S_{\psi}(p)$ the same is true. So we have shown that (for our given choice of
$g\notin G^{00}$), for all $h\in G^{00}$, $gh^{-1}\notin S_{\psi}(p)$. Namely 
$g\notin S_{\psi}(p)\cdot G^{00}$. This proves Claim (v). 

\vspace{2mm}
\noindent
We can now complete the proof of the Lemma. We know that $G$ has the DCC on type-definable subgroups of bounded index (see \cite
{BOPP}). Hence by Claim (v), 
$G^{00}$ is the intersection of finitely many of the $S_{\psi}(p)\cdot G^{00}$. Thus 
(as $G^{00} \subseteq Stab_{\psi}(p) \subseteq S_{\psi}(p)\cdot G^{00}$) $G^{00}$ is the intersection 
of finitely many of the $Stab_{\psi}(p)$. By (iii), $G^{00}$ is definable in ${\bar M}^{*}$.

\begin{Remark} (i) The proof above uses only the existence of an $f$-generic type (rather than a generic type), and hence, by 5.6 and 5.10, Lemma 8.2 goes through assuming only that $T$ has $NIP$, $G$ is a definable definably amenable group, and 
$G/G^{00}$ is a compact Lie group.
\newline
(ii) Otero and Peterzil pointed out to us a more direct proof of 8.2, but which makes more use of the fact that
$G$ has $fsg$: Choose a global generic type $p$ of $G$, and by 4.3 of \cite{NIP} and the fact that $G/G^{00}$ has $DCC$, $G^{00}$ is of the form $Stab_{\cal I}(X_{1})\cap...\cap Stab_{\cal I}(X_{n})$, where the $X_{i}$ are definable sets in $p$, and ${\cal I}$ denotes the ideal of nongeneric definable subsets of $G$. (Here 
$Stab_{\cal I}(X)$ is the set of $g\in G$ such that the symmetric difference of $X$ and $gX$ is nongeneric, and so is type-definable in ${\bar M}$.) But genericity (so also nongenericity) is ``definable" in ${\bar M}^{*}$, because a definable subset $X$ of $G$ is generic iff for some $g\in G$, $gX\in p$. Hence each $Stab_{\cal I}(X_{i})$ is definable in ${\cal M}^{*}$ which suffices. 
\end{Remark}

\vspace{5mm}
\noindent
The next step is given to us by Lemma 4.3 of \cite{O-P}.
\begin{Lemma} Let $I$ be the interval $[0,1)$ in $\bar M$ (or $K$).
Then there are one-one definable continuous functions $\gamma_{1},..,\gamma_{n}$ from $I$ to $G$
such that $G = \gamma_{1}(I) + ... + \gamma_{n}(I)$ (where $\gamma_{j}(I)$  denotes the image of $I$ under 
$\gamma_{j}$). 
\end{Lemma}

\noindent
In fact the $\gamma_{i}$ are generators of the $o$-minimal fundamental group of $G$ and the proof in 
\cite{O-P} has a 
(co)homological character. 

\vspace{2mm}
\noindent
For $j = 1,..,n$, let $I_{j} = \gamma_{j}(I)$. As $\pi$ is a surjective homomorphism we have:
\begin{Corollary} $G/G^{00} = \pi(I_{1}) +.... + \pi(I_{n})$. 
\end{Corollary}

\noindent
We are now in the following interesting situation. $G/G^{00}$ as well as its subsets $\pi(I_{j})$ are
compact spaces under the logic topology, namely in their capacity as bounded hyperdefinable sets in the saturated model ${\bar M}$, but they are 
also, by Lemma 8.2, definable (or rather interpretable) sets in the (non-saturated) weakly $o$-minimal
structure ${\bar M}^{*}$. We will show in the next step that each $\pi(I_{j})$ is a ``definable $o$-minimal set" in ${\bar M}^{*}$ which piecewise is an interval in $\R$.

Let us first be precise about what we mean by a definable set in an ambient structure being $o$-minimal. Let $N$ be a structure, and $X$ a definable (or interpretable) set in $N$. When we say ``definable" we mean definable in $N$ with parameters from $N$. Suppose that $<$ is a definable linear ordering on $X$. We call $X$ $o$-minimal in $N$ with respect to $<$, if any definable subset $Y$ of $X$ is a finite union of intervals
$(a,b)$ (where possibly $a = -\infty$, $b = +\infty$) and singletons.  We will call $X$ strongly $o$-minimal if in addition there is a finite bound on the number of intervals and points in definable families. 

 In the next proposition we will mention the existence of definable Skolem functions on a definable set. So let us again be precise about the meaning. Again let $N$ be a structure in language $L$ and $X$ a set definable (or interpretable) in $N$, for now definable without parameters. We say that $X$ has definable Skolem functions in $N$ (or the formula defining $X$ has definable Skolem functions in $Th(N)$), if for any $L$-formula $\phi(x,y)$ where 
$x$ is a variable ranging over $X$, and $y$ is an arbitrary tuple of variables, there is a partial $\emptyset$-definable function $f_{\phi}(y)$ such that
 \newline
$N \models (\forall y)((\exists x\in X \phi(x,y)) \rightarrow (\phi(f_{\phi}(y),y) \wedge f_{\phi}(y)\in X))$.

\vspace{2mm}
\noindent
We can also speak of $X$ having Skolem functions in $N$ over some set $A$ of parameters from $N$.
The reader should note that if $X$ has definable Skolem functions in $N^{eq}$ (over some parameters), then $X$ has elimination of imaginaries (over some parameters), namely whenever $Z\subseteq X^{n}$ and equivalence relation $E$ on $Z$ are definable in $N$ then there is a definable bijection of $Z/E$ with some definable $W\subseteq X^{m}$.

\vspace{2mm}
\noindent
We now return to the main narrative, with notation following Lemma 8.4.
The following Proposition is fundamental. The appendix is devoted to the proof of part (iii).

\begin{Proposition}  (Work in ${\bar M}^{*}$.) For each $j$, $\pi(I_{j})$ is a finite disjoint union of definable sets
$X_{1},..,X_{r}$ and points $c_{1},..,c_{s}$, such that each $X_{i}$ is equipped with a definable total ordering $<_{i}$, such that (for each $i = 1,..,r$),
\newline
(i) $(X_{i},<_{i})$ is (abstractly) isomorphic to $\R$ with the usual ordering,
\newline
(ii) $X_{i}$ is strongly $o$-minimal (with respect to $<_{i}$),
\newline
(iii) $X_{i}$ has definable Skolem functions after possibly naming some parameters from ${\bar M}^{*}$.
\end{Proposition}
\pf  Let us fix $j$. Note that $I_{j}$ has a canonical linear ordering coming from the map $\gamma$.
Let $E$ be the equivalence relation ``$x-y\in G^{00}$" on $I_{j}$, i.e. coming from $\pi$. As $E$ is definable in 
${\bar M}^{*}$ (which has weakly $o$-minimal theory), each $E$-class is a union of at most $k$ convex sets 
for some $k$. Let $Y$ be the set of elements of $I_{j}$ which are in the ``first" convex subset of their $E$
-class.
Then $Y$ is definable,  so a finite union of convex definable sets, and $\pi(Y) = \pi(I_{j})$.
Let $E_{Y}$ denote the restriction of $E$ to $Y$. 
 Write $Y$ minimally as a finite disjoint union of 
(definable) convex sets $Y_{1},..,Y_{t}$. So each $E_{Y}$-class is convex and contained in a unique $Y_{i}$.
Let $X_{i} = \pi(Y_{i})$ and note that $<$ induces a linear ordering $<_{i}$ on $X_{i}$ (of course definable
in ${\bar M}^{*}$). Let us restrict our attention to some $Y_{i}$ such that $X_{i}$ is infinite. 
\newline
{\em Claim (I).}  $<_{i}$ is dense on $X_{i}$
 (i.e. if $a,b\in X_{i}$ and $a<_{i}b$ then there is $c\in X_{i}$ with $a <_{i} c <_{i} b$). 
\newline
{\em Proof.} We may work inside an interval $I' = (a',b')$ of $Y_{i}$ (i.e. with $\pi(a')\leq_{i} a$
and $\pi(b') \geq_{i} b$). So $\pi^{-1}(a) \cap I'$ and $\pi^{-1}(b)\cap I'$ are convex sets which are type-definable in ${\bar M}$ and disjoint. Moreover the first has no greatest element and the second has no smallest element (as every coset of $G^{00}$ is open in $G$). Hence by compactness (in ${\bar M}$) there is
$c'\in I'$  such that $\pi^{-1}(a)\cap I' < c' < \pi^{-1}(b)\cap I'$. Let $c = \pi(c')$.

\vspace{2mm}
\noindent
From Claim (I) we may assume that $(X_{i},<_{i})$ has no first or last element (by removing them if they exist).

\vspace{2mm}
\noindent
{\em Claim (II).} $(X_{i},<_{i})$ is complete, namely every bounded above (below) subset has a supremum
(infimum).
\newline
{\em Proof.} Again we may work in $\pi(I')\subset X_{i}$ for some interval $I' = (a',b')$ in $Y_{i}$. We consider 
$\pi(I')$ with the logic topology (equivalently as a closed subset of $G/G^{00}$). 
Let 
$(A,B)$ be a Dedekind cut in $(\pi(I'), <_{i})$. For $a\in A$, $A_{a} = \{x\in \pi(I'): a\leq_{i} x\}$ is closed in $\pi(I')$, as it is clearly
the image of a type-definable (in ${\bar M}$) subset of $I$. Likewise for $b\in B$,
$B_{b} = \{x\in \pi(I'): x\leq_{i} b\}$ is closed. Hence by compactness of the space $\pi(I')$, there is a point in the intersection of all the $A_{a}$'s and $B_{b}$'s. This suffices.

\vspace{2mm}
\noindent
{\em Claim (III).}  $(X_{i},<_{i})$ is separable, namely has a countable dense subset with respect to the ordering $<_{i}$.
\newline
{\em Proof.}  We know that $X$ with the logic topology is second countable (has a countable basis), because $E$ is given by countably many formulas. (See  Remark 1.6 of \cite{BOPP} for example.) We will show that every $<_{i}$-interval $(a,b)$ in $X_{i}$ contains an open subset of $X$. So as $X$ has a countable basis, $(X_{i},<_{i})$ will have a countable dense subset.

We work one $i$ at a time. The reader should convince himself/herself that there is no harm in assuming  that $Y = Y_{1}$.  So $X = X_{1}$. We relabel $<_{1}$ as $<_{X}$. We fix an $<_{X}$ interval $(a,b)$ in $X$ and we want to find a subinterval which is open in $X$ (with the logic topology). Let $Z$ be the union of all the second convex components of $E$-classes in $I$. Let $E_{Z}$ be the restriction of $E$ to $Z$. So $Z$ is definable in ${\bar M}^{*}$ and is a disjoint union of finitely many definable convex subsets $Z_{1},..,Z_{m}$ of $I$ such that each $E_{Z}$ class is convex and contained in a unique $Z_{j}$. Consider $\pi(Z)\cap (a,b)$. If it is finite, then after passing to a smaller interval, we may assume that $\pi(Z)\cap[a,b]$ is empty. This means that for $c\in [a,b]$, $\pi^{-1}(c) \subset Y$.
Let $a_{0}\in \pi^{-1}(a)$, and $b_{0}\in \pi^{-1}(b)$. So $(a,b)$ is the set of $E$-classes which are contained in the interval $(a_{0},b_{0})$ in $I$. Thus $(a,b)$ is open in $X$ and we are finished. 

So we may assume (by $o$-minimality and passing to a smaller interval) that $(a,b)$ is contained in $\pi(Z)$. 
The ordering $<$ on $Z$ induces a definable ordering $<_{2}$ on $\pi(Z) \subseteq X$. After ignoring finitely many points, we know (by the description of linear orderings definable on an $o$-minimal structure), that we can break up $(a,b)$ into finitely many $<_{X}$intervals, on each of which $<_{2}$ agrees with $<_{X}$ or $>_{X}$. It follows that we can find a subinterval 
$(c,d)$ of $(a,b)$ and some $j = 1,.., m$ such that $\pi^{-1}([c,d])\cap Z$ is contained in $Z_{j}$, and moreover for some $c_{0}',d_{0}'$ in $Z_{j}$ (preimages of $c,d$), $[c,d]$ is the image under $\pi$ of the closed $<$-interval between $c_{0}'$ and $d_{0}'$, and moreover $<_{2}$ on $[c,d]$ agrees with $<_{X}$ or $>_{X}$. 

Now if $k=2$ we are finished: Let $c_{0},d_{0}\in Y$ be preimages of $c,d$. Then $(c,d)$ is precisely the set of $E$-classes contained in $(c_{0},d_{0})\cup (c_{0}',d_{0}')$  (or in $(c_{0},d_{0}) \cup (d_{0}',c_{0}')$ if $<_{2}$ on $[c,d]$ is $>_{X}$). 

If $k > 2$ we continue, replacing $(a,b)$ by $(c,d)$, and considering now $W$ the union of the third convex components of $E$-classes. Our choice of $(c,d)$ means that passing to smaller subintervals does not disturb the compatibility with the second convex components of $E$-classes.

This finishes the proof of Claim (III). 

\vspace{2mm}
\noindent
It follows from Claims (I), (II) and (III) that $(X_{i},<_{i})$ is isomorphic as an ordered set
to $\R$ with its usual ordering. So we have proved (i).

\vspace{2mm}
\noindent 
(ii) follows quickly. For if $Z$ is a definable (in ${\bar M}^{*}$) subset of $(X_{i})$, then by
weak $o$-minimality of ${\bar M}^{*}$
$\pi^{-1}(Z)\cap Y_{i}$ is a finite union of convex sets. By completeness of $(X_{i},<_{i})$, $Z$
is a finite union of intervals and points. Weak $o$-minimality of $Th({\bar M}^{*})$ implies that there
 is a bound on the number of connected components of $Z$ as it varies in a definable family.
 
 \vspace{2mm}
 \noindent
 (iii). See the appendix.

\vspace{10mm}
\noindent

Let $X_{1}^{j},...,X_{r_{j}}^{j}$ be the sets obtained for $\pi(I_{j})$ in Proposition 8.6. By Proposition 8.6, each $X_{k}^{j}$ is strongly $o$-minimal in ${\bar M}^{*}$. The reader should convince himself or herself and the general machinery of $o$-minimality 
(dimension, cell decomposition, etc.) applies to Cartesian products of the $X_{k}^{i}$ and definable (in ${\bar M}^{*}$) subsets $W$ of such Cartesian products. In fact we will call such a definable set $Z$, a semi-$o$-minimal definable set. By (iii) of Proposition 8.6, and the discussion preceding the statement of 8.6, we have ``elimination of imaginaries" for such definable sets: namely if $Z$ is a subset of a Cartesian product of the $X_{k}^{j}$ and $E$ an equivalence relation on $Z$, both definable in ${\bar M}^{*}$ then $Z/E$ is in definable bijection  with some definable $W$ which is a subset of a Cartesian product of the $X_{k}^{j}$'s. We will apply this to $G/G^{00}$ considered as a group definable (or interpretable) in ${\bar M}^{*}$. In fact let us write $J$ for this group, so as to distinguish it from $G/G^{00}$ as a bounded hyperdefinable group in ${\bar M}$.
By Corollary 8.5 $J$ is in 
the definable closure (uniformly) of the $X_{k}^{j}$'s. Thus there is a definable subset $Z$ of some Cartesian product of the $X_{k}^{j}$, and a definable equivalence relation $E$ on $Z$ such that $J$ is definably isomorphic to $Z/E$. Hence $J$ is definably isomorphic to some definable $W$ which is a subset of a Cartesian product of the $X_{k}^{j}$'s. But $J$ also has a definable group structure, hence is a semi-$o$-minimal group.
Note that $J$ is definably connected (in the sense of having no proper definable subgroup of finite index) as
it is divisible. The
general theory (\cite{Pillay}) of equipping definable groups in $o$-minimal structures with a definable group manifold structure applies to $J$. We 
conclude using (i) of Proposition 8.6 that $J$ with its definable manifold topology is locally Euclidean, and thus
(by the special case of Hilbert's 5th problem due to Pontryagin) is a connected commutative Lie group, whose Lie group 
dimension coincides with its semi-$o$-minimal dimension. By Corollary 8.5 the semi-$o$-minimal dimension of 
$J$ is at most $n$. On the other hand a connected commutative Lie group is a finite product of copies of
$(\R,+)$ and $S_{1}$. As $J$ has the same torsion as $n$-copies of $S_{1}$ it follows that, as a Lie group, $J$ has dimension $\geq n$. It follows that
\begin{Corollary} $J$ is a semi-$o$-minimal connected compact Lie group definable in 
${\bar M}^{*}$ with both semi-$o$-minimal and Lie group dimension $n$.
\end{Corollary}

We will point out later that, as expected, the topology on $J$ coincides with the logic topology on $G/G^{00}$.
But we will be able to complete the proof of Theorem 8.1 without using this. The next step is:
\begin{Lemma} Let $Y$ be a definable (in ${\bar M}$) subset of $G$ of dimension $< n$.
Then 
the semi-$o$-minimal dimension of $\pi(Y)\subseteq J$ is $< n$.
\end{Lemma}
\pf  The proof uses the ideas from the proof of Lemma 10.3 of \cite{NIP}. We give 
some details anyway. We first note that, given $Y$ as in the hypothesis there is a definable subset $Z$ of
$I_{1}\times...\times I_{n}$ such that the map sending $(x_{1},..,x_{n})\to x_{1}+ .. + x_{n}$ induces
a (definable) bijection between $Z$ and $Y$. So $dim(Z) < n$.  Let $E$ be the equivalence relation
on $G$ ``$x-y\in G^{00}$". By the proof of 8.6 we can and will assume that on each of $I_{1},..,I_{n}$ all 
$E$-classes are convex (with respect to canonical total orderings of $I_{1},..,I_{j}$) and that for each 
$i=1,..,n$ $\pi(I_{i})$ is strongly $o$-minimal (in ${\bar M}^{*}$). Let us denote by $\pi_{n}$ the natural map from $I_{1}\times..\times I_{n}$ to $\pi(I_{1})\times..\times \pi(I_{n})$.
\newline
{\em Claim.} $\pi_{n}(Z)$ has dimension $< n$ in the semi-$o$-minimal structure 
$\pi(I_{1})\times ..\times\pi(I_{n})$ (in ${\bar M}^{*}$).
\newline
{\em Proof.} In fact we will prove, by induction on $n$ that if $W$ is any definable 
(in $\bar M$) subset of $I_{1}\times .. \times I_{n}$ of dimension $<n$, then 
$\pi_{n}(W)\subset \pi(I_{1})\times..\times \pi(I_{n})$ has dimension $< n$. And this is an imitation of
the proof of 10.3 in \cite{NIP}: We may assume $W$ is the graph of a continuous definable function 
$f:C\to I_{n}$, where $C$ is an open definable subset of $I_{1}\times..\times I_{n-1}$. Assume for a 
contradiction
that $\pi_{n}(W)$ has dimension $n$, in $\pi(I_{1})\times..\times \pi(I_{n})$. It follows that $\pi_{n}(W)$ contains
the closure of a subset of the form $U \times (a,b)$, where $U$ is an open rectangular box in 
$\pi(I_{1}) \times .. \times \pi(I_{n-1})$ and $(a,b)$ is an interval in $\pi(I_{n})$. 
We may assume that $\pi_{n-1}^{-1}(U)\subseteq C$. Let $c\in I_{n}$ be such
that $a< b < c$, and let $c'\in \pi^{-1}(c)\subseteq I_{n}$. 
For any $x\in U$, there are $x_{1},x_{2}\in \pi_{n-1}(x)$ such that $\pi(f(x_{1})) = a$ and 
$\pi(f(x_{2})) = b$. Then $f(x_{1}) < c' < f(x_{2})$, so by ``definable connectness" of 
$\pi_{n-1}^{-1}(x)$ (and continuity of $f$) there is $x_{3}\in \pi_{n-1}^{-1}(x)$ such that 
$f(x_{3}) = c'$. We have shown that $\pi_{n}(\{z\in C:f(z) = c'\})$ contains $U$, so by induction hypothesis
$\{z\in C:f(z) = c'\}$ has dimension $n-1$. But this is true for infinitely many $c'\in I_{n}$, a 
contradiction. The claim is proved.

\vspace{2mm}
\noindent
We now return to our definable subset $Y$ of $G$ such that 
$Y = \{x_{1} + .. + x_{n}:(x_{1},..,x_{n})\in Z\}$. As $\pi$ is a homomorphism,
$\pi(Y) \subset G/G^{00}  = \{\pi(x_{1}) + ... + \pi(x_{n}):(x_{1},..,x_{n})\in Z\}$
= $\{y_{1} + ... + y_{n}: (y_{1},..,y_{n}) \in \pi_{n}(Z)\}$. By the Claim
$\pi_{n}(Z)$ has dimension $< n$, hence $\pi(Y)$ does too. This completes the proof of Lemma 8.8.

\vspace{5mm}
\noindent
Note that the proof of the above lemma (in fact of  the Claim) actually yields that $dim(\pi(Y)) \leq dim(Y)$ for $Y$ a
definable subset of $G$.

\vspace{5mm}
\noindent
The final point is (the easy):
\begin{Lemma} Let $Z\subset J$ be definable in ${\bar M}^{*}$ of (semi-$o$-minimal) dimension $< n$, and Haar measurable in $G/G^{00}$. Then $Z$ has Haar measure $0$.
\end{Lemma}
\pf  The proof will be by induction on $dim(Z) = k$. For $k = 0$, $Z$ is finite so it is clear.
Now let $Z$ of dimension $k>0$. 
\newline
{\em Claim.} For any $r$ there are $a_{1},..,a_{r}\in J$ such that
$dim(a_{i}Z \cap a_{j}Z) < k$ for all $i\neq j \leq r$. 
\newline
{\em Proof of claim.} Work in a saturated elementary extension of $N$ of ${\bar M}^{*}$. Let $a_{1},..,a_{r}
\in J(N)$ be generic independent (in the $o$-minimal sense) over the base model, namely
$dim(a_{1},..,a_{r}/{\bar M}^{*}) = nr$. Then one sees easily that $dim(a_{i}Z(N)\cap a_{j}Z(N)) < k$ for $i\neq j$. By definability of dimension we can find such $a_{1},..,a_{r}$ in $J$. This proves the claim.

\vspace{2mm}
\noindent
Note that as $Z$ is measurable in $G/G^{00}$, so is any intersection of translates of $Z$.
Hence by the induction hypothesis and the claim, each $a_{i}Z\cap a_{j}Z$ ($i\neq j$) has
Haar measure $0$. But then the measure of the union of the $a_{i}Z$ is $r$ times the measure
of $Z$. So (by choosing $r$ large) this forces $Z$ to have measure $0$. The proof is complete.

\vspace{5mm}
\noindent
{\em Proof of Theorem 8.1.} If $Y\subset G$ is definable in $\bar M$ with dimension $< n$, then
$\pi(Y)$ is closed so measurable in $G/G^{00}$, but by Lemma 8.8, $\pi(Y)$ has dimension $<n$ as a
definable subset of $J$, hence has Haar measure $0$ in $G/G^{00}$ by Lemma 8.9. By 10.5 in \cite{NIP} we 
obtain compact domination of $G$ by $G/G^{00}$ (equipped with its Haar measure). 

\vspace{10mm}
\noindent
We conclude the paper by proving that the topologies on $J$ and $G/G^{00}$ coincide. Our proof will make use of a little more ``theory" some of which is of interest in its own right. 

We begin with an arbitrary complete theory $T$ with $NIP$. For $\phi(x,y)\in L$ by a (complete) global $\phi$-type we mean a maximal consistent collection of the formulas of the form $\phi(x,c)$, 
$\neg\phi(x,c)$ for $c\in \bar M$.

\begin{Lemma} Let $M_{0}$ be a small model, $\phi(x,y)\in L$ and $p_{0}(x)$ a complete global
$\phi$-type which is $Aut({\bar M}/M_{0})$ invariant. Then $p_{0}$ extends to a complete global type $p(x)$ which is $Aut({\bar M}/M_{1})$-invariant for some small model $M_{1}$ (i.e. in earlier terminology $p$ is invariant).
\end{Lemma}
\pf As in the proof of Step I of Lemma 5.8, we can find a small model $M_{1}$ containing $M_{0}$ and a Keisler
measure $\mu_{x}$ over $M_{1}$ which has a unique extension to a Keisler measure over any larger $M_{2}$ which is consistent with $p_{0}$. (To say that $\mu_{x}$ is consistent with $p_{0}$ means that $\mu_{x}\cup p_{0}$ extends to a global Keisler measure. Assuming that $M_{1}$ is $|M_{0}|^{+}$-saturated, which we can and will do, this is equivalent to
requiring that $\mu_{x}$ extends $p_{0}|M_{1}$.) So $\mu_{x}$ has a unique extension to a global
Keisler measure $\mu'$ which extends $p_{0}$. Clearly $\mu'$ is $Aut({\bar M}/M_{1})$-invariant (as $p_{0}$ is). By Proposition 4.6 (see also the proof of 4.5) we obtain some global type $p'(x)$ which is
$M_{1}$-invariant and extends $p_{0}$. 

\begin{Corollary} (Strong Borel definability for invariant $\phi$-types.) Let $p_{0}$ be a complete
global $\phi$-type which is $M_{0}$-invariant. Then $X = \{b:\phi(x,b)\in p_{0}(x)\}$ is a finite
Boolean combination of type-definable over $M_{0}$ sets.
\end{Corollary}
\pf By Lemma 8.10 and Proposition 2.6, $X$ is a finite Boolean combination of type-definable
over $M_{1}$ sets $Y_{i}$ for some small model $M_{1}$ containing $M_{0}$. Let $Y_{i}'$ be
$\{b:\exists c \in Y, tp(c/M_{0}) = tp(b/M_{0})\}$. Then $Y_{i}'$ is type-definable
and $M_{0}$-invariant hence type-definable over $M_{0}$. And $X$ is the same finite
Boolean combination of the $Y_{i}'$.

\vspace{5mm}
\noindent
Here is the application  which will be relevant to our concerns:
\begin{Lemma} 
Let ${\bar M}$ be a saturated model, and $M_{0}$ a small submodel. Let ${\bar M}^{*}$
be the Shelah expansion discussed earlier. Suppose that $X$ is definable in 
${\bar M}^{*}$ and is $Aut({\bar M}/M_{0})$-invariant. Then $X$ is a finite Boolean 
combination of type-definable (over $M_{0}$ in the structure $\bar M$) sets. 
\end{Lemma}
\pf. As $Th({\bar M}^{*})$ has quantifier elimination, there is some complete type
$p(x)$ over $\bar M$ and $L$-formula $\phi(x,y)$ such that $b\in X$ iff $\phi(x,b)\in p$
(for all $b\in {\bar M}$). Let $p_{0} = p|\phi$. So $p_{0}$ is a complete global $M_{0}$-invariant $\phi$
-type and we can apply Corollary 8.11

\vspace{5mm}
\noindent
We now return to the setting and notation of Theorem 8.1 and its proof. In particular $J$ is the set $G/G^{00}$ viewed as a definable (compact Lie) group in ${\bar M}^{*}$ and we just say $G/G^{00}$ for $G/G^{00}$ with the logic topology, another compact Lie group.
\begin{Lemma} (i) Suppose $Z$ is a definable (in ${\bar M}^{*}$) subset of $J$. Then $Z$ is  measurable in $G/G^{00}$. 
\newline
(ii) If moreover $H$ is a definable (in ${\bar M}^{*}$) subgroup of $J$ then $H$ is closed in $G/G^{00}$.
\end{Lemma}
\pf  (i) Clearly $\pi^{-1}(Z)$ is a definable set in ${\bar M}^{*}$ and whether or not $x\in \pi^{-1}(Z)$ 
depends on $tp(x/M_{0})$ in ${\bar M}$. Hence it satisfies the assumptions of Lemma 8.12. Hence by Lemma 8.12, $\pi^{-1}(Z)$ is a finite Boolean combination of 
type-definable sets. We may assume (by multiplying by $G^{00}$) that each of these type-definable sets is a 
union of translates of $G^{00}$. It follows that $Z = \pi(\pi^{-1}(Z))$ is a finite Boolean combination of 
closed subsets of $G/G^{00}$, so Borel and measurable. 
\newline
(ii)
As $G/G^{00}$ is separable, any finite Boolean combination of closed sets is a countable intersection of opens, namely a $G_{\delta}$-set. 
So applying (i) to $H$ we see that $H$ is a $G_{\delta}$ set in $G/G^{00}$. Now the closure 
$\bar H$ of $H$ is a subgroup of $G/G^{00}$. Moreover $H$ and thus each of its translates in $\bar H$ is 
dense in $\bar H$. But any two dense $G_{\delta}$'s must intersect. Hence $H = \bar H$.

\begin{Remark} Again the above lemma holds at various levels of generality: (i) holds assuming just $T$ has $NIP$ and (ii) holds if in addition $G/G^{00}$ is separable. Note that (ii) is saying that any subgroup of $G$ which contains $G^{00}$ and is definable in ${\bar M}^{*}$ is type-definable in 
${\bar M}$. In particular the groups $Stab_{\psi}(p)$ from the proof of Lemma 8.2 are type-definable, so 
using the DCC we obtain another proof that $G^{00}$ is a finite intersection of the $Stab_{\psi}(p)$'s. 
\end{Remark}

\begin{Proposition} The topologies on $J$ and $G/G^{00}$ coincide.
\end{Proposition}
\pf As the group structures on $J$ and $G/G^{00}$ coincide, and both $J$ and $G/G^{00}$ are  compact
(Hausdorff) groups, it suffices to show that any open neighbourhood $U$ of the identity $e$, in the
sense of $J$, contains a neighbourhood of $e$ in the sense of $G/G^{00}$. Let $h$ denote the Haar measure on $G/G^{00}$. Let $W$ be a {\em definable} neighbourhood of $e$ in $J$,  such that
$W^{-1}\cdot W\cdot W^{-1}\cdot W$ is contained in $U$. $W$ is clearly generic in $J$ (finitely many translates
of $W$ cover $J$) by compactness of $J$ for example. By Lemma 8.13, $W$ is measurable in $G/G^{00}$, so 
$h(W) > 0$. As pointed out in the proof of the CLAIM in section 6 of \cite{NIP} it follows that
$W^{-1}\cdot W$ has interior in $G/G^{00}$. (A direct proof: As a measurable set is the union of a
closed set and a measure $0$ set, we may assume that $W$ is closed, hence so is $Z = W^{-1}\cdot W$.
If by way of contradiction $Z$ has no interior then the same holds of any finite union of
translates of $Z$. So we find $a_{1}, a_{2},...$ in $G/G^{00}$ such that the $a_{i}W$ are disjoint, contradicting $h(W) >0$.) So $W^{-1}\cdot W\cdot W^{-1}\cdot W$ contains a neighbourhood of $e$ (in $G/G^{00}$)
as required.

\vspace{5mm}
\noindent
It would be interesting to give a more explicit proof of Proposition 8.15. For example if $dim(G) = 1$ then 
the analysis in the proof of Proposition 3.5 of \cite{Pillay-compact Lie}, yields directly (in hindsight) that 
the logic topology equals the $o$-minimal topology on $G/G^{00}$ and coincides with  $S^{1}$.

\vspace{2mm}
\noindent
Note that once we know Proposition 8.15, then for any definable in ${\bar M}^{*}$ subset $Z$ of $G/G^{00}$,
$Z$ has Haar measure $0$ iff it has dimension $< dim(G/G^{00})$. In any case the proof of Theorem 8.1 
together with Proposition 8.15 says that $G$ is $o$-minimally compactly dominated, i.e.
\newline 
(*) Let ${\bar M}'$ be the expansion of ${\bar M}$ obtained by adding a predicate for $G^{00}$. Then 
$G/G^{00}$ is a semi-$o$-minimal compact Lie group in ${\bar M}'$ with topology coinciding with its topology 
as a bounded hyperdefinable group in ${\bar M}$, AND for any definable, in ${\bar M}$ subset $X$ of $G$, the 
set of $b\in G/G^{00}$ such that $\pi^{-1}(b)$ meets both $X$ and its complement (which is of course
definable in ${\bar M}'$) has dimension $< dim(G/G^{00})$, EQUIVALENTLY, has Haar measure $0$.

\vspace{5mm}
\noindent
Let us first remark that (*) also holds for $X$ definable in ${\bar M}^{*}$. This is because $X = Y\cap G({\bar M})$
for some subset $Y$ of $G(N)$ definable (with parameters) in a saturated elementary extension $N$ of ${\bar M}$.
As $(G/G^{00})(N) = (G/G^{00})({\bar M})$ and (*) is valid working in $N$, we deduce that
the set of $b\in G/G^{00}$ such that $\pi^{-1}(b)$ meets both $X$ and its complement (in $G(M)$), which
again note is a subset of $G/G^{00}$ definable in ${\bar M}^{*}$, has Haar measure $0$. We conclude:
\begin{Corollary} Let $G$ be definably compact and definable in the saturated $o$-minimal expansion ${\bar M}$ of
a real closed field. Let ${\bar M}^{*}$ be the Shelah expansion of $\bar M$. Let $J = G/G^{00}$ as a definable group in ${\bar M}^{*}$. Then $G$ is dominated by $J$ in the weakly $o$-minimal theory
$Th({\bar M}^{*})$: namely, working even in a saturated elementary extension of ${\bar M}^{*}$, if $X$ is a definable subset of $G$ then the set of $b\in J$ such that $\pi^{-1}(b)$ meets both $X$ and its complement, has ($o$-minimal)
dimension $< dim(J)$. 
\end{Corollary}

\noindent
One can ask what the formal content and implications of the compact domination statement
(*) are (either in  general, or restricted to the $o$-minimal context). For example, from \cite{NIP} we know that if a definable group is compactly dominated then it has the $fsg$ property and a unique invariant Keisler measure. Of course the proof of Theorem 8.1 (or statement (*)) depends on $G$ having the $fsg$ property as well as the knowledge of torsion points from \cite{Edmundo-Otero}(for definably compact $G$). It would be interesting to try to recover the torsion points statement directly from compact domination. Namely
\newline
{\em Question.} ($o$-minimal context.) Suppose the commutative definably connected definable group $G$ is $o$-minimally compactly
dominated, i.e. statement (*) holds. Can one prove directly that
$dim(G) = dim(G/G^{00})$ (and so conclude using the torsion-freeness and divisibility of $G^{00}$ that for each $p$ the $p$-torsion of $G$ is $(\Z/p\Z)^{dim(G)}$)?

\vspace{5mm}
\noindent
A final question is whether given a definable group G in a $NIP$ theory $T$, there are 
some reasonable assumptions which imply that $G/G^{00}$ is semi-$o$-minimal (in the Shelah expansion) or at 
least o-minimally analysable. A possible assumption would be that $G$ has the DCC on type-definable subgroups of bounded index (equivalently, $G/G^{00}$ with its logic topology is a compact Lie group).

\appendix

\section{Appendix: On Skolem functions for $o$-minimal definable sets}
We prove a general result, Proposition A.2 below, and show that it applies to the situation in section 8 to yield Proposition 8.6(iii).

Let us fix a saturated structure $N$, and a $\emptyset$-definable set $X$ in $N$ such that for some $\emptyset$-definable dense linear ordering without endpoints $<$ on $X$, $X$ is $o$-minimal in $N$ with respect to $<$. 
We will freely adapt $o$-minimality results for the absolute case (where $X$ is the universe of $N$) to this relative case. 

By $X^{eq}$ we mean $N^{eq}\cap dcl(X)$. It is known that the $o$-minimal dimension theory on $X$ extends smoothly to $X^{eq}$ (see for example section 3.1 of \cite{PPSII}). Namely for any set $A$ of parameters from $N$, and $c\in X^{eq}$ we have a natural number $dim(c/A)$, such that 
$dim$ is subadditive  ($dim((c,d)/A) = dim(c/dA) + dim(d/A)$),
$dim(c/A) = 0$ iff $c\in acl(A)$, for $c$ an element (rather than tuple) of $X$, $dim(c/A) = 1$ if $c\notin acl(A)$.
\begin{Definition} We say that $X$ is {\em untrivial} if whenever $(a_{i}:i<\omega)$ is a sequence of elements of $X$ such that
$a_{i}\notin acl(a_{0},..,a_{i-1})$ for all $i$ then there is some $n$ and $b\in X^{eq}$ such that 
$b\in acl(a_{0},..,a_{n})\setminus acl(a_{0},..,a_{n-1})$ and $a_{n}\notin acl(b)$.
\end{Definition}

The main result here is:
\begin{Proposition} Suppose $X$ is untrivial. Then for any elementary substructure $N_{0}$ of $N$, $X$ has definable Skolem functions in $N$ over $N_{0}$. 
\end{Proposition}

\noindent
Before entering the proof of Proposition A.2 we give a lemma which makes use of the fundamental work of Peterzil and Starchenko \cite{Peterzil-Starchenko}

\begin{Lemma}
Suppose $X$ is untrivial. Then for any $a\in X$ with $a\notin acl(\emptyset)$ there is some interval $(c,d)$ containg $a$ such that the definable set $(c,d)$ has definable Skolem functions in $N$ after naming some parameters from $X$.
\end{Lemma}
\noindent
{\em Proof of Lemma.}  Let $a_{i}$ for $i<\omega$ be realizations of $tp(a)$ such that $a_{i+1}\notin acl(a_{0},..,a_{i})$ for all $i$. Let $n$ and $b$ be given by untriviality of $X$. Then a dimension (or independence) calculation shows that $a_{n}\in acl(a_{0},..,a_{n-1},b)$. We also have $a_{n}\notin acl(a_{0},...,a_{n-1})$ and $a_{n}\notin acl(b)$. Let ${\bar b}$ be a finite tuple of elements of $X$ such that $b\in dcl({\bar b})$ and
${\bar b}$ is independent from $a_{n}$ over $b$  (i.e.  $dim({\bar b}/ba_{n}) = dim({\bar b}/b)$).
So we have that $a_{n}\in acl(a_{0},..,a_{n-1},{\bar b})$ but $a_{n}\notin acl(a_{0},..,a_{n-1})$ and $a_{n}\notin acl({\bar b})$. We may assume $a = a_{n}$. By Claim 1.25 of \cite{PPSI}, $a$ is ``PS-nontrivial", namely for some infinite open interval $I$ containing $a_{n}$ and a definable continuous function $F:I\times I \to X$, strictly monotone in each argument.
By Theorem 1.1 of \cite{Peterzil-Starchenko} there is a convex type definable ordered divisible abelian group $H\subset X$, containing $a$ (where the ordering is the restriction of $<$ to $H$. Let $c,d\in H$ such that $c<a<d$, and let $A$ be some set of parameters from $X$ containing $a,c,d$ and over which the group operation on $H$ is type-definable.
Then it is clear that $(c,d)$ has definable Skolem functions in $N$ over $A$.

\vspace{5mm}
\noindent
{\em Proof of Proposition A.2.}  We begin with some reductions. Let $N_{0}$ be a (small) definably closed substructure of $N$ (or even $N^{eq}$). We will say that $N_{0}$ satisfies Tarski-Vaught with respect to $X$ if any formula over $N_{0}$ which is satisfied in $N$ by some tuple from $X$ is satisfied by a tuple from $N_{0}\cap X$. It is clearly enough to restrict ourselves to formulas $\phi(x)$ over $N_{0}$ where $x$ is a single variable ranging over $X$. In any case clearly any elementary substructure of $N$ satisfies 
Tarski-Vaught with respect to $X$. So by using compactness, in order to prove A.2 it is enough to prove: 
\newline
(*) whenever $N_{0}$ satisfies Tarski-Vaught with respect to $X$, then for any tuple ${\bar c}$ from $N$, $dcl(N_{0},{\bar c})$ satisfies Tarski-Vaught with respect to $X$.

\vspace{2mm}
\noindent
If in the context of (*) $\phi(x)$ is a formula over $N_{0}{\bar c}$ with $x$ ranging over elements of $X$, then $\phi(x)$ defines a finite union $X_{0}$ of intervals and points from $X$. The boundary points of $X_{0}$ are in $dcl(M_{0},{\bar c})$, and $X_{0}$ is defined over the set of these boundary points. So in order to prove (*) it suffices to prove:

\vspace{2mm}
\noindent
(**) whenever $N_{0}$ satisfies Tarski-Vaught with respect to $X$, and ${\bar c}$ is a finite tuple of elements of $X$, then $dcl(N_{0},{\bar c})$ satisfies Tarski-Vaught with respect to $X$.

\vspace{2mm}
\noindent
We can of course prove (**) by adding one element from ${\bar c}$ at a time. Hence it suffices to prove:
\newline
(***) whenever $N_{0}$ satisfies Tarski-Vaught in $N$ with respect to $X$, and $a$ is an element of $X$, then $dcl(N_{0},a)$ satisfies Tarski-Vaught in $N$ with respect to $X$.

\vspace{2mm}
\noindent
The rest of the proof is devoted towards proving (***). 

We will suppose that (***) fails and aim for a contradiction. Using $o$-minimality of $X$, the failure of (*) is equivalent to the existence of an element $a\in X$, and $b\in dcl(N_{0},a)$ (where $b$ might be 
$+\infty$ or $-\infty$)  such that the formula $a < x < b$ (if $b>a$) or $b<x<a$ (if $b < a$) isolates a complete type
over $N_{0}a$. 
There is no harm in assuming that $b>a$. We can write $b = g(a)$ for some $N_{0}$-definable (possibly constant) function $g$ on $P$. There are two cases depending whether or not $g(a)\in dcl(N_{0})$. 
\newline
{\em Case (i).}  $g(a)\notin dcl(N_{0})$. 
\newline
Let $p = tp(a/N_{0})$ and $P$ the set of realizations of $p$ in $N$. Note that $g(a)$ realizes $p$ too, and that  $g$ is an $N_{0}$-definable strictly monotone increasing function from $P$ onto itself. 
\newline
{\em Claim I.}  There is no ${\bar e}$ from $X$ and $N_{0}{\bar e}$-definable function $f_{\bar e}$ such that for all $a'\in P$, $a' < f_{\bar e}(a') < g(a')$.
\newline
{\em Proof.} Otherwise by compactness there is $\theta(x)\in p$, such that
\newline
$N\models \forall x(\theta(x) \rightarrow (x < f_{\bar e}(x) < g(x)))$.  As $N_{0}$ satisfies Tarski-Vaught in $N$ with respect to $X$, there is ${\bar e}'$ from $N_{0}$, such that
\newline
$N\models \forall x(\theta(x) \rightarrow x < f_{{\bar e}'}(x) < g(x))$.
\newline
But then $a< f_{{\bar e}'}(a) < g(a)$, contradicting the fact that $a<x<g(a)$ isolated a complete type over $N_{0}$.
Claim I is proved.

\vspace{2mm}
\noindent
{\em Claim II.}  The interval $(a,g(a))$ has definable Skolem functions in $N$ over some set of parameters from $X$. 
\newline
{\em Proof.}  As $p$ is a complete nonalgebraic $1$-type of $X$ over $N_{0}$, it follows from Lemma A.3 that for some $c$ with $a<c<g(a)$, $(a,c)$ has definable Skolem functions in $N$ after naming some parameters from $X$. It follows that for any $a'$ realizing $p$ and $c'\in (a',g(a'))$, $(a',c')$ has definable Skolem functions in $N$ over parameters from $X$. Now (with $c\in (a,g(a))$), $a < c < g(a) < g(c)$. So both $(a,c)$ and $(c,g(a))$ have definable Skolem functions (after naming parameters). Hence so does $(a,g(a))$.

\vspace{5mm}
\noindent
{\em Claim III.} There are realizations $a_{i}$ of $p$ for $i < \omega$ such that, writing $I_{k}$ for the interval $(a_{k},g(a_{k}))$, we have that for all $k$, $dcl_{N_{0}}(I_{0}\cup...\cup I_{k-1}\cup\{a_{0},a_{1},..,a_{k-1},a_{k}\})\cap I_{k} = \emptyset$. (Where $dcl_{N_{0}}(A)$ denotes definable closure in $N$ of $A\cup N_{0}$.)
\newline
{\em Proof.}  Suppose we have already constructed $a_{0},..,a_{k-1}$, and suppose for a contradiction that 
\newline
($\sharp$):  for all
$a\in P$, $dcl_{N_{0}}(I_{0}\cup..\cup I_{k-1} \cup \{a_{0},..,a_{k-1},a\})\cap (a,g(a)) \neq \emptyset$. 
\newline
By Claim II, let ${\bar e}$ be a tuple from $X$ such that each of $I_{0},..,I_{k-1}$ (and so their union) has definable Skolem functions in $N$ over ${\bar e}$, and we may assume that $a_{0},..,a_{k-1}$ are in ${\bar e}$.  Note that $\cup_{j=0,..,k-1}I_{j}$ has definable Skolem functions over $N_{0}{\bar e}$. Now fix $a\in P$. By 
($\sharp$) there is an $N_{0}$-definable function $f(w,z_{0},..,z_{k-1},x)$ such that there are tuples $c_{0},..,c_{k-1}$ from $I_{0}$,..,$I_{k-1}$ respectively such that $f({\bar e},c_{0},..,c_{k-1},a)\in (a,g(a))$. Hence there are such $c_{0},..,c_{k-1}$ which are in addition contained in $dcl_{N_{0}}({\bar e},a)$. So we have shown that for every $a\in P$, $dcl_{N_{0}}({\bar e},a)\cap (a,g(a)) \neq \emptyset$. By compactness there is a $N_{0}{\bar e}$-definable function $f_{\bar e}(-)$ such that for all $a\in P$, $f_{\bar e}(a) \in (a,g(a))$. This contradicts Claim I. Claim III is proved.

\vspace{5mm}
\noindent
It is rather easy to see that in Claim III, $a_{k}\notin acl_{N_{0}}(a_{0},..,a_{k-1})$. (Alternatively the construction in Claim  III goes through with this additional constraint). So we can apply the untriviality of $X$ to find $n$ and $b\in X^{eq}$ such that 
$b\in acl_{N_{0}}(a_{0},...,a_{n})\setminus acl_{N_{0}}(a_{0},..,a_{n-1})$ 
(whence $a_{n}\in acl_{N_{0}}(a_{0},..,a_{n-1},b)$), and $a_{n}\notin acl_{N_{0}}(b)$. 
This leads quickly to a contradiction as we now show.
At this point we will for notational simplicity work over $N_{0}$. 

First choose $c_{i},d_{i}$ for $i< n$ such that $g^{-1}(a_{i}) < c_{i} < a_{i} < d_{i} < g(a_{i})$, and 
$(c_{0},d_{0},....,c_{n-1},d_{n-1})$ is independent from $(a_{n},b)$ (in the $o$-minimal sense). 
Then as $a_{n}\notin acl(b)$, also $a_{n}\notin acl(c_{0},d_{0},...,c_{n-1},d_{n-1},b)$. On the other hand, as $a_{n}\in dcl(a_{0},..,a_{n-1},b)$ there is a $\emptyset$-definable function $f$ such that
$\models \exists x_{0},..,x_{n-1}((\wedge_{i<n} c_{i}< x_{i} < d_{i})\wedge f(b,x_{0},..,x_{n-1}) = a_{n})$.
As $tp(a_{n}/c_{0},d_{0},..,c_{n-1},d_{n-1},b)$ is not algebraic its set of realizations contains an open interval around $a_{n}$. Hence we can find some $b_{n}\in I_{n}$, and $b_{0},..,b_{n-1}\in X$ such that
for each $i< n$, $c_{i} < b_{i} < d_{i}$ and $f(b_{0},..,b_{n-1},b) = b_{n}$.  Now put $b_{i}' = b_{i}$ if $b_{i}\geq a_{i}$, and $b_{i}' = g(b_{i})$ if  $b_{i} < a_{i}$. So $b_{i}'\in I_{i}\cup\{a_{i}\}$ for each $i$. As
$b\in dcl(a_{0},..,a_{n})$ we conclude that $b_{n}\in dcl(b_{0}',..,b_{n-1}',a_{0},..,a_{n-1},a_{n})$. As $b_{n}\in I_{n}$ this contradicts the construction of the $a_{i}$. This contradiction completes the proof under Case (i).

\vspace{5mm}
\noindent
{\em Case (ii).}  $g(a) \in dcl(N_{0})$. 
\newline
So $g(a)$ is either a point of $X$ in $N_{0}$ or $+\infty$. Let $d = g(a)\in N_{0}\cup\{\infty\})$. So clearly $p(x) = tp(a/N_{0})$ is the complete type over $N_{0}$ saying that $x\in X$, $x<d$ and $x>c$ for all $c\in X(N_{0})$ such that $c<d$. 
\newline
{\em Claim IV.} $P$ (the set of realizations of $p$) is indiscernible over $N_{0}$. Namely for each $n$, 
$p(x_{1})\cup p(x_{2})\cup..\cup p(x_{n}) \cup\{x_{1}<x_{2}<.... < x_{n}\}$ extends to a unique complete type over $N_{0}$.
\newline
{\em Proof.} By induction. The case $n=2$ is given to us, as $a<x<d$ isolates a complete $1$-type over $N_{0}a$.
Now assume true for $n\geq 2$, and prove for $n+1$. 
Let $a_{1}< a_{2} <.. < a_{n}$ realize $p$. It suffices to show that $a_{n} < x < d$ isolates a complete $1$-type
over $N_{0}\cup\{a_{1},..,a_{n}\}$, and for that it is enough to prove that 
$dcl(N_{0},a_{1},..,a_{n})\cap (a_{n},d) = \emptyset$. If not $a_{n}< f(a_{1},..,a_{n})< d$ for some $N_{0}$-definable function $f$. But by induction hypothesis, $tp(a_{n}/N_{0}a_{1}..a_{n-1})$ is isolated by $a_{n-1} < x < d$. Hence $N\models  (\forall x)((a_{n-1} < x < d) \rightarrow (x < f(a_{1},..,a_{n-1},x) < d))$.
Now we use the hypothesis that $N_{0}$ satisfies Tarski-Vaught in $N$ with respect to $X$ to find 
$a_{1}'<...<a_{n-1}'$ in $X(N_{0})$, all less than $d$ such that $x < f(a_{1}',..,a_{n-1}',x) < d$ for any $x$ with $a_{n-1}'< x< d$.  But our realization $a$ of $p$ is such an $x$, and we get a contradiction to our Case(ii) hypothesis.  This finishes the proof of Claim IV.

\vspace{2mm}
\noindent
However Claim IV is clearly incompatible with the untriviality of $X$. So the proof in Case (ii) is also complete, as is the proof of Theorem A.2.

\vspace{10mm}
\noindent
We can now give:
\newline
{\em PROOF OF PROPOSITION 8.6(iii).} We return to the context and notation of 8.6. We already know that $X_{i}$ with its ordering $< _{i}$ is strongly $o$-minimal in ${\bar M}^{*}$, and is also a definable subset of $G/G^{00}$. As in section 8, we use $J$ to denote $G/G^{00}$ as a definable (or intepretable) group in ${\bar M}^{*}$. We will show that $X_{i}$ is untrivial (as a definable strongly $o$-minimal set). 

Let $N_{0}$ denote the structure ${\bar M}^{*}$, and let $N$ be a saturated elementary extension. Then $X_{i}(N)$ is $o$-minimal in $N$ and a definable subset of the definable group $J(N)$. In fact we will work over the parameter set $N_{0}$ over which all the data are anyway defined.

Let $(a_{i}:i<\omega)$ be elements of $X_{i}(N)$ which are algebraically independent over $N_{0}$. For each $i$, let $b_{i} = a_{0}\cdot .. \cdot a_{i}$ where the product is in the sense of the group $J(N)$. So $b_{i}\in J(N)$. On the other hand $b_{i}\in dcl(a_{0},..,a_{i})$ so can be viewed as (is interdefinable with) an element of $X_{i}(N)^{eq}$. So we can talk about $dim(b_{i})$. Let $n$ be as in 8.6, namely the dimension of the original $o$-minimal group $G$. By Corollary 8.5 and Proposition 8.6, each element of $J(N)$ is in the definable closure of $n$ elements of $G$, each of which is a member of some $o$-minimal definable set (defined over $N_{0}$). It follows 
easily that 
$dim(b_{i})$ is bounded by $n$. It is easy to see that $dim(b_{i})$ is nondecreasing (as the $a_{i}$ are independent) and $\leq i+1$. Hence for some $m\leq n$, we have $dim(b_{m-1}) = dim(b_{m}) = m$.
Then $b_{m}\in acl(a_{0},..,a_{m})$, and $b_{m}\notin acl(a_{0},..,a_{m-1})$ (for otherwise we conclude that
$a_{m}\in acl(a_{0},..,a_{m-1})$). Finally $a_{m}\notin acl(b_{m})$. For otherwise, also $b_{m-1}$ is in $acl(b_{m})$ whereby $b_{m}$ is interalgebraic with $(b_{m-1},a_{m})$ hence has dimension $m+1$, which it doesn't.
We have proved that $X_{i}(N)$ is untrivial. So we can apply Proposition 8.2 to obtain 8.6(iii).

\vspace{5mm}
\noindent
\begin{Remark} (i) The argument above yields: Suppose $X$ is a definable strongly $o$-minimal set in a structure $N$ and $X$ definably embeds in a definable group $G$ in $N$ where $G$ has finite thorn rank. Then $X$ is untrivial, so has definable Skolem functions in $N$ after naming parameters.,
\newline
(ii) A recent preprint by Hasson and Onshuus proves that  a strongly $o$-minimal definable set $X$ in a structure $N$ is ``stably embedded" in $N$. So in Proposition A.2 one may assume $X$ to be the universe of $N$, and the set-up of the proof, although not its content, can be a little simplified. 
\end{Remark}


\begin{thebibliography}{99}
\bibitem{Adler} H. Adler, Introduction to theories without the independence property, to appear in Archive Math. Logic.

\bibitem{B-P} Y. Baisalov and B. Poizat, Paires de structures $o$-minimales, Journal of Symbolic Logic 63 (1998), 570-578.

\bibitem{Berarducci} A. Berarducci, $o$-minimal spectra, infinitesimal subgroups and cohomology, Journal of Symbolic Logic, 72 (2007), 1177-1193.

\bibitem{BOPP} A. Berarducci, M. Otero, Y, Peterzil, and A. Pillay, A descending chain condition for groups definable in $o$-minimal structures, Annals of Pure and Applied Logic, 134(2005), 303-313.


\bibitem{Dolich} A. Dolich, Forking and independence in $o$-minimal theories, JSL 69 (2004), 215-240.

\bibitem{Edmundo-Otero} M. Edmundo and M. Otero, Definably compact abelian groups, Journal of Mathematical Logic 4 (2004), 163-180.

\bibitem{Gismatullin} J. Gismatullin, $G$-compactness and groups II, preprint 2007.

\bibitem{HHM1} D. Haskell, E. Hrushovski, and D. Macpherson, Definable sets in  algebraically closed valued 
fields: elimination of imaginaries, J. Reine und Angew. Math. 597(2006), 175-236.

\bibitem{HHM2} D. Haskell, E. Hrushovski, and D. Macpherson, {\em Stable domination and independence
in algebraically closed valued fields}, Lecture Notes in Logic 30, CUP, 2007.

\bibitem{HKP} B. Hart, B. Kim, and A. Pillay, Coordinatization and canonical bases in simple
theories, Journal of Symbolic Logic, 65(2000), 293-309.

\bibitem{metastable} E. Hrushovski, Valued fields, metastable groups, draft 2003.

\bibitem{Hrushovski-Kazhdan} E. Hrushovski and D. Kazhdan, Definable integration in valued fields, in Algebraic Geometry and number theory, Progress Math. 253, Birkhauser, Boston MA, 261-405.

\bibitem{NIP} E. Hrushovski, Y. Peterzil, and A. Pillay, Groups, measures and the $NIP$, Journal AMS, 21 (2008), 563-596. 

\bibitem{centralextensions} E. Hrushovski, Y. Peterzil, and A. Pillay, On central extensions and definably compact groups in $o$-minimal structures, preprint, October 2008.

\bibitem{Ivanov} A. Ivanov, Strongly determined types and G-compactness, Fund. Math. 191(2006), 227-247.

\bibitem{Ivanov-Macpherson} A. Ivanov and D. Macpherson, Strongly determined types, Annals of
Pure and Applied Logic, 99 (1999), 197-230.

\bibitem{Keisler1} H. J. Keisler, Measures and forking, Annals of Pure and Applied Logic 45 (1987), 119-169.

\bibitem{Keisler2} H. J. Keisler, Choosing elements in a saturated model, {\em Classification Theory, Proceedings, Chicago 1985}, ed. J. Baldwin, Lecture Notes in Math. 1292, 1987. 

\bibitem{Lascar-Pillay} D. Lascar and A. Pillay, Hyperimaginaries and automorphism groups, Journal
of Symbolic Logic, 66(2001), 127-143.

\bibitem{NP} L. Newelski and M. Petrykowski, Weak generic types and coverings of groups I, 
Fund. Math. 191(2006), 201-225.

\bibitem{OP} A. Onshuus and A. Pillay, Definable groups and compact $p$-adic Lie groups, Journal LMS, 78(2008). 

\bibitem{O-P} M. Otero and Y. Peterzil, $G$-linear sets and torsion points in definably compact groups, to appear in Archive Math. Logic.

\bibitem{Peterzil-Pillay} Y. Peterzil and A. Pillay, Generics in definably compact groups, Fund.
Math. 193(2007), 153-170.

\bibitem{PPSI} Y. Peterzil, A. Pillay, and S. Starchenko, Definably simple groups in $o$-minimal structures, Transactions AMS, 352(2000), 4397-4419.

\bibitem{PPSII} Y. Peterzil, A. Pillay and S.Starchenko, Simple algebraic and semialgebraic groups over real closed fields, Transactions AMS, 352(2000), 4421-4450.

\bibitem{Peterzil-Starchenko} Y. Peterzil and S. Starchenko, A trichotomy theorem for $o$-minimal structures, Proc. LMS, 77 (1998), 481-523.

\bibitem{Peterzil-Steinhorn} Y. Peterzil and C. Steinhorn, Definable compactness and definable
subgroups of $o$-minimal groups, J. London Math. Soc. 59(1999), 769-786.

\bibitem{Pillay} A. Pillay, On groups and fields definable in $o$-minimal structures, J. Pure and
Applied algebra 53 (1988), 239-255.

\bibitem{Pillay-book} A. Pillay, {\em Geometric Stability Theory}, Oxford University Press 1996.

\bibitem{Pillay-compact Lie} A. Pillay. Type-definability, compact Lie groups, and $o$-minimality,
Journal of Mathematical Logic 4 (2004), 147-162.

\bibitem{Poizat} B. Poizat, {\em A Course in Model Theory}, Springer 2000.


\bibitem{Shelah715} S. Shelah, Classification theory for theories with NIP - a modest beginning, Sci.Japon. 59 (2004), 265-316.

\bibitem{Shelah783} S. Shelah, Dependent first order theories, continued, to appear in Israel J. Math.

\bibitem{Shelah-groups} S. Shelah, Minimal bounded index subgroup for dependent theories, Proc. AMS, 136(2008), 1087-1091.

\bibitem{Shelah-groups} S. Shelah, Definable groups and 2-dependent theories, preprint 2007.

\bibitem{Vapnik-Chervonenkis} V.N. Vapnik and A.Y. Chervonenkis, On the uniform convergence of relative frequencies of events to their probabilities, Theory Probab. Appl., 16 (1971), 264-280.

\end{thebibliography}
\end{document}